\documentclass{aims}
\usepackage{amsmath}
  \usepackage{paralist}
  \usepackage{graphics} 
  \usepackage{epsfig} 
\usepackage{graphicx}  \usepackage{epstopdf}
 \usepackage[colorlinks=true]{hyperref}
\usepackage[final]{changes}
\usepackage{changes}
\usepackage[font=footnotesize]{caption}
\usepackage[utf8]{inputenc}
\usepackage[T1]{fontenc}
\usepackage{mathtools}
\usepackage{amssymb}
\usepackage{amsmath, amsthm, amssymb}
\usepackage{amsthm}
\usepackage{graphicx}
\usepackage{float}
\usepackage{color}

\hypersetup{urlcolor=blue, citecolor=red}

\def\currentvolume{X}
 \def\currentissue{X}
  \def\currentyear{200X}
   \def\currentmonth{XX}
    \def\ppages{X--XX}
     \def\DOI{10.3934/xx.xx.xx.xx}

\theoremstyle{definition}

\newcommand{\ep}{\varepsilon}

\newcommand{\N}{\mathbb{N}}

\newcommand{\R}{\mathbb{R}}

\newcommand{\C}{\mathbb{C}}

\newcommand{\fettk}{\textbf{k}}
\newcommand{\absk}{\abs{\fettk}}

\newcommand{\trk}{\text{tr}  \hat{L} (\absk )}

\newcommand{\detlk}{\text{det}  \hat{L}(\absk )}

\newcommand{\sss}{\sum_{i=1}^3}
\newcommand{\wn}{k}

\newcommand{\rkl}{\right)}
\newcommand{\lkl}{\left(}

\newcommand{\e}{\text{e}}
\providecommand{\abs}[1]{\left\vert#1\right\vert}
\providecommand{\norm}[1]{\left\Vert#1\right\Vert}

\providecommand{\Real}[1]{\text{Re}[#1]}
\providecommand{\ali}[1]{\begin{align}#1\end{align}}
\providecommand{\alinon}[1]{\begin{align*}#1\end{align*}}
\providecommand{\figref}[1]{Fig.\ref{#1}}
\providecommand{\figrefa}[1]{Fig.\ref{#1}(a)}
\providecommand{\figrefb}[1]{Fig.\ref{#1}(b)}
\providecommand{\figrefc}[1]{Fig.\ref{#1}(c)}

\providecommand{\secref}[1]{Section \ref{#1}}

\title[Pattern analysis in a benthic bacteria-nutrient system] 
      {Pattern analysis in a benthic bacteria-nutrient system}

\author[Daniel Wetzel]{}

\subjclass{MBE0623}

 \email{d.wetzel@uni-oldenburg.de}

\begin{document}
\maketitle

\centerline{\scshape Daniel Wetzel$^*$}
\medskip
{\footnotesize
 \centerline{Institut f\"ur Mathematik, Universit\"at Oldenburg, 26111 Oldenburg, Germany}
 \centerline{Fachbereich 3 - Mathematik, Universit\"at Bremen, 28359 Bremen, Germany}
} 

\medskip

\bigskip

 \centerline{(Communicated by the associate editor name)}

\begin{abstract}

We study steady states in a reaction-diffusion system for a benthic bacteria-nutrient model in a marine sediment over 1D and 2D domains by using Landau reductions and numerical path following methods. We point out how the system reacts to changes of the strength of food supply and ingestion. We find that the system has a stable homogeneous steady state for 
\replaced{relatively large}{ great} rates of food supply and ingestion, while this state becomes unstable if one of these rates decreases and Turing patterns \replaced{such as}{ like} hexagons and stripes start to exist. One of the main results of the present work is a global bifurcation diagram for solutions over a bounded 2D domain. This bifurcation diagram includes branches of stripes, hexagons, and mixed modes. Furthermore, we find a number of snaking branches of stationary states, which are spatial connections between homogeneous states and hexagons, homogeneous states and stripes as well as stripes and hexagons in parameter ranges, where both corresponding states are stable.
The \deleted{considered }system \added{under consideration} originally contains some spatially varying coefficients and with these \replaced{exhibits}{ shows some} layerings of patterns.
The existence of spatial connections between different steady states in bistable ranges shows that spatially varying patterns are not necessarily due to spatially varying coefficients.  
 
\deleted{All these results have already been observed for other systems, but most were not proven generally. }\added{The present work} \deleted{This work }gives \replaced{another}{ an other} example, where these effects arise and shows how the analytical and numerical observations can be used to detect signs that a marine bacteria population is in danger to die out or on its way to recovery, respectively.

\deleted{The following results are new. }We find a type of hexagon patches on a homogeneous background\replaced{, which seems to be new discovery.}{ with no prior mention in the literature.} We show the first numerically calculated solution-branch, which connects two different types of hexagons in parameter space. \deleted{We call states on this branch rectangles. }We check numerically \added{for bounded domains} whether the stability changes for hexagons and stripes, which are \replaced{extended}{ continued} homogeneously into the third dimension. We find that stripes and one type of hexagons have the same stable range over bounded 2D and 3D domains. \replaced{This does not hold for the other type of hexagons. Their stable range is shorter for the bounded 3D domain, which we used here. 
We find a snaking branch, which bifurcates when the hexagonal prisms loose their stability. Solutions on this branch connects spatially between hexagonal prisms and a genuine 3D pattern (balls).}{ while the other type of hexagons becomes unstable earlier. Here we find a snaking branch of solutions, which are spatial connects between hexagonal prisms and a genuine 3D pattern (balls).}

\end{abstract}
\section{Introduction}\label{intro}
In this paper we study Turing patterns for a bacteria-nutrient system. It was shown by Turing \cite{turing} in 1952 that nonhomogeneous steady states arise in reaction-diffusion systems, when a homogeneous state is unstable for the full system and stable for the kinetics. This discovery was followed by a large number of works, where systems of different scientific disciplines \replaced{such as}{ like} biology \cite{gierermeinhardt,murray}, chemistry \cite{epstein03,epstein06,epstein09,epstein3d11}, ecology \cite{meron01,meron04,meron05}, and physics \cite{swift77,goswkn84} are studied for so-called Turing patterns. 
\replaced{As}{ It is} discussed in \cite{bac95,bac02,bac03,bac031} \deleted{that }also microorganisms form patterns, \added{and} it is well known that microorganisms play an important role in marine sediments. For example, it is pointed out in \cite{stal03,underwood03} and \cite{gerbersdorf} that extracellular-polymeric-substance secretions of diatoms and benthic bacteria stabilize sediments, respectively. Other experimental investigations of benthic microorganisms can be found in \cite{sed00,sed03,sed03'}.  

\added{In order} to find Turing patterns \added{in standard reaction-diffusion systems}, two different species with different rates of diffusion are required. This makes it difficult to find Turing patterns in chemical experiments, since simple chemicals have almost the same diffusion coefficients, and it took almost 40 years to find the first chemical Turing pattern experimentally \cite{erstesTP}.

During the last 30 years, a great interest arose in localized Turing patterns. 
It was already understood by Pomeau \cite{pomeau} in 1986 that standing fronts which connect a Turing pattern with a homogeneous state and also standing pulses on homogeneous backgrounds which pass near a Turing pattern can be found in reaction-diffusion systems, when both corresponding states are stable. Interestingly\added{,} these states do not only exist\deleted{,} if the conserved quantities of both states are equal, but their branches move back and forth in parameter space and pass stable and unstable ranges. This scenario is referred to as \added{homoclinic} snaking \cite{woods}. There are a lot of works, which investigate this effect over 1D domains (see e.g. \cite{burke,bukno2007,BKLS09}). For a detailed analysis by using the Ginzburg-Landau formalism and beyond all order asymptotics see \cite{chapk09,dean11}.
Fronts and pulses correspond to unbounded domains so that one cannot find these states on bounded domains. What remains are stationary states, which are periodic in space and for which the corresponding orbits pass near the homogeneous state and the Turing pattern. We call such states periodic connections. Their branches also show a snaking behavior (see \cite{bbkm2008,dawes08,dawes09,KAC09,hokno2009} for further details).

The most famous 2D Turing patterns are stripes and hexagons. It \replaced{was}{ is} also understood by Pomeau \cite{pomeau} that standing fronts and pulses should exist in bistable ranges between hexagons and stripes. Periodic connections between hexagons and homogeneous states and between stripes and hexagons are observed in \cite{hilaly} by using numerical time integrations. Investigations of snaking for stationary connections between hexagons and homogeneous states and between hexagons and stripes can be found in \cite{hexsnake,lloyd13} and \cite{schnaki}, respectively.

\subsection{The Model}
The system, which is discussed in the present work, is a reaction-diffusion system for a simplified benthic bacteria-nutrient model in a marine sediment. It \replaced{was}{ is} set up\deleted{ and analyzed} in \cite{ulrike},\deleted{ It is dimensionless and its correctness is not compared to experimental studies, but the authors of \cite{ulrike} claim that this system} models some realistic features, which are motivated by experimental studies, and \added{in dimensionless form} is given by
\begin{equation}\label{dgl}
\begin{aligned}
\partial_t u&=\bigg(\gamma+(1-\gamma)\frac{u}{k+u}\bigg)u\frac{v}{1+v} -m u +\varepsilon +\delta_{u} \Delta u, \\
\partial_t v&=-\bigg(\gamma+(1-\gamma)\frac{u}{k+u}\bigg)u\frac{v}{1+v}+\sigma(v_0-v)+\delta_{v} \Delta v.
\end{aligned}
\end{equation} 
Here $u=u(t,\tilde{x},\tilde{y})$ denotes the population density of one bacteria population and $v=v(t,\tilde{x},\tilde{y})$ the concentration of its (only) nutrient, where $\tilde{x}$ and $\tilde{y}$ are the horizontal and vertical spatial coordinates in the sediment, respectively. $t$ is the time.
$\sigma$, $\gamma, \ k, \ m, \ \varepsilon, \ v_0 \ \delta_u, \ \delta_v$ are parameters, which are all positive.
The terms $\delta_{u} \Delta u$ and $\delta_{v} \Delta v$ are used to describe the diffusion of bacteria and the nutrient, respectively. The bacteria are larger and heavier than its nutrient, so the nutrient diffuses faster than the bacteria. This can be modeled by setting $\delta_u<\delta_v$.
The ratio of active bacteria is described by
\alinon{
\bigg(\gamma+(1-\gamma)\frac{u}{k+u}\bigg).
}
By this approach the following two features are modeled\replaced{:}{.} a certain part $\gamma\in[0,1]$ of bacteria is always active, which means that they search for nutrients and the active bacteria send signal molecules to activate dormant bacteria. This communication works well for a high population density, i.e., the half saturation term $u/(k+u)$ is approximately one so that almost all bacteria are active, while $u/(k+u)$ is smaller than one for a small density of bacteria so that not all dormant bacteria can be activated. The half saturation term $v/(1+v)$ tells us that not all bacteria are able to find nutrients for small nutrient concentrations. That a communication between bacteria and the use of half saturation terms for describing the growth of bacteria is realistic \replaced{was}{ is} already pointed out in \cite{look,look2} and \cite{monod}, respectively. The linear term $mu$ models the mortality of bacteria and $\varepsilon$ is the rate of bacteria inflow. 
The nutrient concentration in the sea water is given by $v_0$. To understand how the term $\sigma(v_0-v)$ enriches the model, we consider the system
\ali{\partial_t v=\sigma(v_0-v).\label{sigi}
} 
The unique solution of \eqref{sigi} is given by
\alinon{v=(v(0)-v_0) e^{-\sigma t}+v_0,}
which converges to $v_0$ for $t\rightarrow \infty $.
Thus the parameter $\sigma$ determines the rate at which the nutrient concentration in the sediment adapts to $v_0$. In \cite{ulrike} it is pointed out that this adaptation of the nutrient concentration comes from a transport of nutrients, which occurs as a result of burrow and pump activities by worms, shells and other animate beings in the sediment. This process is referred to as bioirrigation. In the following we call $\sigma$ the balancing rate.

Since Turing's fundamental paper \cite{turing} it has been known that some reaction-diffusion systems \replaced{possess spatially}{ are solved by} non-homogeneous solutions. The most famous ones in 2D are stripe and hexagonal spot patterns. Sometimes spot patterns are classified \replaced{into}{ in} cold and hot, which means that they have a minimum and maximum in the center of every spot, respectively. Because of the predator prey structure of \eqref{dgl} we have a hot-spot pattern for $u$, when we have a cold-spot pattern for $v$ and vice versa. Thus we always present the pattern of $u$ only, and when we are saying that a solution of \eqref{dgl} is hot or cold, this means that this is the case for $u$. 

\begin{figure}[H]
\begin{center}
\includegraphics[width=0.7\textwidth]{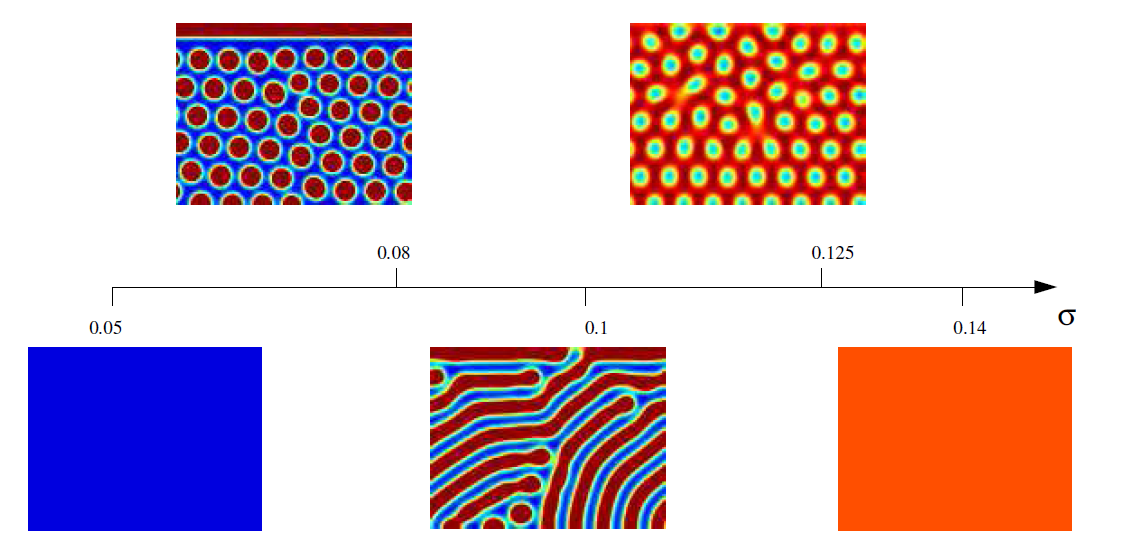}
\end{center}
\caption{This figure is taken from \cite{ulrike}. Shown are quasi-stable solutions of \eqref{dgl} for $\gamma=0.25$, parameter set \eqref{para}, and different balancing rates $\sigma$.}\label{feudelfig}
\end{figure}

Such so-called Turing patterns \replaced{have also been}{ are also} found in \cite{ulrike} for \eqref{dgl}. In \figref{feudelfig} we see the following quasi-stable patterns \footnote{\added{With this we mean solutions that appear to be stationary in time-domain simulations, but may actually change very slowly.}}\deleted{, which are found via time-integrations}: A homogeneous pattern with a low density of bacteria for $\sigma=0.05$, hexagonal hot-spots for $\sigma=0.08$, stripes for $\sigma=0.1$, hexagonal cold-spots for $\sigma=0.125$, and a homogeneous pattern with a high density of bacteria for $\sigma=0.14$.

\deleted{We give cause for concern that the}\added{The} model \eqref{dgl} is set up for 2D domains, where the horizontal and vertical dimension are considered, but the \added{third} spatial dimension is neglected. When we \replaced{extend}{ continue} a stable\deleted{ 2D} solution\added{ on a bounded 2D domain} homogeneously into the third dimension, it is also a solution for \deleted{the system over }\added{this extended}\deleted{ a} 3D domain, but can be unstable for this domain. If this is the case, it is not observable in nature. Clearly, the stability results\added{ for this solution} hold, when one considers a sediment, which is very thin in the third dimension\added{ and bounded for the first and second dimension}, but this does not hold for most marine sediments.\\
However, not all marine bacteria populations live inside the sediment and form 3D colonies. For instance phototrophic bacteria live on thin films on the marine sediment for practicing photosynthesis. For such films one can also use system \eqref{dgl} to model population densities of such bacteria.
The parameter $v_0$ can be seen as the nutrient concentration in the sea water, sediment, or both. 
\replaced{Another}{ An other} way to vindicate the investigation of system \eqref{dgl} over 2D domains is that we
\deleted{actually }treat 3D domains and consider 2D domains as a first step for a better understanding.

\deleted{For reaction-diffusion systems over bounded and unbounded domains it holds that a \added{stationary }solution is spectrally, linearly, or nonlinearly stable in $L^2$ (in the following we simply say stable) if all real parts of \replaced{points in the spectrum}{ the eigenvalues} for the linearization are negative and unstable if at least one real part is positive. If one real part vanishes and the other ones are negative, then we see this solution as a possible bifurcation point and do not care about its stability.}

The goal of this paper is to continue the investigations of \cite{ulrike} by understanding the bifurcation scenarios of \eqref{dgl} and to find and investigate more stationary patterns via bifurcation \replaced{analyses}{ analysis}. We do this over 1D domains\deleted{ as an introduction}, before we start to consider 2D domains. \\
From the discussion above one can see that $\sigma$ and $\gamma$ correlate with the strength of food supply and ingestion for the bacteria-nutrient system, respectively, which clarifies that these parameters play an important role for the system. In this paper we study how the system reacts, when $\sigma$ and $\gamma$ change their values. 
Mostly we treat $\gamma$ as a given and fixed parameter, while we use $\sigma$ as a bifurcation parameter, i.e., we examine how solutions and their types of stability change and what kind of new solutions bifurcate by varying $\sigma$.
For all other parameters we use the main parameter set of \cite{ulrike}, which is given by 
\ali{k=1, \quad v_0=4.125,\quad \varepsilon=0.005,\quad m=0.3175,\quad \delta_u=2 \cdot 10^{-5}, \quad \delta_v=10^{-3}. \label{para}
}

We use the Landau reduction to understand bifurcation scenarios for unbounded domains locally. To get a \added{more} global bifurcation diagram, we use the continuation and bifurcation software {\tt pde2path} \cite{p2p}. It uses numerical methods \replaced{such as}{ like} the finite element method, so we are not able to treat unbounded domains. Instead we consider bounded domains with Neumann boundary conditions, \added{so all found solutions can be \replaced{extended}{ continued} periodically over the unbounded domain.} 

\added{ A stationary solution of a reaction-diffusion system over bounded and unbounded domains is spectrally, linearly, and nonlinearly stable in $L^2$ if all real parts of points in the spectrum for the linearization are negative. We use the discretized Jacobian to find results about the stability for solutions which we find via the finite element method and call such a solution stable if all real parts of the Jacobian are negative. We call them unstable if at least one real part is positive. In the following stability refers to this setting. For more mathematical information about stability for PDEs we refer the reader to \cite{sandstede2002}. We do not consider stability at bifurcation points, where the real part of eigenvalues vanishes.}\deleted{Since we use Neumann boundary conditions, all found solutions can be continued periodically over the unbounded domain.} A solution\deleted{,} which is stable over a bounded domain, is not necessarily stable over larger domains.\deleted{ This does not hold for unstable solutions.} If a solution is unstable over a bounded domain, then it is also the case for larger domains\added{, for which the unstable modes fit into the domain}.   \\
One of the main results of this paper is a global bifurcation diagram (see \figref{hsbbifu}), which can be seen as a continuation of \figref{feudelfig} and from which we can read \replaced{off}{ out} the existence and stability of hexagon, stripe, mixed mode, and homogeneous solutions. Thus we can determine if the bacteria population is in danger of extinction.  \\
Furthermore, we find bistable ranges between two different \replaced{types}{ Types} of solutions. For such a bistable range one can show analytically that a necessary condition for a heteroclinic connection is fulfilled (see \cite{pomeau}). Homo- and heteroclinics are solutions on unbounded domains, so we cannot find those by using the finite element method, but \added{as an approximation} we find \added{spatially }periodic connections \added{with large periods} between both stable solutions. This means that a layering of two patterns is not necessarily an effect of space dependent parameters but can occur for homogeneous balancing rates.

\subsection{Outline}
\deleted{It is shown in \cite{ulrike} that the problem of determining the homogeneous solutions can be reduced to computing the zeros of a cubic polynomial. The coefficients are given with some typos: the first $\gamma$ of $a_1$ must be $\sigma$. The zeros of this polynomial are not derived analytically in \cite{ulrike} because of their complex terms. In \secref{sechom} we show the reduction to the cubic polynomial again and sort the coefficients with respect to $\sigma$ and $\gamma$, because we investigate \eqref{dgl} in the following for different values of $\sigma$ and $\gamma$, while we use the set \eqref{para} for all other parameters of \eqref{dgl}. Furthermore, we apply analytical formulas for the three zeros of the cubic polynomial to determine the regions of a bounded domain of the $\sigma$-$\gamma$-plane, where the homogeneous solutions are real to see for which combinations of $\sigma$ and $\gamma$ exist only one or rather three homogeneous states.
After this follows an introduction of Turing instabilities and bifurcations. By using a fine discretization of the $\sigma$-$\gamma$-domain, we determine Turing unstable ranges. With this information we are able to predict the bifurcation locations and wavenumbers of spatially periodic 1D and 2D solutions.}  
\added{In \secref{sechom} we give analytical formulas for the homogeneous solutions of \eqref{dgl} and study their stabilities.} 
In \secref{seclandau} we recall how to reduce a general two species reaction-diffusion system to the Landau amplitude equation system on a hexagonal lattice. 
In \secref{sec1d} we study non trivial patterns over 1D domains
and show global bifurcation diagrams for $\gamma=0.3$. We find that some stripe solutions bifurcate from the homogeneous branch and terminate in a different bifurcation on the homogeneous branch by holding their wavelength. Moreover, we find stripe branches of the same wavelength which are not connected, but change their wavelength and connect to other stripe branches. Using the Landau reduction, we choose a $\gamma$-value where stripes bifurcate subcritically, to \replaced{generate}{develop} bifurcation diagrams with snaking branches of localized stripes on homogeneous backgrounds.

In \secref{2d} we investigate 2D patterns and \replaced{generate}{ develop} some global bifurcation diagrams. As in \cite{schnaki} we find solutions and solution branches of localized hexagonal spots with a planar interface to striped backgrounds. Furthermore, we find patches of localized hexagons on homogeneous backgrounds. Some of these patches are already shown in \cite{hexsnake}, while others are not mentioned in the literature before. 
\deleted{Some patches which we find have a flat overlying interface such that the snake does not show a strong snaking behavior. We use the Landau reduction to find $\gamma$-values, where the interface is steeper and thus the snaking behavior is more pronounced.}

In \cite{ulrike} a main point is to consider the system for space dependent parameters.  In \secref{seclayering} we discuss what this means and show some layering of patterns for such a system. 
\added{In \secref{discussion} we discus the relevance of the stationary states which we found for the bacteria-nutrient system \eqref{dgl}. }
\deleted{Comparing this to localized patterns of the section before we point out in \secref{discussion} that a layering of patterns is not necessarily an effect of space dependent parameters, but can occur for homogeneous ones.
Furthermore, we discuss the relevance of the different patterns, which we found, for the bacteria-nutrient system. The main point is that we can determine from the pattern if a system is in danger to die out.}

\deleted{We call hexagons and stripes, which are extended homogeneously into the third dimension hexagonal prisms and lamellae, respectively. 
We check the stability of hexagonal prisms and lamellae over bounded 3D domains in \secref{sec3d}.  Furthermore, we show a snaking branch of spatial connections between hexagonal prisms and genuine 3D patterns.}

\section{Homogeneous solutions and Turing instabilities}\label{sechom}
\added{It is shown in \cite{ulrike} that the problem of determining the homogeneous solutions can be reduced to computing the zeros of a cubic polynomial. The coefficients are given with some typos: the first $\gamma$ of $a_1$ must be $\sigma$ in \cite[p.115]{ulrike}. The zeros of this polynomial are not derived analytically in \cite{ulrike} because of their complex terms. In this section we show the reduction to the cubic polynomial again and sort the coefficients with respect to $\sigma$ and $\gamma$, because we investigate \eqref{dgl} in the following for different values of $\sigma$ and $\gamma$, while we use the set \eqref{para} for all other parameters of \eqref{dgl}. Furthermore, we apply analytical formulas for the three zeros of the cubic polynomial to determine the regions of a bounded domain of the $\sigma$-$\gamma$-plane, where the homogeneous solutions are real to see for which combinations of $\sigma$ and $\gamma$ exist only one or rather three homogeneous 	states.}

\begin{figure}[h]
\begin{center}
\boxed{
\begin{minipage}{0.3\textwidth}
(a)\\
\includegraphics[width=1\textwidth]{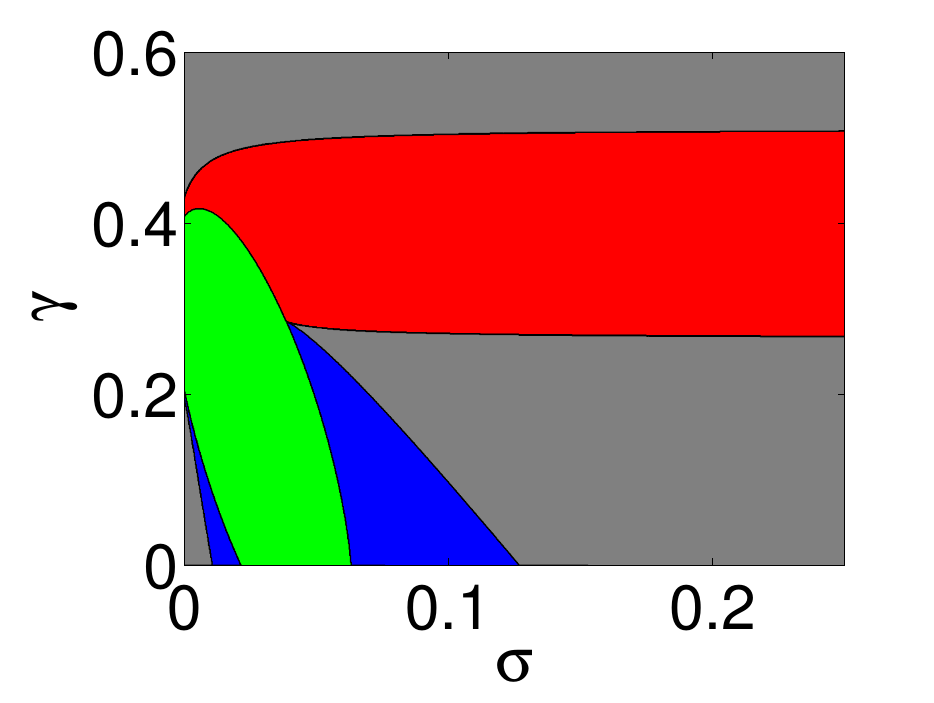}
\end{minipage} }
\boxed{
\begin{minipage}{0.3\textwidth}
(b)\\
\includegraphics[width=1\textwidth]{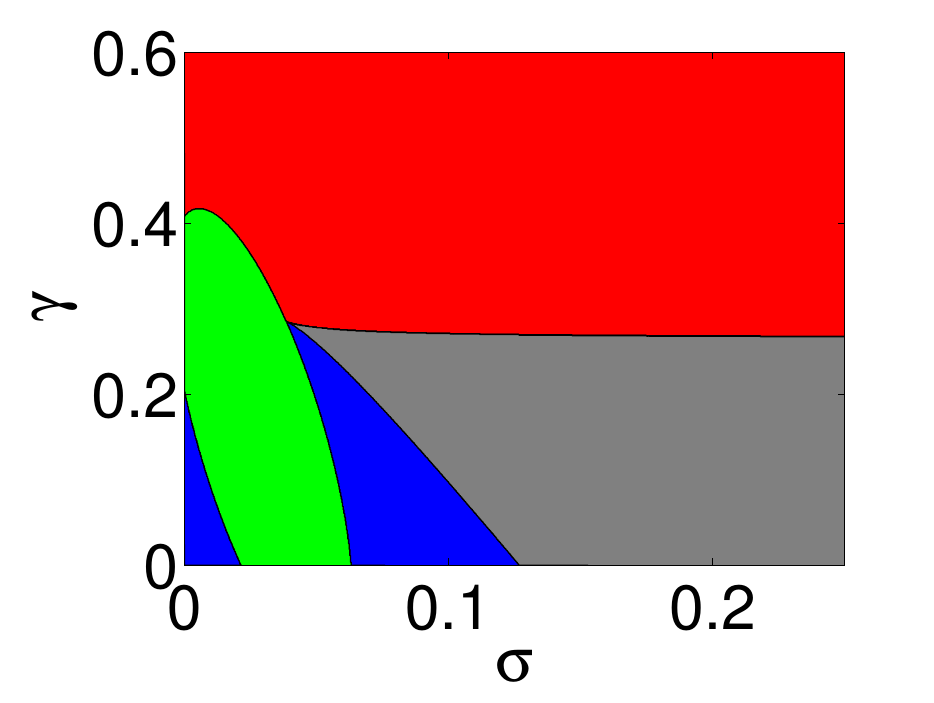}
\end{minipage} }
\end{center}
\caption{(a) $(u_1,v_1)$, $(u_2,v_2)$, and $(u_3,v_3)$ are real in red, green, and blue regions, respectively. All three are real in the gray regions. This also holds for (b), where, in addition, the three homogeneous states are positive in the gray regions. To \replaced{compute}{ perform} these figures, we discretize the $\sigma$-$\gamma$-domain $(0,0.25)\times(0,0.6)$ in $1000\times 1000$ points and checked the corresponding conditions for every point.}
\label{homplot}
\end{figure}

Rescaling \eqref{dgl} with $x=\frac{\tilde{x}}{\sqrt{\delta_u}}$ and $y=\frac{\tilde{y}}{\sqrt{\delta_u}}$, yields
\begin{equation}\label{dgl2}
\begin{aligned}
\partial_t u&=\bigg(\gamma+(1-\gamma)\frac{u}{k+u}\bigg)u\frac{v}{1+v} -m u +\varepsilon +\Delta u,\\
\partial_t v&=-\bigg(\gamma+(1-\gamma)\frac{u}{k+u}\bigg)u\frac{v}{1+v}+\sigma(v_0-v)+\delta \Delta v,
\end{aligned}  
\end{equation}
where $\delta=\frac{\delta_v}{\delta_u}=50$. 
In order to simplify the notation we set $w=(u,v)$, 
$D=\lkl\begin{smallmatrix}1 & 0\\ 0& \delta\end{smallmatrix}\rkl$, and
\ali{ f(u,v)=\left(\begin{array}{c} g(u,v) \\ h(u,v) \end{array}\right)=\left(\begin{array}{l} \bigg(\gamma+(1-\gamma)\frac{u}{k+u}\bigg)u\frac{v}{1+v} -m u +\varepsilon \\ -\bigg(\gamma+(1-\gamma)\frac{u}{k+u}\bigg)u\frac{v}{1+v}+\sigma(v_0-v) \end{array}\right).  \label{kin}
}
The function $f$ is called the reaction term or kinetic of the full system \eqref{dgl2}. 
With this simplifications we can write \eqref{dgl2} as
\ali{\partial_t w= f(w)+D\Delta w \text.\label{rds}}
Homogeneous steady states are solutions of $f(u,v)=0$ which are space and time independent. 
For a nonzero rate of bacteria inflow $\ep$ and nonzero product of the balancing rate $\sigma$ and the nutrient concentration $v_0$ in the sea water it holds \added{that}
\ali{u \neq 0, \quad v \neq 0\label{nonzero}}
for homogeneous solutions of \eqref{dgl2}.
Adding the first to the second equation of \eqref{kin}, yields the linear relationship 
\ali{v=v_0-\frac{mu-\ep}{\sigma}.  \label{linrelation}} 
From \eqref{nonzero}, \eqref{linrelation}, and  the fact that the population density of the bacteria $u$ and the concentration of the nutrient $v$ cannot be smaller than zero follows
\alinon{u\in\lkl 0,\frac{\ep+\sigma v_0}{m}\rkl,\quad v\in \lkl 0, v_0+\frac{\ep}{\sigma}\rkl .}
Substituting \eqref{linrelation} into $f(u,v)=0$, reduces the problem of finding homogeneous solutions of \eqref{dgl2} to the problem of finding the zeros of the polynomial 
\ali{\label{gleichung}
u^3+bu^2+cu+d=0
}
with
\alinon{
b=b_g \gamma +b_s \sigma+b_0,\quad
c=c_g \gamma+c_s \sigma+c_{s g}\sigma \gamma+c_0,\quad
d=d_s \sigma+d_0,
}
where 
\alinon{
&b_g=-\frac{k}{m-1} \approx 1.47, \
b_s=\frac{v_0-v_0m-m}{m(m-1)}\approx -11.53,\
b_0=\frac{m^2k+\ep-2m\ep}{m(m-1)}\approx-0.47,\\
&c_g=\frac{\ep k}{m(m-1)}\approx -0.02, \ 
c_s=\frac{(v_0+1)(\ep-mk)}{m(m-1)}\approx 7.39, \
c_{sg}=\frac{v_0 k}{m(m-1)}\approx -19.04, \\
&c_0=\frac{-2m \ep k+\ep^2}{m(m-1)}\approx 0.01,\
d_s=\frac{\ep k(v_0+1)}{m(m-1)}\approx -0.12,\
d_0=\frac{\ep^2 k}{m(m-1)}\approx -10^{-4}.
}
The zeros of \eqref{gleichung} are given by 
\alinon{u_1&=\sqrt{-\frac{4p}{3}} \ \cos \lkl\frac13 \arccos\lkl-\frac{q}{2} \sqrt{-\frac{27}{p^3}}\rkl \rkl  -\frac{b}{3}, \\
u_2&=-\sqrt{-\frac{4p}{3}} \ \cos \lkl\frac13 \arccos\lkl-\frac{q}{2} \sqrt{-\frac{27}{p^3}}\rkl +\frac{\pi}{3} \rkl  -\frac{b}{3}, \\
u_3&=-\sqrt{-\frac{4p}{3}} \ \cos \lkl\frac13 \arccos\lkl-\frac{q}{2} \sqrt{-\frac{27}{p^3}}\rkl -\frac{\pi}{3} \rkl   -\frac{b}{3},
 }
 where 
\alinon{
 p=c-\frac{b^2}{3}, \quad q=\frac{2b^3}{27}+d-\frac{bc}{3}.
 }
The polynomial \eqref{gleichung} can have one or three real zeros for fixed $\sigma$ and $\gamma$. The regions of the $\sigma$-$\gamma$-plane, where a homogeneous solution is real, can be calculated analytically. Determining analytically, where the population density of the bacteria and the nutrient concentration of such a solution is also positive, seems not so trivial. In \figrefa{homplot} we show the regions on a bounded domain of the $\sigma$-$\gamma$-plane, where homogeneous solutions are real. They are also positive for the regions shown in \figrefb{homplot}. Comparing (a) and (b), we see that the positivity condition is not fulfilled for all homogeneous solutions $(u,v)$ of \eqref{gleichung}. 
Here $(u_2,u_3)$ and $(v_1,v_2)$ are negative in the upper horizontal gray band and the lower left gray region of \figref{homplot}(a), respectively.  

We verified that the conditions 
\alinon{u_1>u_2>u_3 \qquad \text{and} \qquad v_1<v_2<v_3 }    
hold for all discretization points used to generate \figref{homplot} if all three homogeneous states are real. 
Thus $(u_1,v_1)$ is the state with the highest population density of bacteria and lowest concentration of nutrient of these three homogeneous states, while this is \deleted{the }opposite\deleted{ case} for $(u_3,v_3)$. The equilibrium $(u_2,v_2)$ lies in the middle.  

The Jacobian of $f$ is given by 
\alinon{
J_f(u,v)=\lkl \begin{matrix} g_u & g_v \\ h_u & h_v \end{matrix} \rkl = \lkl \begin{matrix} \xi-m & \vartheta \\ -\xi & -\vartheta-\sigma \end{matrix} \rkl,
}
where
\alinon{
\xi=\lkl \gamma +(1-\gamma) \frac{2ku+u^2}{(k+u)^2} \rkl \frac{v}{1+v}, \qquad
\vartheta = \lkl \gamma +(1+\gamma)\frac{u}{k+u} \rkl \frac{u}{(1+v)^2}.
}
The linearization of \eqref{rds} in a homogeneous state $w^*$ is given by
$\partial_t w=L(\Delta)(w-w^*)$, 
where $L(\Delta)=J_f+D\Delta$. 
It holds 
\alinon{L(\Delta) e^{i (x,y)\cdot \fettk}=\hat{L}(\absk)e^{i(x,y)\cdot \textbf{k}} \quad \text{with} \quad 
 \hat{L}(\absk)=J_f-D \absk^2 \quad \text{and} \quad \fettk \in \R^2 .}
 This yields the eigenvalue problem 
 \begin{align}\label{eigprob}\hat{L}(\absk) \phi(\absk)=\mu(\absk) \phi(\absk),\end{align}
where 
\begin{align}
\mu_\pm(\absk)=\frac{\trk}{2} \pm \sqrt{\left(\frac{\trk}{2} \right)^2-\detlk}. \label{mu}
\end{align}
In \secref{intro} we already recalled that a homogeneous solution $w^*$ is stable if $\Real{\mu_\pm(\absk)}$ $<0$ for all $\absk$ and unstable if there is a $\fettk \in \R^2$ so that $\Real{\mu_+(\absk)}>0$ or $\Real{\mu_-(\absk)}>0$.
Furthermore, $w^*$ is called Turing-unstable if $w^*$ is unstable in the full system \eqref{rds}, but stable in $\partial_t w=f(w)$. It can be shown easily that $w^*$ is Turing-unstable if the following two conditions are fulfilled: \\
\begin{tabular}{l l}
i) & $\Real{\mu_+(0)}<0$ and $\Real{\mu_-(0)}<0$. \\ 
 ii) & There is a $\fettk$ such that $\Real{\mu_+(\absk)}>0$ or $\Real{\mu_-(\absk)}>0$. \\ \\
\end{tabular}

\noindent Notice that the conditions i) and ii) are not necessary for Turing instabilities. We call $w^*$ space-independent unstable if $w^*$ is unstable in $\partial_t w=f(w)$. 
Let
\alinon{
b_1=-g_u-h_v,\quad
b_2=g_u h_v- g_v h_u,}
\alinon{
b_3=\delta g_u+h_v, \quad
b_4=(\delta g_u+h_v)^2-4 \delta (g_u h_v - g_v h_u). }
It holds that a homogeneous state of \eqref{rds} is 
\ali{
&\bullet \text{stable if $b_1>0 \text{ and } b_2>0 \text{ and } (b_3<0 \text{ or } b_4<0)$,}\label{turing1}\\
&\bullet \text{space-independent unstable if $b_1<0 \text{ or } b_2<0$,}\label{turing2}\\
&\bullet \text{Turing unstable if $b_1>0 \text{ and } b_2>0 \text{ and } b_3>0 \text{ and } b_4>0$.}\label{turing3}
}
The first and second conditions can be shown easily. The third one is shown in \cite{murray}.
We call a solution a Turing endpoint and Turing bifurcation point if it lies at a transition from Turing-unstable to space-independent unstable and from Turing-unstable to stable, respectively. Turing bifurcation points occur if $\detlk$ vanishes. We call the corresponding $\sigma$ and $\absk$ critical balancing rate $\sigma_c$ and critical wavenumber $k_c$, respectively.

\begin{figure}[h]
\boxed{
\begin{minipage}{0.27\textwidth}
(a)\qquad\quad $(u_1,v_1)$\\
\includegraphics[width=1\textwidth]{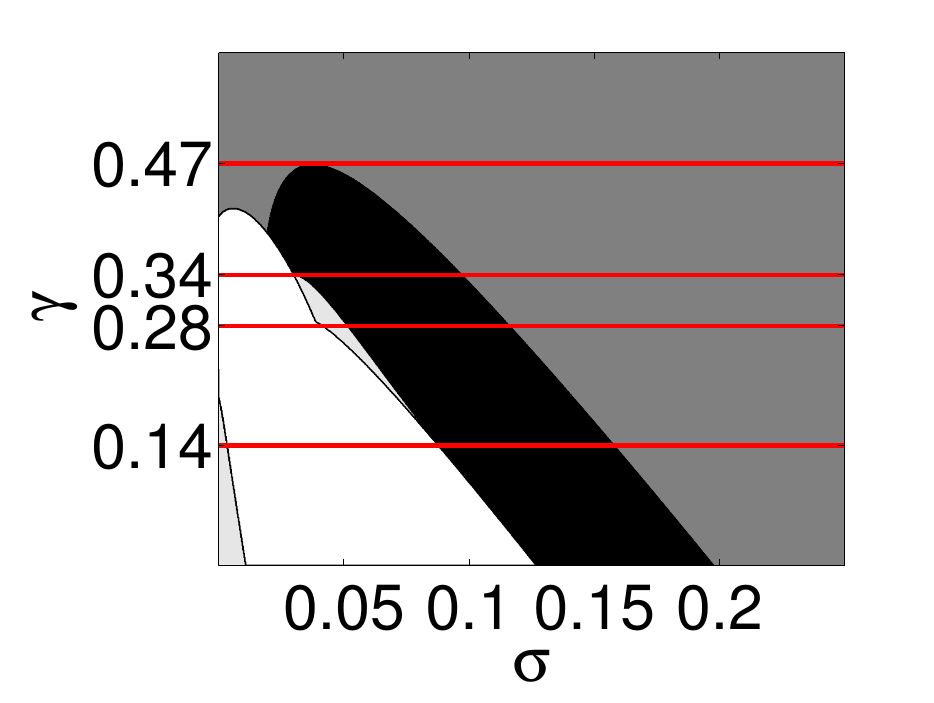}
\end{minipage} }
\boxed{
\begin{minipage}{0.27\textwidth}
(b)\qquad\quad $(u_2,v_2)$\\
\includegraphics[width=1\textwidth]{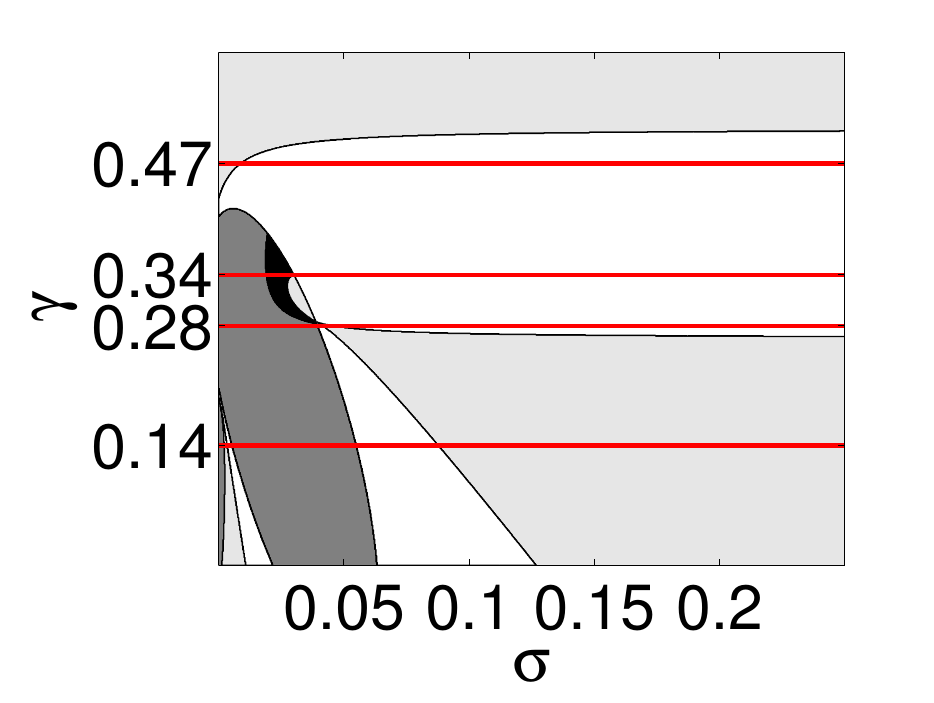}
\end{minipage} }
\boxed{
\begin{minipage}{0.27\textwidth}
(c)\qquad\quad $(u_3,v_3)$\\
\includegraphics[width=1\textwidth]{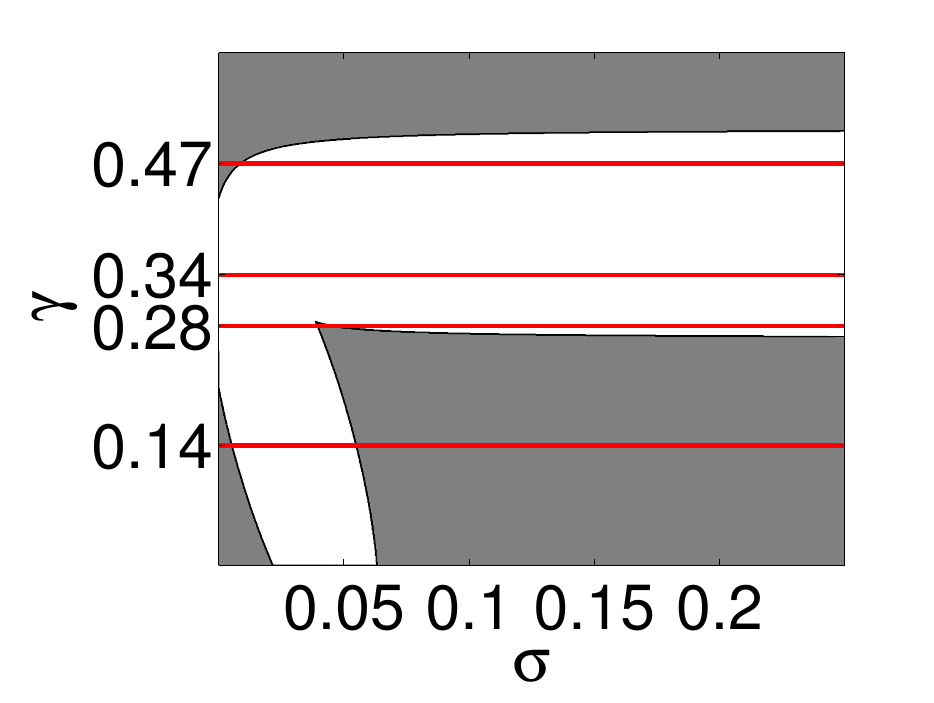}
\end{minipage} }
\boxed{
\begin{minipage}{0.07\textwidth}
(d)\\
\includegraphics[width=1\textwidth]{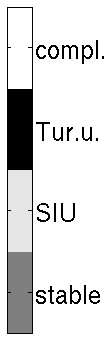}
\end{minipage} }

\caption{Illustrated are the stabilities for $(u_1,v_1)$, $(u_2,v_2)$, and $(u_3,v_3)$ in (a), (b), and (c), respectively. We use the same $\sigma$-$\gamma$-domain and discretization as in \figref{homplot} to check the corresponding conditions \eqref{turing1}, \eqref{turing2}, and \eqref{turing3}. The red lines are the boundaries of the intervals $I_1$, $I_2$, $I_3$, $I_4$, and $I_5$. The color bar for (a), (b), and (c) is shown in (d). The abbreviations compl., Tur.u., and SIU stand for complex, Turing-unstable, and space-independent unstable, respectively.   }

\label{stabplot}
\end{figure}

  In \figref{stabplot} we illustrate the stabilities of the three different homogeneous states in the $\sigma$-$\gamma$-domain of \figref{homplot}. First of all, we can see in \figrefc{stabplot} that $(u_3,v_3)$ is always stable if it is real in our chosen region. 
 By using $\sigma$ as bifurcation parameter and the activity stimulation $\gamma$ as a fixed parameter, we can classify the bifurcation scenarios into five different types. We do this by partitioning the $\gamma$-interval into the following five sections 
\alinon{
I_1=(0,0.14], \ I_2=(0.14,0.28] , \ I_3=(0.28,0.34] , 
\ I_4=(0.34,0.47] ,  \ I_5=(0.47,0.6].
}



The least interesting interval is $I_5$. Here the state $(u_1,v_1)$ is always real, positive, and stable, while the other two are not real or not positive. For all other intervals we have Turing-unstable ranges. 






For  $I_3$ and $I_4$ we always have one homogeneous solution. The Turing-unstable range is continuous for $I_4$, while it is not continuous for $I_3$ so that we have two Turing-unstable ranges, which are separated by a space-independent-unstable range, which is bounded by Turing endpoints.\\
For $I_1$ and $I_2$ we always have ranges, where three homogeneous solutions exist. On the left boundary of these ranges is a fold, where $(u_1,v_1)$ equals $(u_2,v_2)$. The Turing-unstable range for $I_1$ starts at this fold, while the Turing-unstable range begins for $I_2$ in a Turing endpoint on the right side of this fold. 
\\
Turing patterns branch from Turing bifurcation points as discussed below. At Turing endpoints we have $\Real{\mu_\pm(0)}=0$, while $\text{Im}[\mu_\pm(0)]\neq 0$ such that a necessary condition for Hopf bifurcations is fulfilled at Turing endpoints. Currently\added{ the software} {\tt pde2path} \cite{p2p} does not handle Hopf bifurcations and \replaced{we do not consider solution branches, which bifurcate from these Turing endpoints.}{thus we postpone this analysis.}

\begin{figure}[h]
\boxed{
\begin{minipage}{0.22\textwidth}
(a)\qquad $\gamma=0.01$\\
\includegraphics[width=1\textwidth]{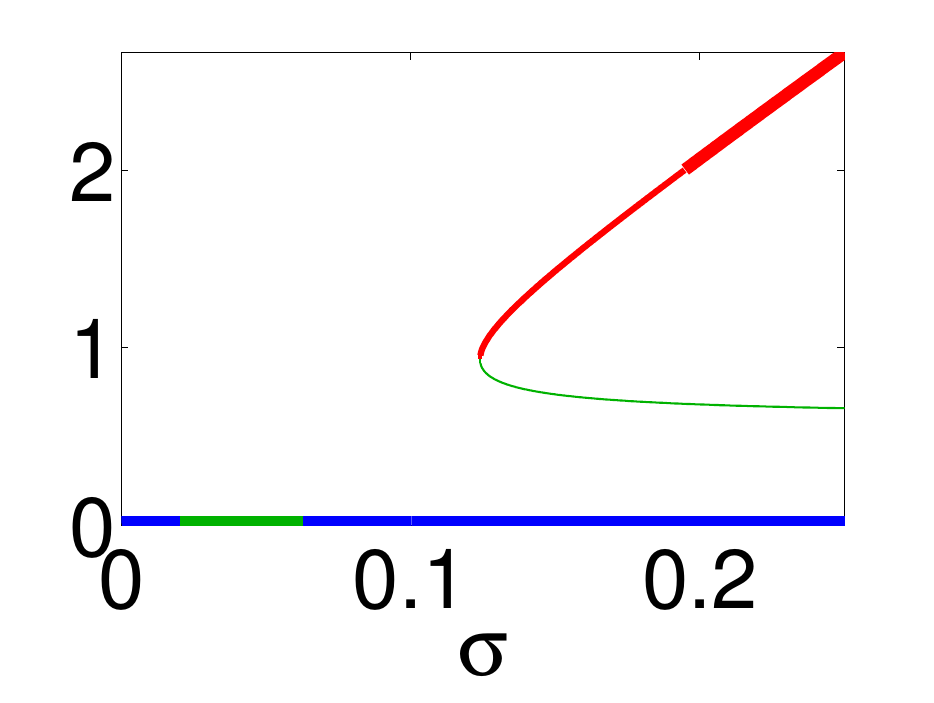}\\
\includegraphics[width=1\textwidth]{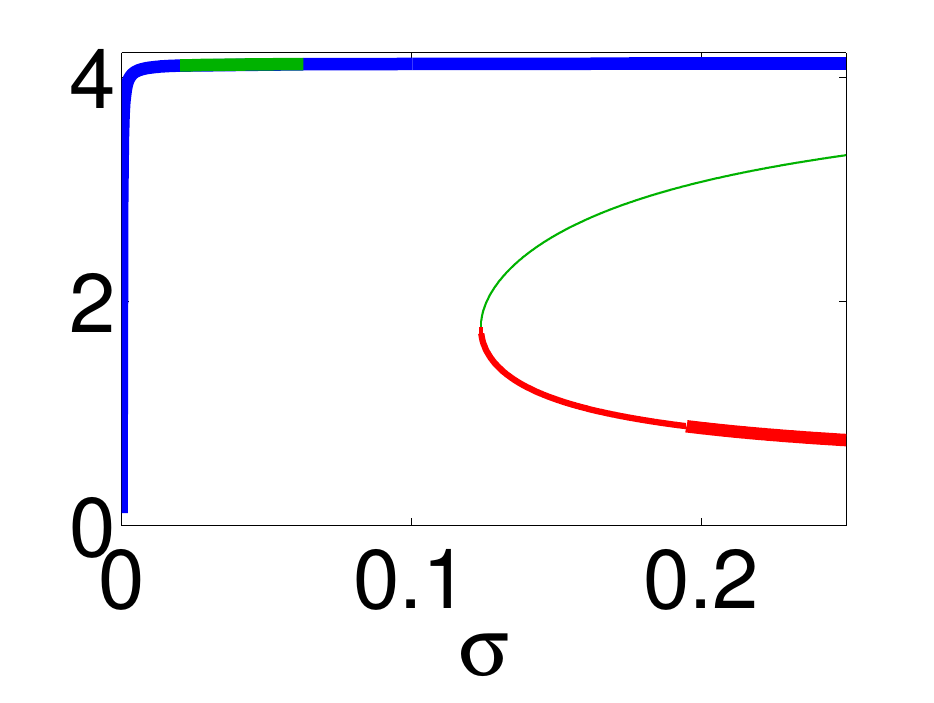}
\end{minipage} }
\boxed{
\begin{minipage}{0.22\textwidth}
(b)\qquad $\gamma=0.25$\\
\includegraphics[width=1\textwidth]{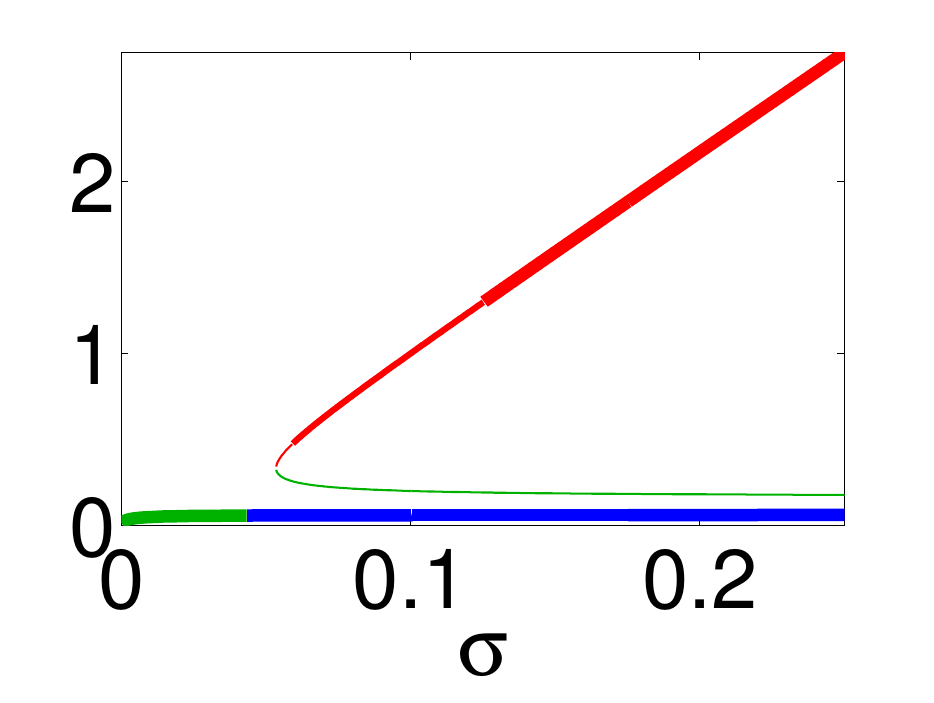}\\
\includegraphics[width=1\textwidth]{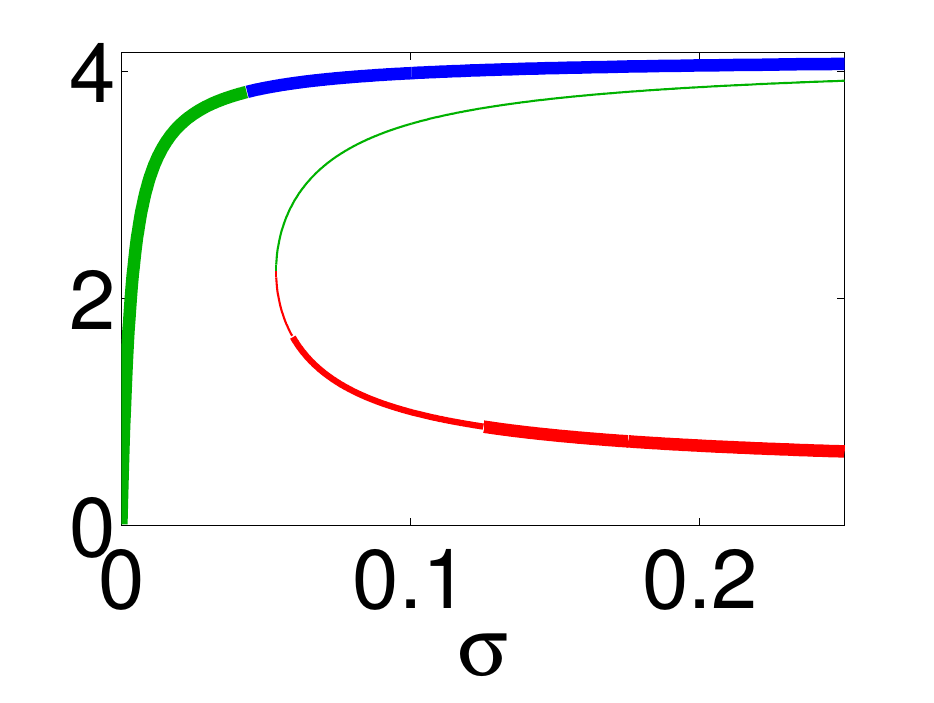}
\end{minipage} }
\boxed{
\begin{minipage}{0.22\textwidth}
(c)\qquad $\gamma=0.275$\\
\includegraphics[width=1\textwidth]{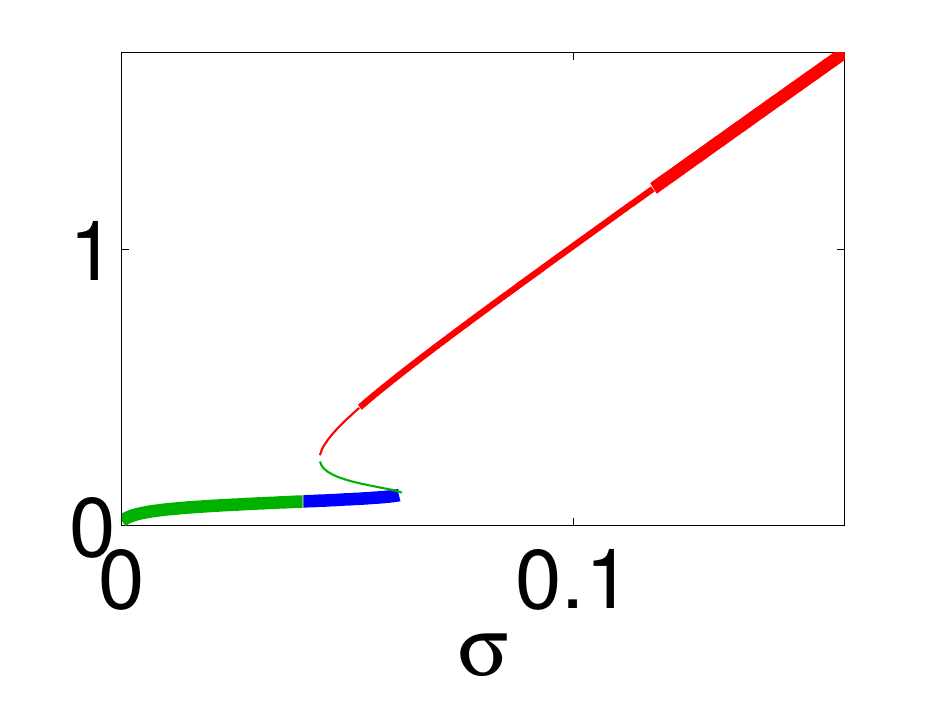}\\
\includegraphics[width=1\textwidth]{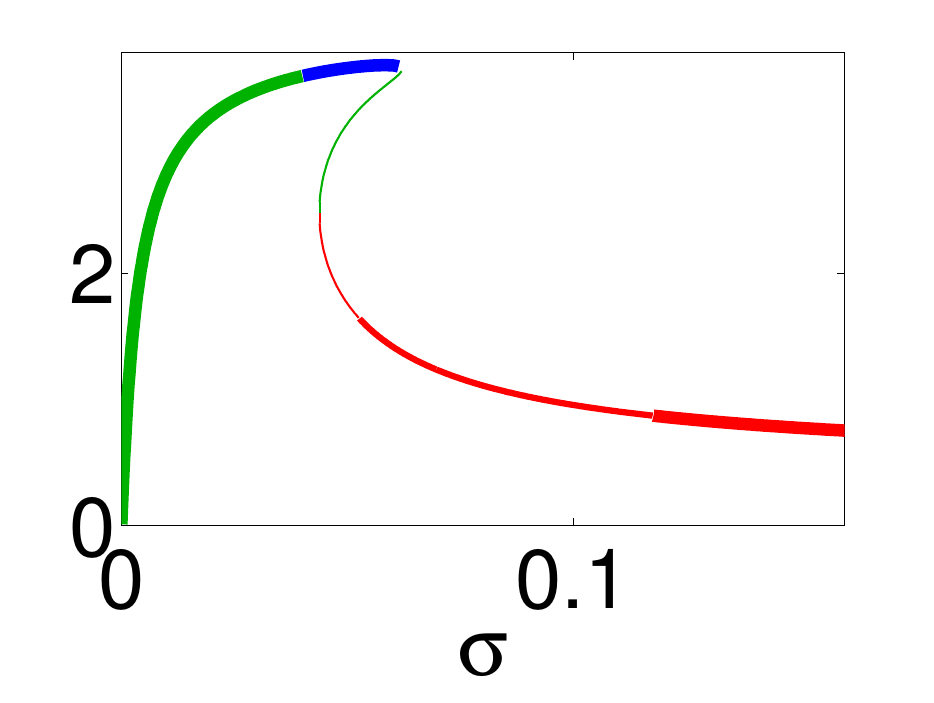}
\end{minipage} }
\boxed{
\begin{minipage}{0.22\textwidth}
(d)\qquad $\gamma=0.3$\\
\includegraphics[width=1\textwidth]{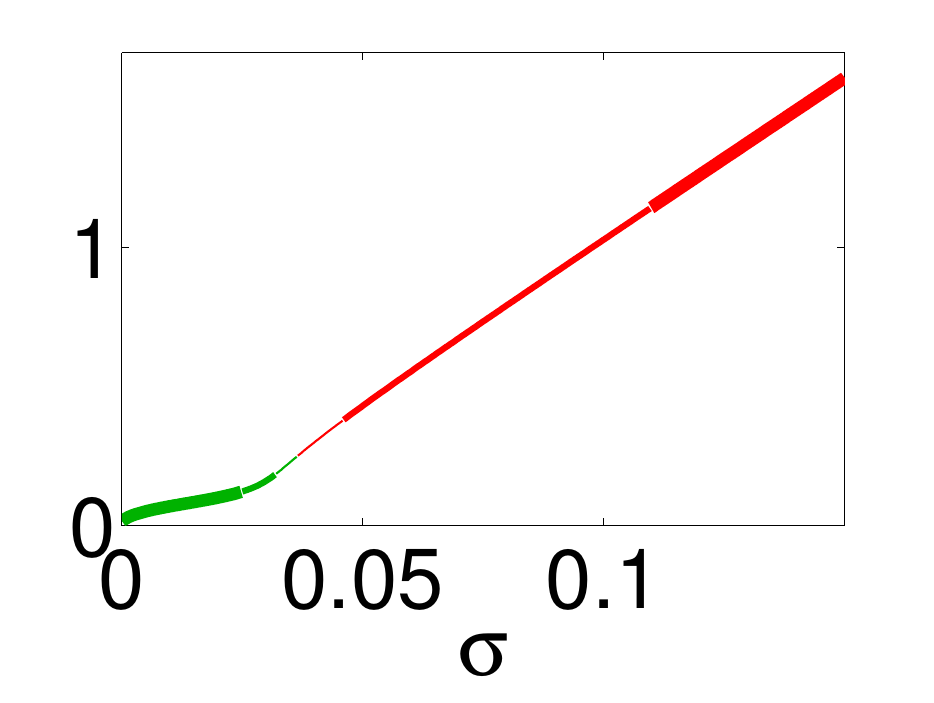}
\includegraphics[width=1\textwidth]{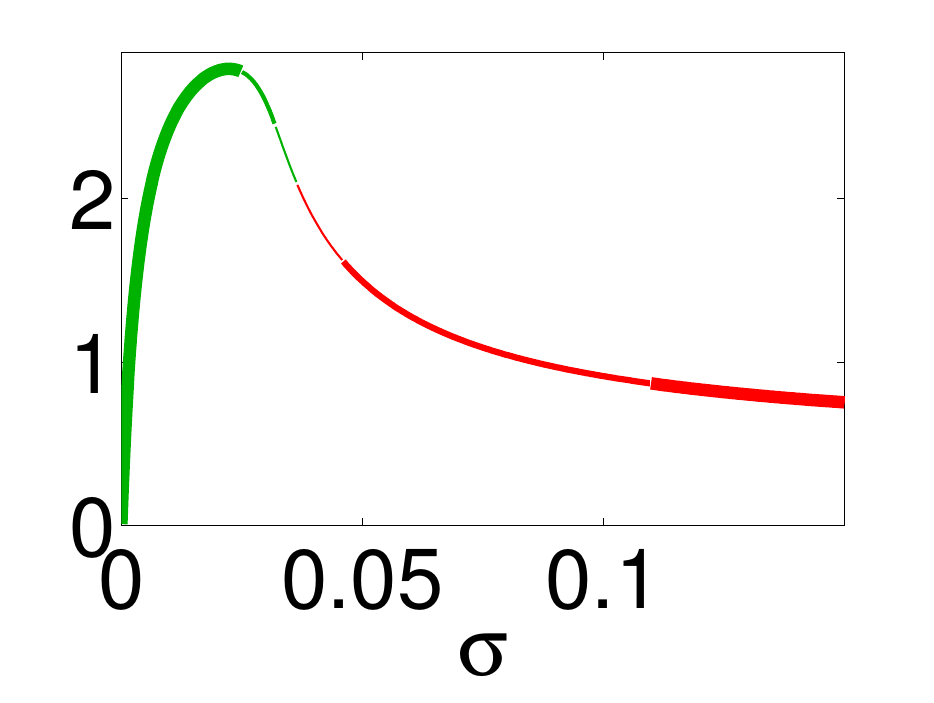}
\end{minipage} }

\caption{\deleted{Shown are the }Bifurcation diagrams of homogeneous positive real solutions for $\gamma=0.01$, 0.25, 0.275, and 0.3 in (a), (b), (c), and (d), respectively. Red, green, and blue lines represent $u_1$, $u_2$, and $u_3$ in the upper diagrams and $v_1$, $v_2$, and $v_3$ in the lower ones, respectively. Thick, medium, and thin lines represent stable, Turing-unstable, and space-independent-unstable solutions, respectively.}
\label{bifuhomplot}
\end{figure}

Example bifurcation diagrams for $I_1$ and $I_3$ can be seen in \figref{bifuhomplot} (a) and (d), respectively. \figref{bifuhomplot} (b) and (c) show diagrams for $I_2$.
\section{Landau reduction} \label{seclandau}
First we recall how the Landau reduction, which can also be found in \cite{pismen06,Hoyle} and is a \deleted{generalized }center manifold reduction \added{on a lattice}, works in general for systems of the form \eqref{rds}.  

 Let $w^*=(u^*,v^*)$ be a homogeneous solution of $f$ and $i,j$ two non-negative integers. We write $\partial^i_u \partial^j_v f(w^*)$  as $f_{ \underbrace{u\dotsc u}_{i \text{ times}} \underbrace{v\dotsc v}_{j \text{ times}}}$. We use an analogous notation for $g$ and $h$. 
Taylor-expanding $f$ around $(u^*,v^*)$ to third order and setting $(\tilde{u},\tilde{v})=(u,v)-(u^*,v^*)$, we obtain 
\begin{align*}
f(w)\approx & \sum_{ a,b\in\mathbb N_0 \atop  a+b\le 3 }  \partial^a_u \partial^b_v f(w^*) \frac{\tilde{u}^a \tilde{v}^b}{a! \ b!} \\
=&\underbrace{f_u \tilde{u}+f_v \tilde{v}}_{J_f(w^*)w^T}
+\frac{1}{2} f_{uu} \tilde{u}^2+f_{uv} \tilde{u}\tilde{v}+\frac{1}{2} f_{vv}\tilde{v}^2
+\frac{1}{6} (f_{uuu}  \tilde{u}^3+f_{vvv}  \tilde{v}^3)\\
&+\frac{1}{2} (f_{uuv} \tilde{u}^2  \tilde{v} + f_{uvv} \tilde{u}  \tilde{v}^2).
\end{align*}
Here we stop at order three, because we use the Landau reduction only up to third order. Substituting the expansion above into \eqref{rds}, the system becomes
\begin{align}
\partial_t  w = L(\Delta) w +B(w,w)+C(w,w,w),\label{dglredu}
\end{align} 
where $L(\Delta)=J_f(w^*)       
 +\left( \begin{smallmatrix} \Delta&0 \\ 0&\delta \Delta \end{smallmatrix} \right)$. $B$ and $C$ are symmetric bilinear and trilinear forms, respectively. For $p,q,r \in \R^2 $ they have the form:
\begin{align*}
&B(p,q)=\frac{1}{2} f_{uv}(p_1q_2+p_2q_1)+\frac{1}{2}(f_{uu}p_1 q_1+f_{vv}p_2q_2),\\
&C(p,q,r)=\frac{1}{6}(f_{uuu}p_1q_1r_1+f_{vvv}p_2q_2r_2)\\
&\quad\quad\quad\quad\quad+\frac{1}{6} \big( f_{uuv}(p_1q_1r_2+r_1p_1q_2+q_1r_1p_2)+f_{uvv}(p_1q_2r_2+r_1p_2q_2+q_1r_2p_2)\big).
\end{align*}
The eigenvalue $\mu_-$ is negative in Turing-unstable ranges, while the eigenvalue  $\mu_+$ has \replaced{non-negative}{ positive} parts. Thus we always consider $\mu_+$ in the following and write $\mu$ for simplicity.

The spot patterns shown in \figref{feudelfig} have a hexagonal structure, because every spot has six direct neighbors. Thus we start our Landau reduction on hexagonal lattices with the following ansatz
\begin{align}
w=\sss A_i e_i\Phi+c.c.=(A_1e_1+A_2e_2+A_3e_3) \Phi+c.c.\label{ansatz}
\end{align}
to reduce \eqref{dglredu} to a system of space-independent amplitudes.
Here $\Phi=\phi(\wn)$ is the eigenvector of $\hat{L}(k)$ which correspond to $\mu$, $e_j=\e^{i (x,y)\cdot \mathbf{k_j}}$, $A_j=A_j(t)\in\C$ for $j=1,2,3$,
\[ \mathbf{k_1}= \wn\left( \begin{smallmatrix} 1 \\ 0 \end{smallmatrix} \right), \quad \quad
  \mathbf{k_2}=\tfrac{\wn}{2}\left( \begin{smallmatrix} -1 \\ \sqrt3 \end{smallmatrix} \right), \quad \text{and} \quad
 \mathbf{k_3}=\tfrac{\wn}{2}\left( \begin{smallmatrix} -1 \\ -\sqrt3 \end{smallmatrix} \right).\]
\begin{figure}[h]
\centering{
\begin{minipage}{0.4\textwidth}
\includegraphics[width=1\textwidth]{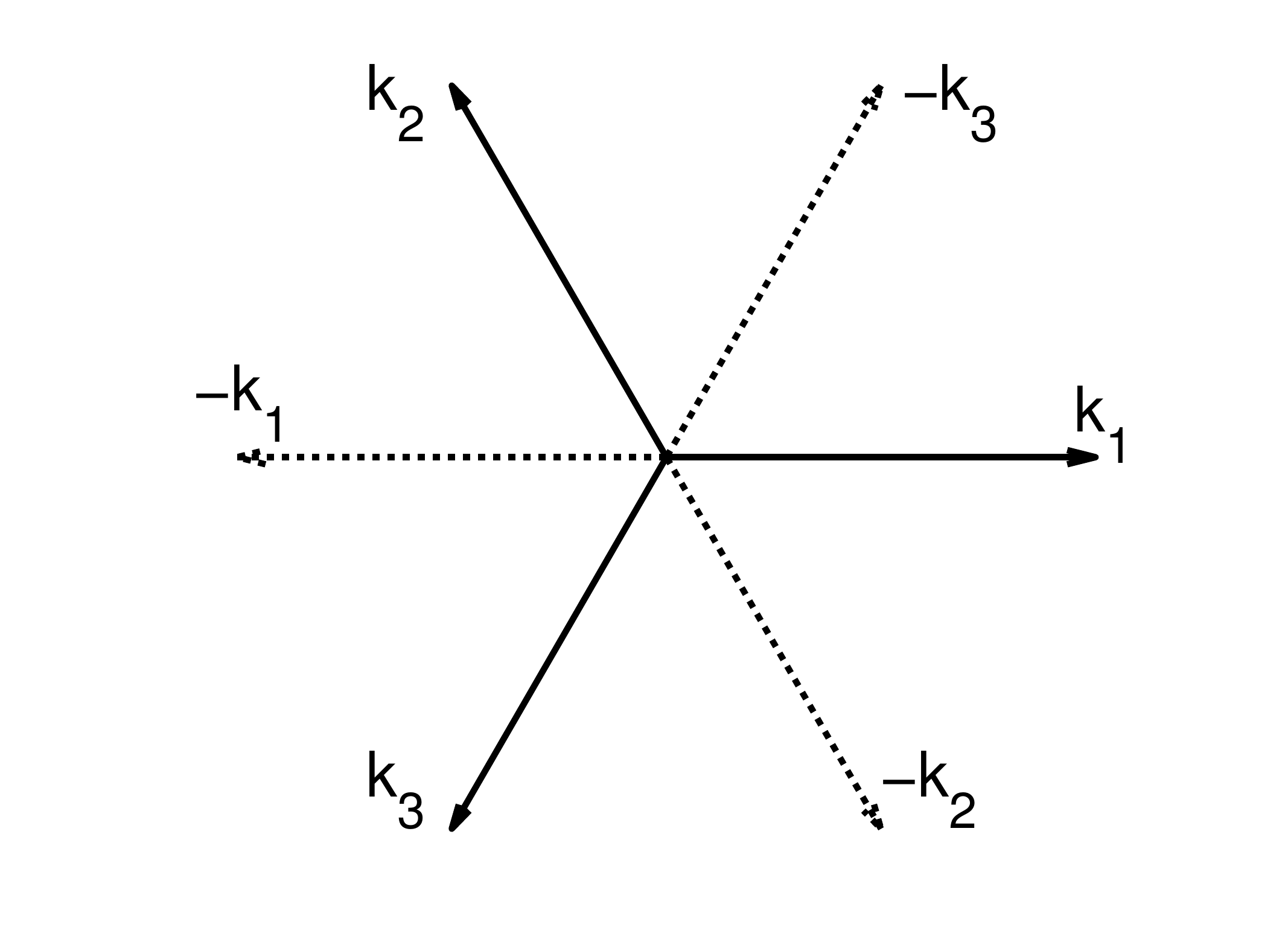}
\end{minipage} 
}
\caption{Sketch of the vectors $\mathbf{k_1}, \ \mathbf{k_2}, \ \mathbf{k_3}$.     }
\label{vec}
\end{figure}

We choose the length (wavenumber $\wn$) of the wave vectors $\mathbf{k_1}$, $\mathbf{k_2}$, $\mathbf{k_3}$ such that there is a bifurcation parameter value $\sigma_\wn$ for which the curve of eigenvalues has a zero in $\wn$. 
It holds
 \begin{align*}
 \partial_t w=\sss \partial_t A_ie_i \Phi+c.c., \quad
 L(\Delta)w= \mu(\wn) \sum_{i=1}^3 A_i  e_i \Phi+c.c.,
\end{align*}
 \begin{align*}
 B(w,w)
 = \ & \bigg(\Big(\sss A_i^2 e_i^2\Big)B(\Phi,\Phi)+\big(2 \overline{A_1} \overline{A_2}e_3+2 \overline{A_1} \overline{A_3}e_2+2 \overline{A_2} \overline{A_3} e_1\big) B(\overline{\Phi},\overline{\Phi}) \\
 &+\Big(\sss |A_i|^2 +\sum_{1\le i<j \le 3} 2 A_i\overline{A_j} e_i \overline{e_j}\Big)B(\Phi,\overline{\Phi})\bigg)+c.c.,
 \end{align*}
 and
 \begin{align*}
 C(w,w,w)=  \ &\Big((A_1e_1+A_2e_2+A_3e_3)^3  C(\Phi,\Phi,\Phi)\\
 &+3(A_1e_1+A_2e_2+A_3e_3)(\overline{A_1e_1+A_2e_2+A_3e_3})^2  C(\Phi,\overline{\Phi},\overline{\Phi})\Big)+c.c..
 \end{align*}
 To eliminate quadratic terms from the residual $\text{Res}(w)=-\partial_t  w + L(\Delta) w +B(w,w)+C(w,w,w)$ which do not correspond to modes $e_1$, $e_2$, or $e_3$, we extend the ansatz \eqref{ansatz} to
 \begin{align}\label{newansatz}
 w=& \Big( \sss A_i e_i \Phi+  \sss A_i^2 e_i^2\phi_{ii}+\frac12 \sss |A_i|^2 \phi_0+ \sum_{1\le i < j \le 3} A_i \overline{A}_j e_i \overline{e_j} \phi_{ij}\Big)+c.c..
 \end{align}
 Notice that $\phi_{ii}$ and $\phi_{ij}$ are independent of $i$ and $j$. Using this notation, one can see that $\phi_{ii}$ and $\phi_{ij}$ correspond to $A_i^2$ and $A_i\overline{A_j}$, respectively. 
 Substituting \eqref{newansatz} into \eqref{dglredu} and sorting with respect to $e_m^n$, yields
 \begin{alignat*}{4}
 & A_i^2 e_i^2: & \qquad & B(\Phi,\Phi)+L(2\wn)\phi_{ii}=0 & \qquad & \Rightarrow & \qquad & \phi_{ii}=-L(2\wn)^{-1}B(\Phi,\Phi)\\
 & |A_i|^2:& \qquad &2B(\Phi,\overline{\Phi})+L(0)\phi_0=0 &\qquad & \Rightarrow &\qquad & \phi_0 =-2L(0)^{-1}B(\Phi,\overline{\Phi}) \\
 & A_i \overline{A_j} e_i \overline{e_j}: & \qquad &  2B(\Phi,\overline{\Phi})+L(\sqrt{3}\wn) \phi_{ij}=0 & \qquad & \Rightarrow & \qquad & \phi_{ij} =-2L(\sqrt{3}\wn)^{-1}B(\Phi,\overline{\Phi})
\end{alignat*}

\noindent To remove terms of order $e_i$ from the residual, we extend the ansatz \eqref{newansatz} to $\tilde{w}=w+\sum_{i=1}^3 \phi_{3i} e_i$. Substituting $\tilde{w}$ into \eqref{dglredu} and sorting with respect to $e_1, \ e_2, \ e_3$ yields
\begin{align}
  e_1: \quad -\hat{L}(\wn)\phi_{31}=&-\partial_t A_1+ d_{1} A_1+d_2 \overline{A_2}\overline{A_3} +d_3 A_1 |A_1|^2 \nonumber \\
  &+d_4 A_1(|A_2|^2 +|A_3|^2)+R_1, \nonumber \\
  e_2:  \quad   -\hat{L}(\wn)\phi_{32}=&-\partial_t A_2+d_1 A_2+d_2 \overline{A_1}\overline{A_3} +d_3 A_2 |A_2|^2 \label{glpr} \\
  &+d_4 A_2(|A_1|^2 +|A_3|^2)+R_2, \nonumber \\
 e_3:  \quad   -\hat{L}(\wn)\phi_{33}=&-\partial_t A_3+d_1 A_3+d_2 \overline{A_1}\overline{A_2} +d_3 A_3 |A_3|^2 \nonumber \\
 &+d_4 A_3(|A_1|^2 +|A_2|^2)+R_3, \nonumber 
 \end{align}
with
 \begin{align}
d_1&= \mu(\wn) \Phi,\nonumber  \\
d_2&= 2   B(\overline{\Phi},\overline{\Phi}), \nonumber  \\
d_3&= 3C(\Phi,\Phi,\overline{\Phi})+2B(\overline{\Phi},\phi_{ii})+2B(\Phi,\phi_{0}), \nonumber \\
d_4&= 6C(\Phi,\Phi,\Phi)+2B(\Phi,\phi_{ij})+2B(\Phi,\phi_{0}). \nonumber
 \end{align}
 The summands $R_1, \  R_2, \  R_3$ represent all higher order terms, e.g., $A_1|A_2|^4$ is a term of $R_1$. By the Fredholm alternative there exists a solution for \eqref{glpr} iff every equation of \eqref{glpr} is an element of ker($\hat{L}(\wn)^H)^\bot$. Let $\Phi^*$ be the adjoint eigenvector of $\hat{L}(\wn)$ to the eigenvalue $\mu(\wn)$ evaluated in $\sigma_\wn$, i.e., $\hat{L}(\wn)^H \Phi^*=\overline{\mu(\wn)} \Phi^*$, and let $\Phi^*$ be normalized such that $\langle \Phi, \Phi^* \rangle=1$. Multiplying \eqref{glpr} with $\Phi^*$ and setting $R_1= R_2=R_3=0$, yields
\begin{align}
 &\ & \quad & \partial_t A_1=c_1 A_1+c_2 \overline{A_2}\overline{A_3} +c_3 A_1 |A_1|^2+c_4 A_1(|A_2|^2 +|A_3|^2), \nonumber \\
 &\ & \quad &  \partial_t A_2=c_1 A_2+c_2 \overline{A_1}\overline{A_3} +c_3 A_2 |A_2|^2+c_4 A_2(|A_1|^2 +|A_3|^2), \label{gl}\\
 &\ & \quad &  \partial_t A_3=c_1 A_3+c_2 \overline{A_1}\overline{A_2} +c_3 A_3 |A_3|^2+c_4 A_3(|A_1|^2 +|A_2|^2), \nonumber 
\end{align}
where $c_i=\langle d_i, \Phi^* \rangle$. The classical Landau reduction evaluates $c_2$, $c_3$, $c_4$, $\Phi$, $\phi_{ii}$, $\phi_{ij}$, $\phi_0$ in $\sigma_c$ and $c_1$ in $\sigma$. One can see in \cite{schnaki} that \replaced{including the $\sigma$-dependence}{ evaluating them all in $\sigma$} can give better approximations. However, in the following we use the Landau reduction to predict the existence and bifurcation directions of stationary states which branch from homogeneous solutions. To approximate solutions and follow their branches we use {\tt pde2path}. 

This method yields the same system \eqref{gl} and Landau coefficients if we use phase shifted space coordinates in the ansatz, i.e., using $(x+\psi,y+\psi)$ with $\psi \in (0,2\pi)$ instead of $(x,y)$. 
Later we will use numerical methods to find solutions of \eqref{dgl} on bounded domains with Neumann boundary conditions. Because of the Neumann boundary conditions we are not able to find all phase shifts of a solution.\\
\added{It is possible to perform a reduction to \eqref{gl}, which is valid in the sense of the center manifold theorem if we find a point in the parameter space for which $c_2=0$. In this case one is able to give consistent results for small $c_2$ via a codimension-two bifurcation. The problem is that this is an unnatural case in applications. Normally one is interested in a codimension-one bifurcation for which the coefficient $c_2$ is not small.}\\
It is not shown yet that this method gives an approximation of a solution of the full PDE \eqref{rds}. Furthermore, it is unclear whether one can conclude the stability from the reduced system \eqref{gl}. Comparisons between solutions found via this presented Landau reduction and the finite element method can be found in \cite{schnaki}. One can see there that the Landau reduction gives acceptable approximations for the specific reaction diffusion system, which is considered there if the amplitudes are small. In the following we will use the Landau reduction to predict existence and stability of states near the onset and will see that these predictions fit well to the numerical results.\\
\deleted{It is possible to perform a reduction to \eqref{gl}, which is valid in the sense of the center manifold theorem if we find a point in the parameter space for which $c_2=0$. In this case one is able to give consistent results for small $c_2$ via a codimension-two bifurcation. The problem is that this is an unnatural case in applications. Normally one is interested in a codimension-one bifurcation for which the coefficient $c_2$ is not small.} 
\section{1D Patterns}\label{sec1d}
First we consider solutions over one dimensional domains. Here the modes $e_2$ and $e_3$ do not exist. Hence, the system \eqref{gl} reduces to
\ali{
\partial_t A=c_1 A +c_3 |A|^2A. \label{1dredut}
}
Stationary amplitudes solve 
\ali{
c_1 A +c_3 |A|^2A=0. \label{1dredu}
}
Clearly, $A=0$ solves \eqref{1dredu}. Inserting this solution into the ansatz \eqref{newansatz}, yields the homogeneous solution. More interesting are the second type of solutions which we obtain from \eqref{1dredu}. They fulfill 
\ali{|A|=\sqrt{- \frac{\mu(\wn)}{c_3}} \label{ff} }
and generate periodic solutions by substituting  $(A_1,A_2,A_3)=(A,0,0)$ into \eqref{newansatz}. This type of solution exists also in 2D and we call these solutions stripes because of their 2D-density plot. When we use numerical methods to determine the stripes in the following, we use Neumann boundary conditions. Stripes which fulfill Neumann boundary conditions over a domain $(-l\pi/ \wn,l\pi/ \wn)$ with $l\in \N$ correspond to the amplitudes 
\alinon{S_\pm=\pm \sqrt{- \frac{\mu(\wn)}{c_3}}. }
If the bacteria density $u$ has its maximum (minimum) in $x=0$ for such a stripe solution, we call it hot (cold) stripes. Notice that we do not automatically have hot and cold stripes for $S_+$ and $S_-$, respectively. All other amplitudes which fulfill \eqref{ff} generate phase shifts of the hot resp. cold stripes.

Let $\sigma_\wn$ be a balancing rate and $\wn\in\R^+$ a wavenumber such that the curve of eigenvalues $\mu$ for $\sigma_\wn$ has a single zero in $\wn$. Let $\sigma_s$ and $\sigma_l$ be two balancing rates, which are \replaced{sufficiently}{ infinitesimally} smaller and larger than $\sigma_k$, respectively. It holds that the curve of eigenvalues $\mu$ is positive in $\wn$ for $\sigma_s$ or $\sigma_l$, while it is negative for the other. The stripes can only exist, where \eqref{ff} is fulfilled such that the algebraic sign of $c_3$ evaluated in $\sigma_\wn$ tells us in which direction the stripes bifurcate.
For $\sigma_c$ we have a double zero in $k_c$. For a bifurcation of stripes it is necessary that $\mu(k_c)>0$ in $\sigma_l$ or $\sigma_r$.  

\subsection{Changing wavelength} \label{secchanging}
Let us consider the system \eqref{dgl2} for $\gamma=0.3$. From our analysis above we already know that we have only one homogeneous solution for all $\sigma\in (0,0.25)$ with an unstable range bounded by Turing bifurcation points. Furthermore, we know that there are two Turing endpoints in this unstable range. We use our data set  of \figref{stabplot} to find out that the right and left Turing bifurcation points are given by $\sigma_c^r\approx 0.11$ and $\sigma_c^l \approx 0.025$ with corresponding critical wave numbers $k_c^r\approx 0.187$ and $k_c^l\approx 0.067$, respectively. The curve of eigenvalues $\mu_\pm$ (see \eqref{mu}) has a zero in $k$, when $\text{det}L(k)=0$. This is the case for 
\alinon{
k_\pm= \sqrt{\frac{d g_u+h_v}{2d} \pm \sqrt{\left(\frac{d g_u+h_v}{2d}\right)^2+\frac{g_vh_u-g_uh_v}{d}}}.
} 
\begin{figure}[h]
\begin{minipage}{0.49\textwidth}
(a) global 1D bifurcation diagram\\
\includegraphics[width=1\textwidth]{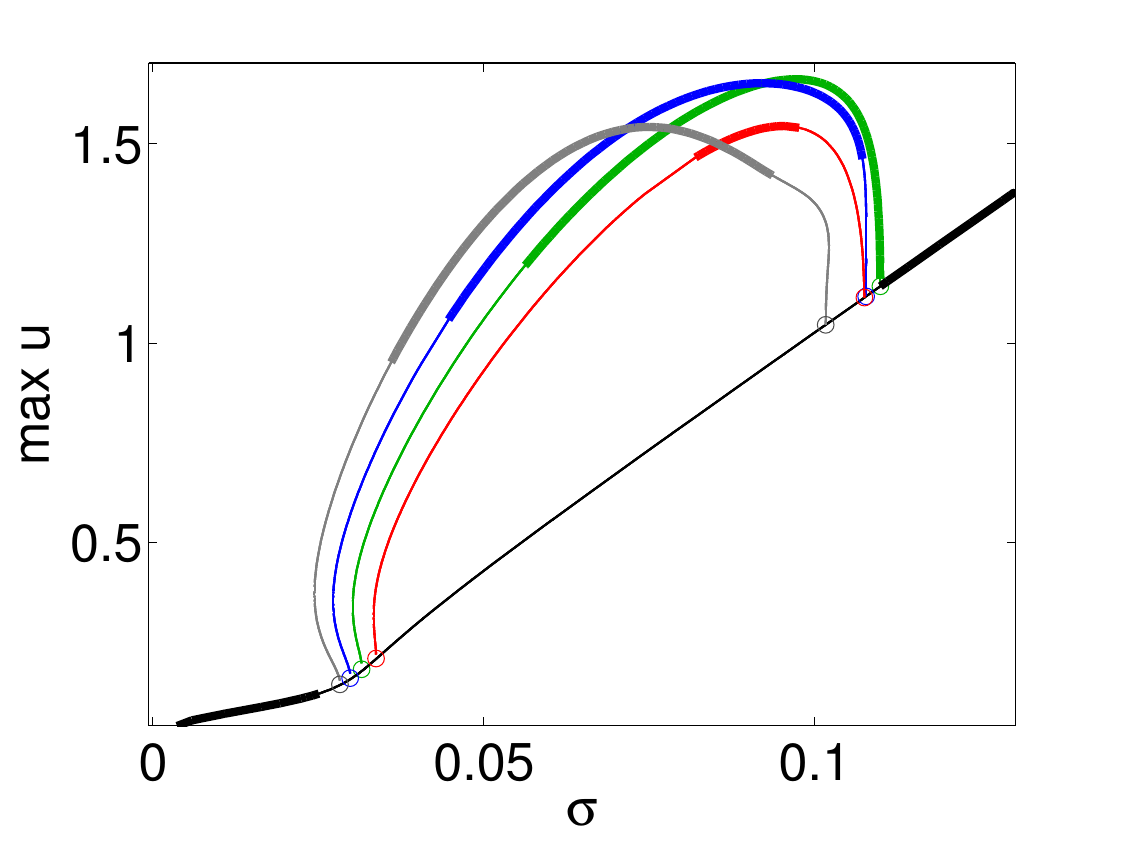}
\end{minipage} 
\begin{minipage}{0.49\textwidth}
(b) $u$ of {\tt R4} for $\sigma=0.08$ \\
\includegraphics[width=1\textwidth]{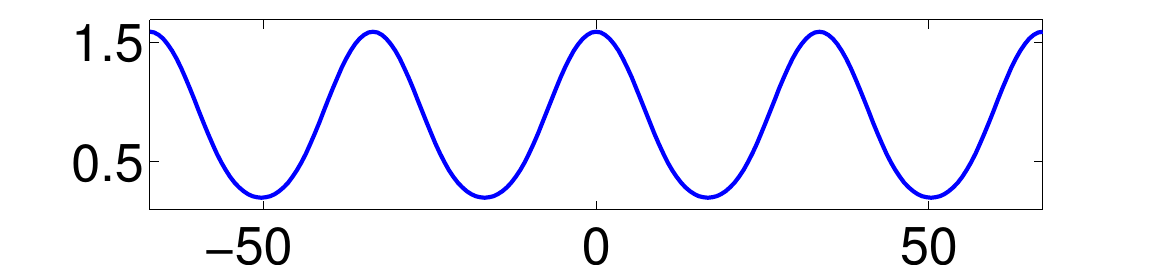}\\
\vspace{0.1cm}

(c) $u$ of {\tt R35} for $\sigma=0.09$ \\
\includegraphics[width=1\textwidth]{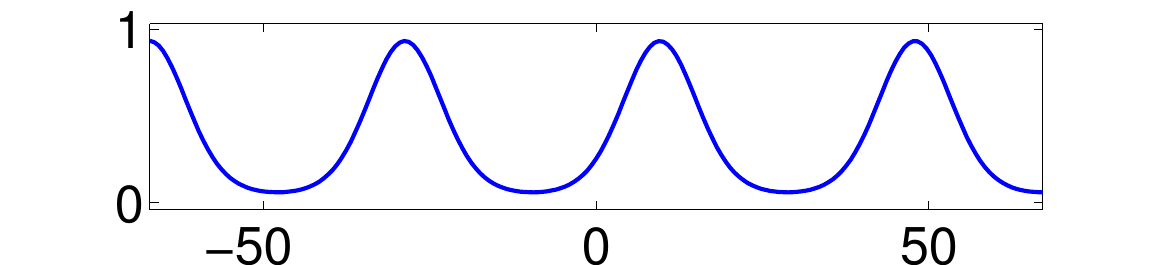}
\end{minipage}

\begin{minipage}{0.49\textwidth}
(d) $u$ of {\tt R45} for $\sigma=0.04$\\
\includegraphics[width=1\textwidth]{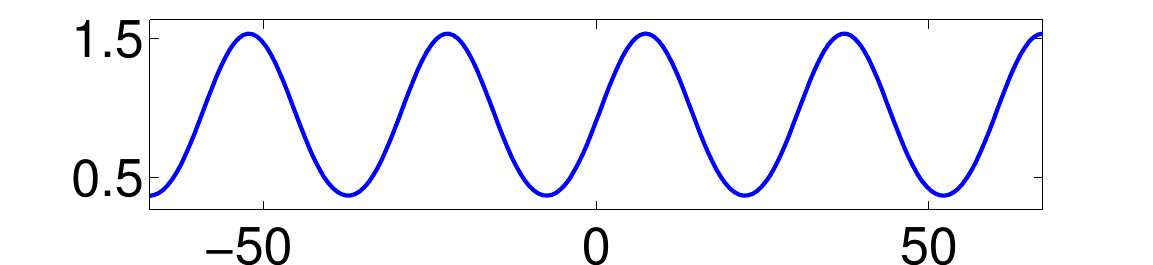}
\end{minipage}
\begin{minipage}{0.49\textwidth}
(e) $u$ of {\tt R3} for $\sigma=0.06$\\
\includegraphics[width=1\textwidth]{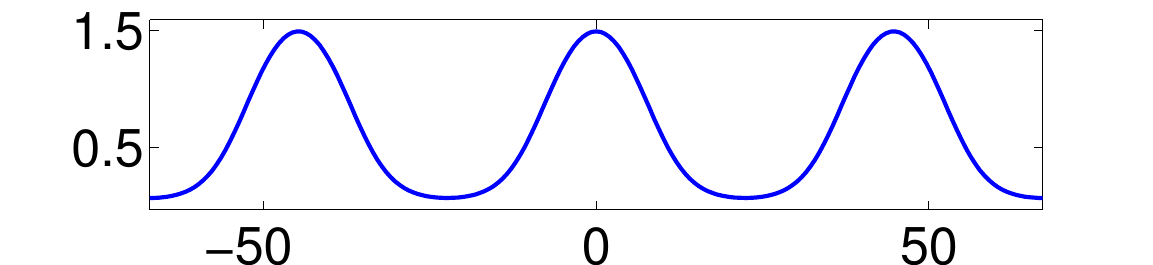}
\end{minipage} 
\caption{All plots are for $\gamma=0.3$ and the domain $\Omega^r=(-4\pi/k_c^r,4\pi/k_c^r)$. (a) From the first, second, third, and fourth bifurcation point (counting from the right side) of the homogeneous branch (black) bifurcate branches of periodic solutions with 4 (green), 3.5 (blue), 4.5 (red), and 3 (gray) periods, which we call {\tt R4}, {\tt R35}, {\tt R45}, and {\tt R3}, respectively. Example solutions of {\tt R4}, {\tt R35}, {\tt R45}, and {\tt R3} are shown in (b), (c), (d), and (e), respectively. All these example solutions are stable. Here and for all coming calculations for which we used {\tt pde2path} we always use Neumann boundary conditions, and thick and thin parts of the branches represent stable and unstable solutions, respectively.  }
\label{fromright}
\end{figure}

Clearly, in the Turing unstable range it holds $k_\pm\in \R$. We checked that this is also the case between the Turing endpoints. 
The Landau formalism above predicts that periodic solutions of the type 
\ali{(u,v)=(u^*,v^*)+2A\cos(k x)\Phi+\text{h.o.t.}\label{stripecos}}
bifurcate from $(u^*,v^*)$ at $\sigma_c^r$ with $k=k_c^r$ if we consider the problem over the 1D domain $\Omega^r=(-4\pi/k_c^r,4\pi/k_c^r)$. We are able to prove analytically that there are balancing rates $\sigma_{35}$, $\sigma_{45}$, $\sigma_{3}$ with $\sigma_c^r>\sigma_{35}>\sigma_{45}>\sigma_{3}$, where branches $R35,$ $R45,$ and $R3$ of periodic solutions of the type \eqref{stripecos} bifurcate with $k=3.5 k_c^r/4$, $4.5 k_c^r/4$, and $3 k_c^r/4$, respectively. These branches and some example solutions are shown in \figref{fromright}.

\begin{figure}[h]
\centering{
\begin{minipage}{0.249\textwidth}
(a) \\
\includegraphics[width=1\textwidth]{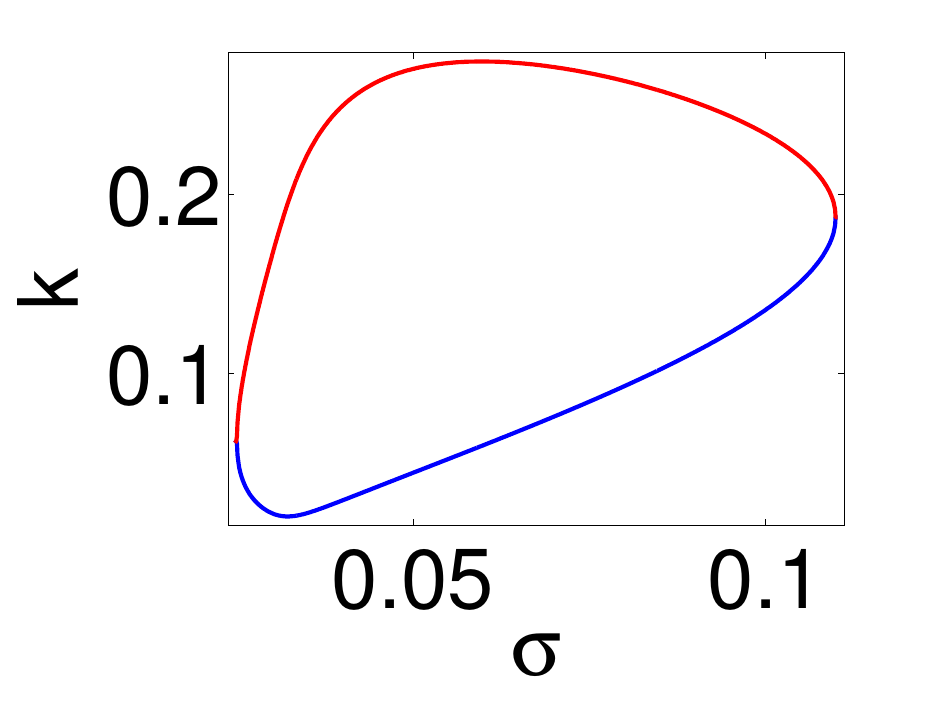}
\end{minipage} 
\begin{minipage}{0.249\textwidth}
(b) \\
\includegraphics[width=1\textwidth]{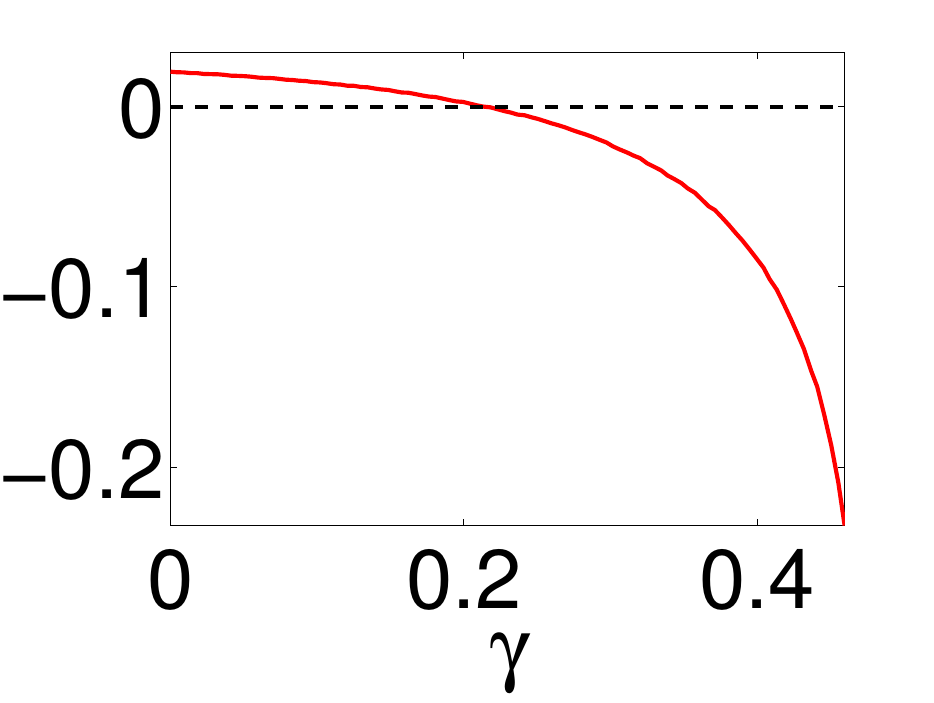}
\end{minipage} 
\begin{minipage}{0.249\textwidth}
(c) \\
\includegraphics[width=1\textwidth]{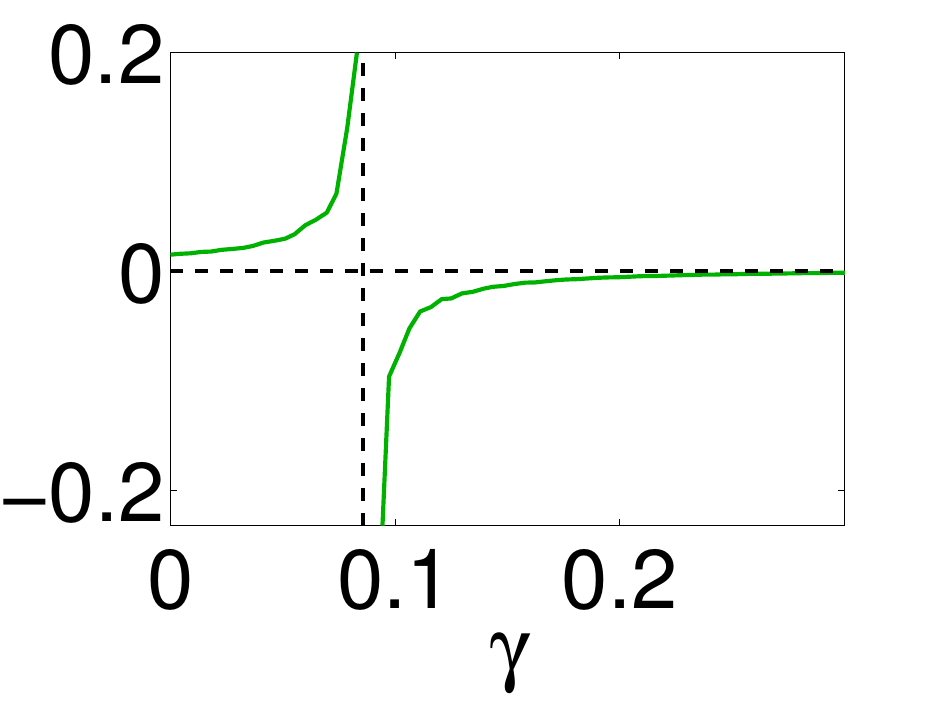}
\end{minipage} 
}
\caption{(a) The blue and red curves represent $k_-$ and $k_+$ for $\gamma=0.3$. The coefficients $c_3$ and  $c_f=c_2^2/(4(c_3+2c_4)^2)$ evaluated in $\sigma_c$ as function of $\gamma$ are shown in (b) and (c). They are used to predict the strength of the subcriticality of stripes and hexagons, respectively. Hexagon patterns and the role of $c_f$ will be introduced in \secref{2d}.    }
\label{kandc3}
\end{figure}

\begin{figure}[h]
\begin{minipage}{0.49\textwidth}
(a)\\
\includegraphics[width=1\textwidth]{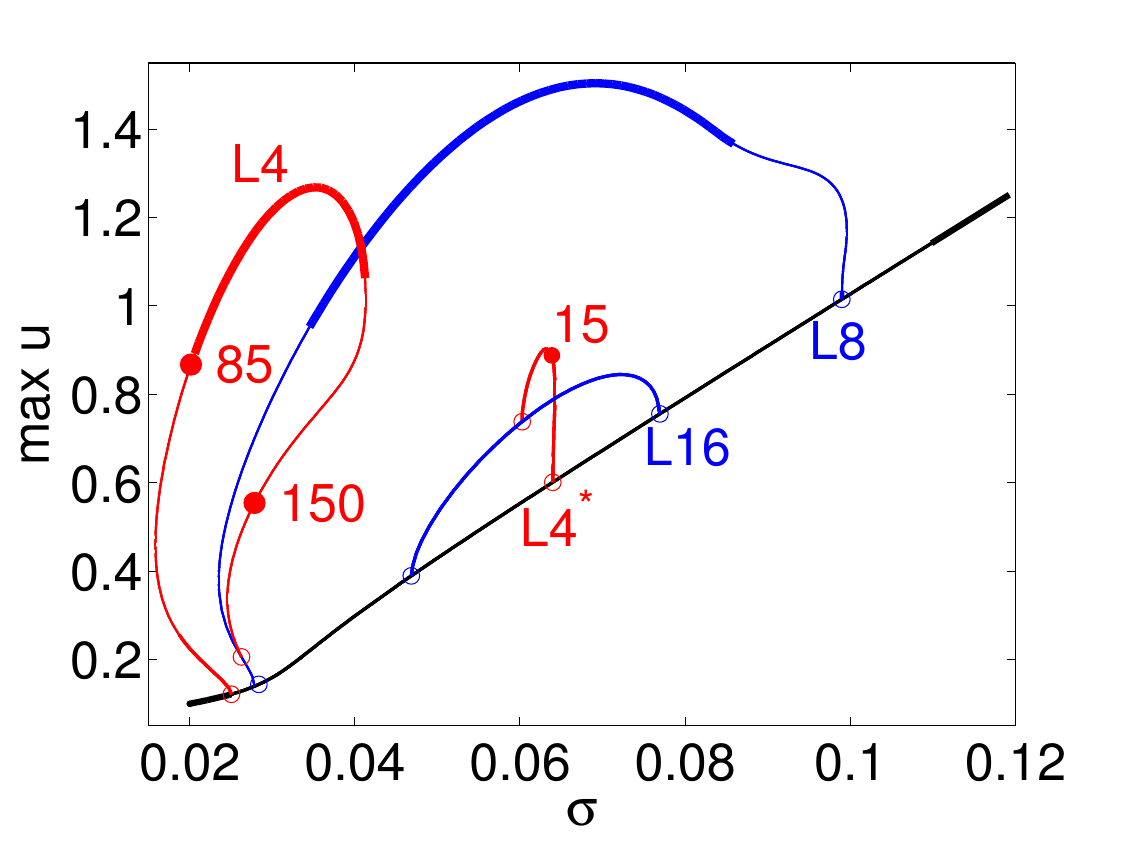}
\end{minipage} 
\begin{minipage}{0.49\textwidth}
(b) u at points 85 and 150 of {\tt L4}\\
\includegraphics[width=1\textwidth]{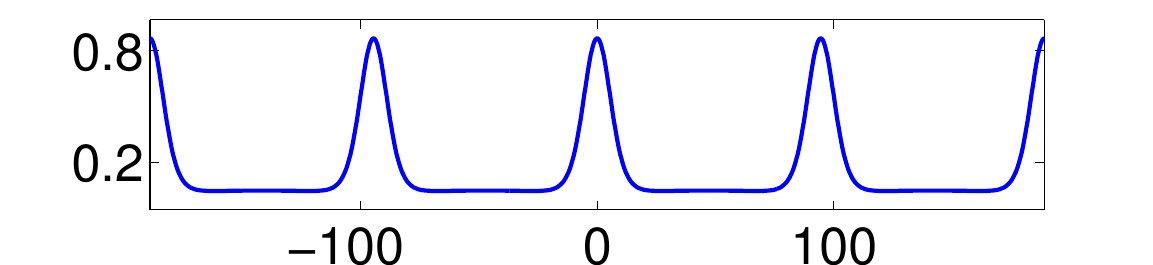}\\
\includegraphics[width=1\textwidth]{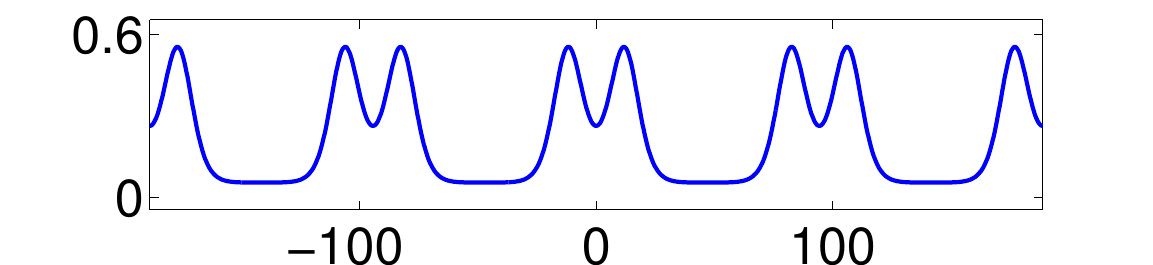}\\
(c) u at point 15 of {\tt L4$^*$}\\
\includegraphics[width=1\textwidth]{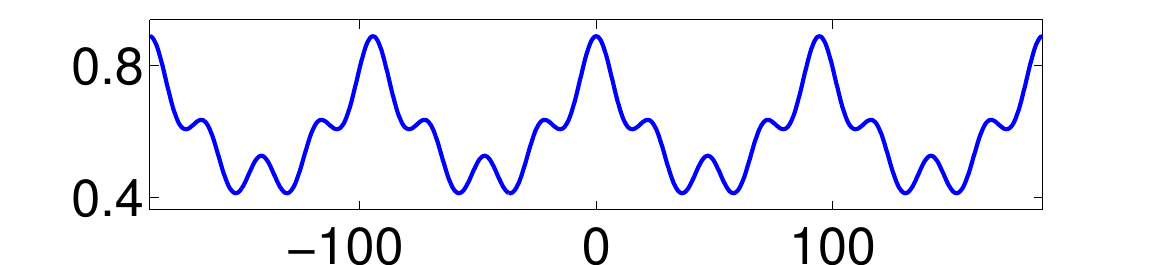}
\end{minipage} 
\caption{All plots are for $\gamma=0.3$ and the domain $\Omega^l=(-4\pi/k_c^l,4\pi/k_c^l)$. (a) On {\tt L8} and {\tt L16} we always have solutions of 8 and 16 periods, respectively. Both branches bifurcate and terminate on the homogeneous solution branch (black). {\tt L4} and {\tt L4$^*$} bifurcate from the homogeneous solution branch as periodic solutions with 4 periods and terminate on {\tt L8} and {\tt L16}, respectively. Example solutions of {\tt L4} and {\tt L4$^*$} are shown in (b) and (c), respectively.      }
\label{loop}
\end{figure}

Let $s\in\{\sigma_c^r, \ \sigma_{35}, \ \sigma_{45}, \ \sigma_3 \}$ and $\kappa$ be the corresponding wavenumber. 
In \figrefa{kandc3} we see that a balancing rate $s_2\neq s$ in the unstable range exists for which the eigenvalue curve $\mu_\pm$ has a real zero in $\kappa$ such that a stripe solution with wavenumber $\kappa$ also branches in $s_2$. 
By using numerical methods to follow the branches, which bifurcate in $s$, we see that they terminate in $s_2$ (see \figref{fromright}). 
One may conjecture that stripe branches which correspond to the same wavenumber are connected, but this is not always the case. We also computed the branch {\tt L4} which bifurcates in $\sigma_c^l$ with the critical wavenumber $k_c^l$ over the domain $\Omega^l=(-4\pi/k_c^l,4\pi/k_c^l)$. In \figrefa{kandc3} we see that there is a $\sigma\neq\sigma_c^l$, where stripe solutions of the wavelength $k_c^l$ bifurcate. We call the corresponding solution branch {\tt L4$^*$}. In \figref{loop} we see that {\tt L4} and {\tt L4$^*$} are not connected, but they connect to bifurcations on {\tt L8} and {\tt L16}, which are stripe solutions of 8 and 16 periods, respectively. One might guess that this depends on the Turing endpoints. However, we also compute the branches for $\gamma=0.4$ (not shown). Here we have no Turing endpoints, but the same effects.

\subsection{Localized patterns and snaking}
\begin{figure}[h]
\begin{minipage}{0.49\textwidth}
(a)\\
\includegraphics[width=1\textwidth]{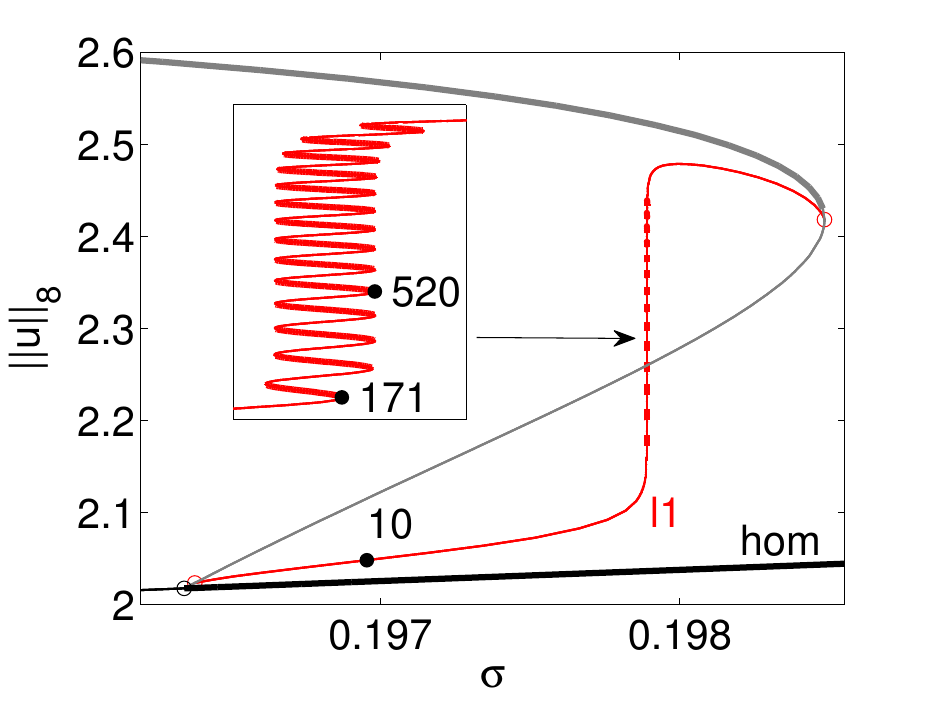}
\end{minipage} 
\begin{minipage}{0.49\textwidth}
(b)\small{ $u$ at points 10, 171, and 520 as indicated in (a) }\\
\includegraphics[width=1\textwidth]{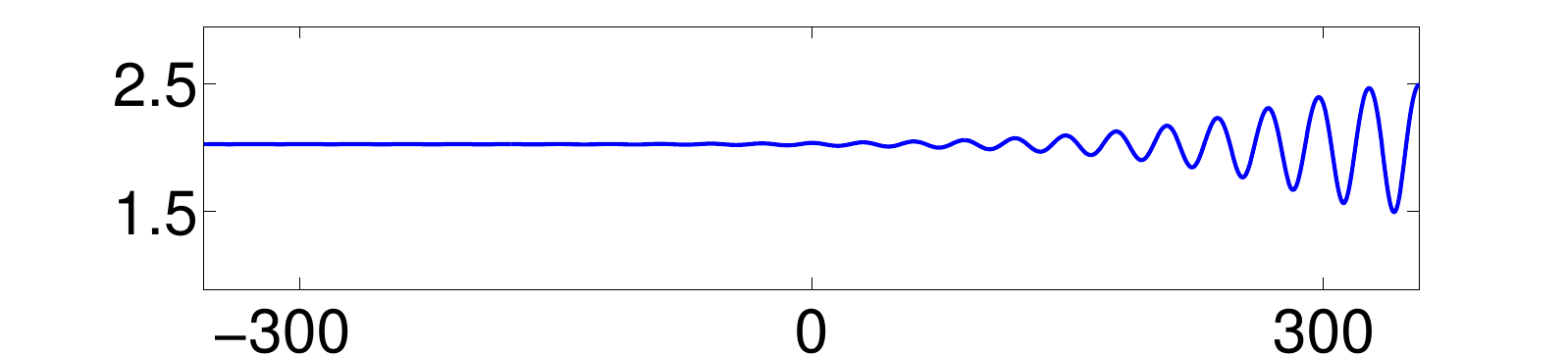}
\includegraphics[width=1\textwidth]{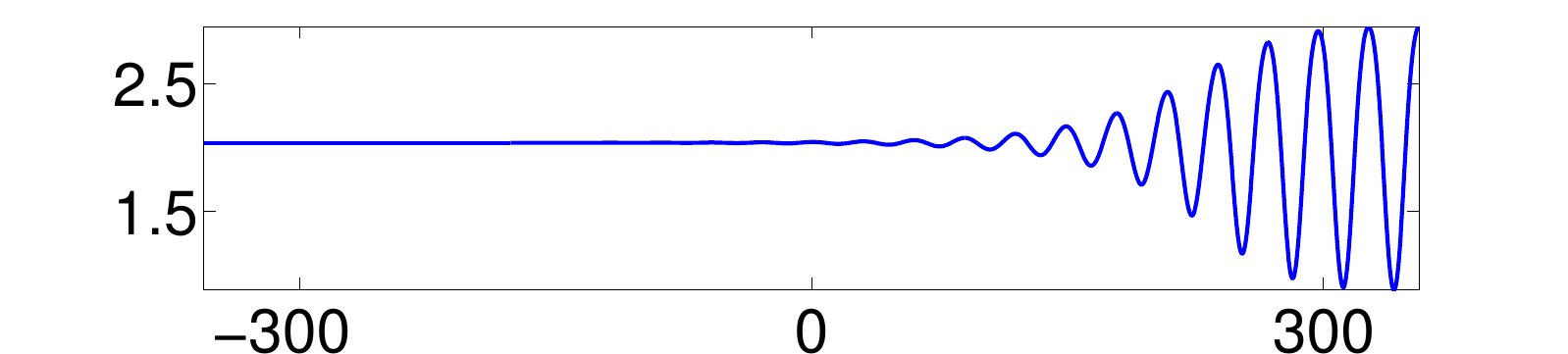}
\includegraphics[width=1\textwidth]{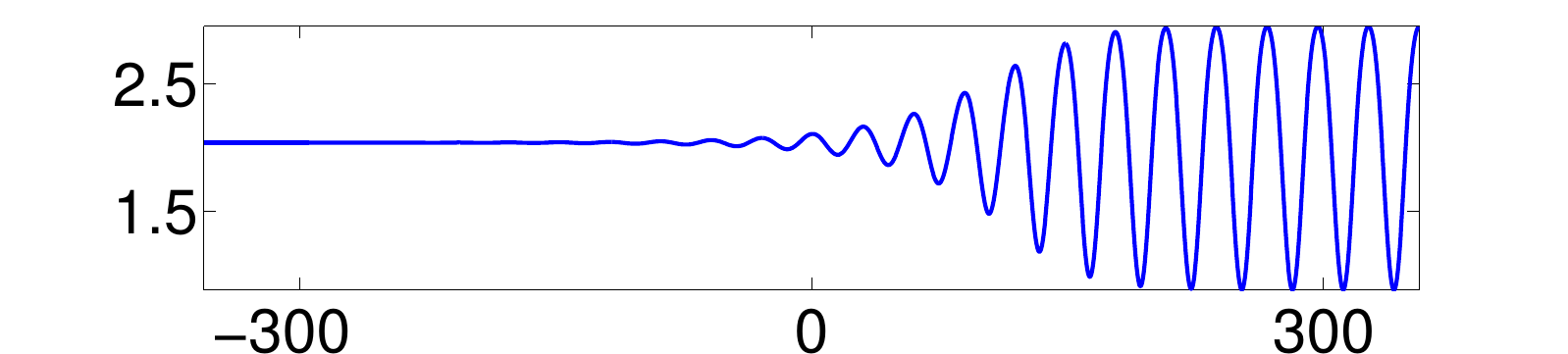}
\end{minipage} 
\caption{(a) Bifurcation diagram for 1D patterns over the domain $\Omega=(-24 \pi /k_c,24 \pi /k_c)$ for $\gamma= 0.004$ including the homogeneous solution (black), hot stripes (gray), and a branch of periodic connections between stripes and the homogeneous state (red) which bifurcates from the first bifurcation point of the hot stripes. The right and left horizontal boundaries of the zooming-in of the snake are 0.197889 and 0.197891.  (b) Plots of $u$ for solutions which are labeled in (a). The other branches which bifurcate from the hot stripe branch are illustrated \figref{1dl2} and \figref{1dl345}. Here and in the following it holds that $\norm{u}_8=(\frac{1}{\Omega}\int_\Omega |u(z)|^8 \text{d}z)^{1/8}$, where $\Omega$ is the considered domain and $z$ represents the spatial coordinates. This norm is used to obtain bifurcation diagrams with separated branches.         }
\label{1dl1}
\end{figure}

The Landau coefficient $c_3$ evaluated in $\sigma_c$ is positive for $\gamma \in (0,0.209)$ (see \figrefb{kandc3}) and thus stripes bifurcate subcritically for these active stimulations. We calculate the branch {\tt s} of hot stripes $(u_s,v_s)$ by using {\tt pde2path} for $\gamma\approx 0.004$ over the domain $\Omega=(-24 \pi /k_c,24 \pi /k_c)$ with $k_c=0.212$ (see \figref{1dl1}). {\tt s} bifurcates subcritically from the homogeneous solution $(u^*,v^*)=(u_1,v_1)$ at $\sigma_c\approx0.196$, as predicted by the Turing and Landau analysis above. A fold occurs at $\sigma\approx 0.1985$ and the stripes become stable such that there is a bistable range between $(u^*,v^*)$ and $(u_s,v_s)$.

\begin{figure}[h]
\begin{minipage}{0.49\textwidth}
(a)\\
\includegraphics[width=1\textwidth]{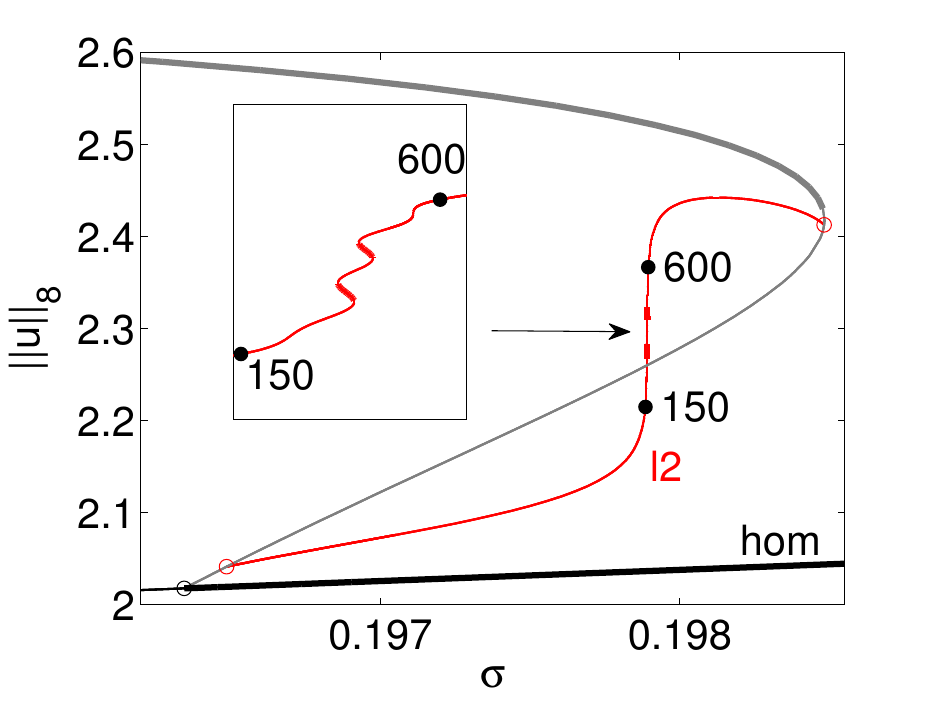}
\end{minipage} 
\begin{minipage}{0.49\textwidth}
(b)\hspace{2.5cm} \small{150}\\
\includegraphics[width=1\textwidth]{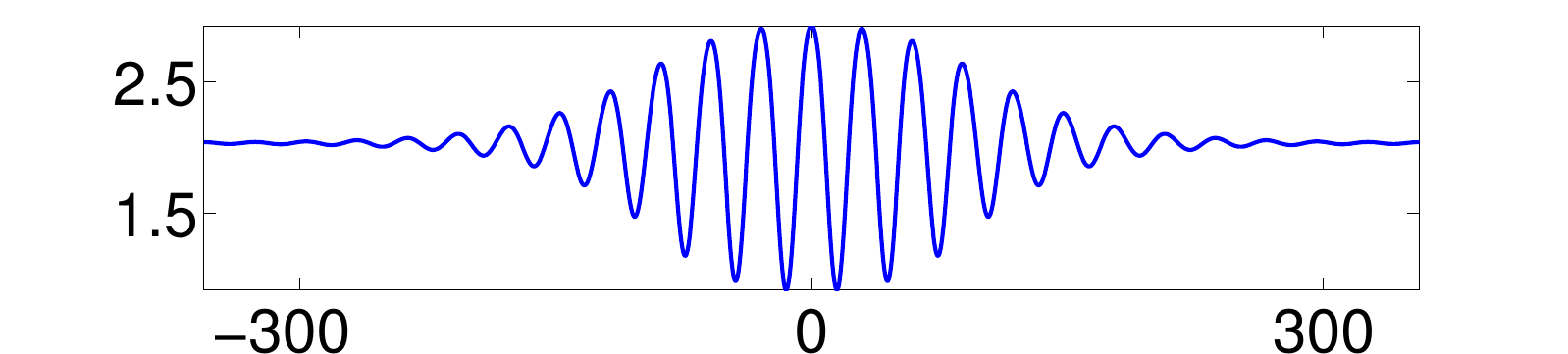}
(c)\hspace{2.5cm} \small{150} \\ \includegraphics[width=1\textwidth]{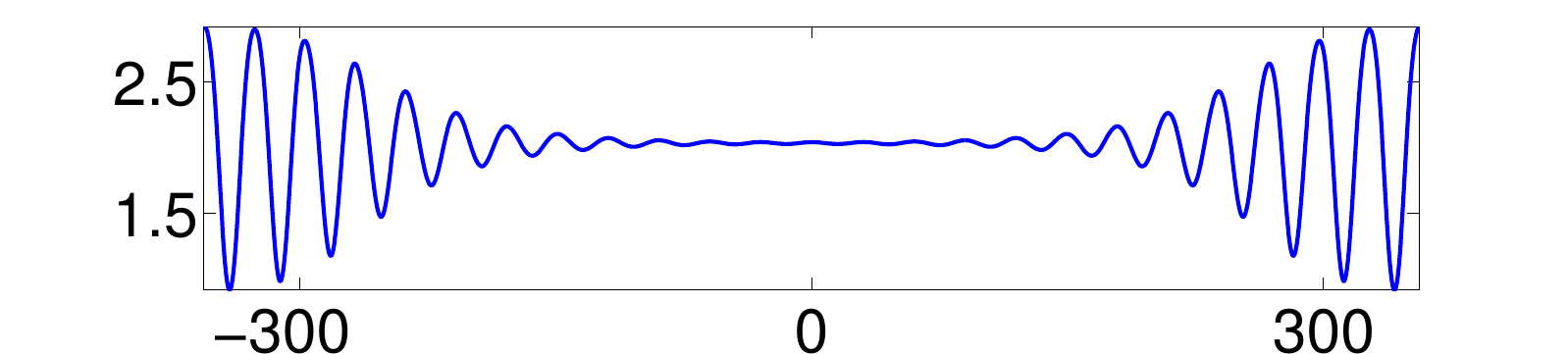}
(d)\hspace{2.5cm} \small{600} \\ \includegraphics[width=1\textwidth]{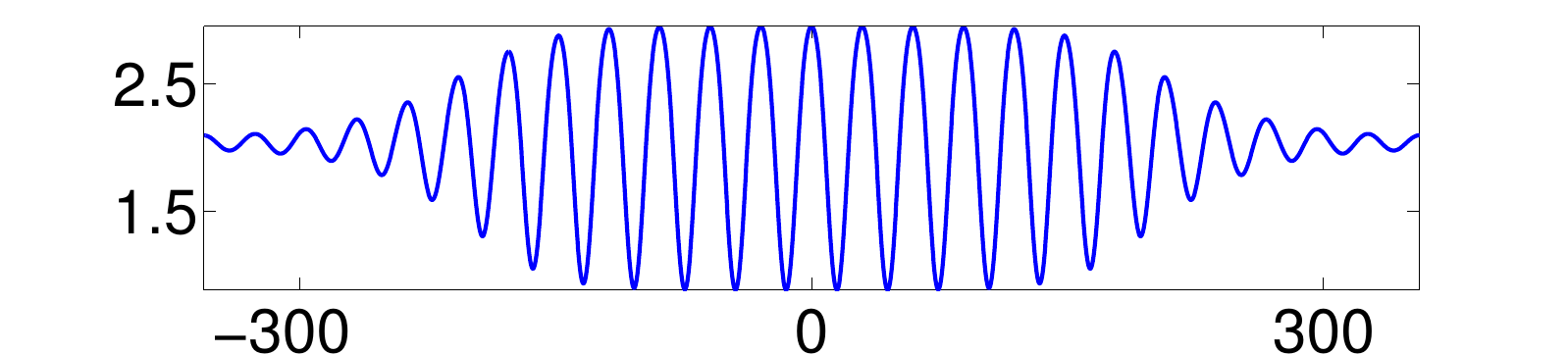}
\end{minipage} 
\caption{(a) Bifurcation diagram for the branch {\tt l2}, which bifurcates from the second bifurcation point of the hot stripes. Plots of $u$ for solutions of {\tt l2} at points 150 and 600 are shown in (b) and (d), respectively. (c) u at point 150 by using the negative tangent of the one we used to bifurcate on {\tt l2}.}
\label{1dl2}
\end{figure}

\begin{figure}[h]
\begin{minipage}{0.49\textwidth}
(a)\\
\includegraphics[width=1\textwidth]{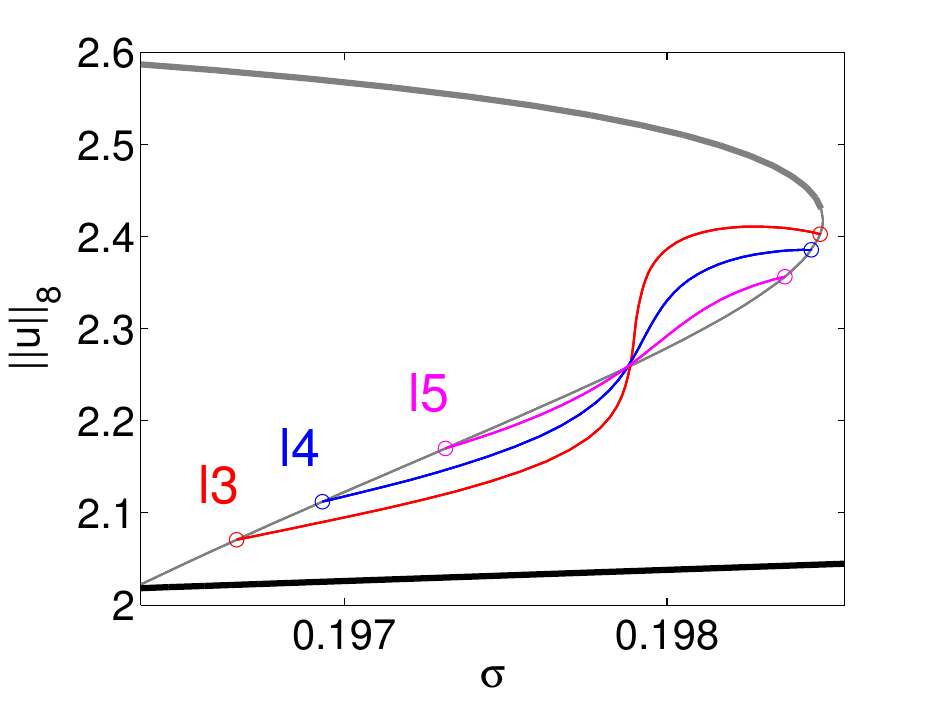}
\end{minipage} 
\begin{minipage}{0.49\textwidth}
(b)  \small{$u$ at any point in the middle of {\tt l3}, {\tt l4}, {\tt l5}}\\
\includegraphics[width=1\textwidth]{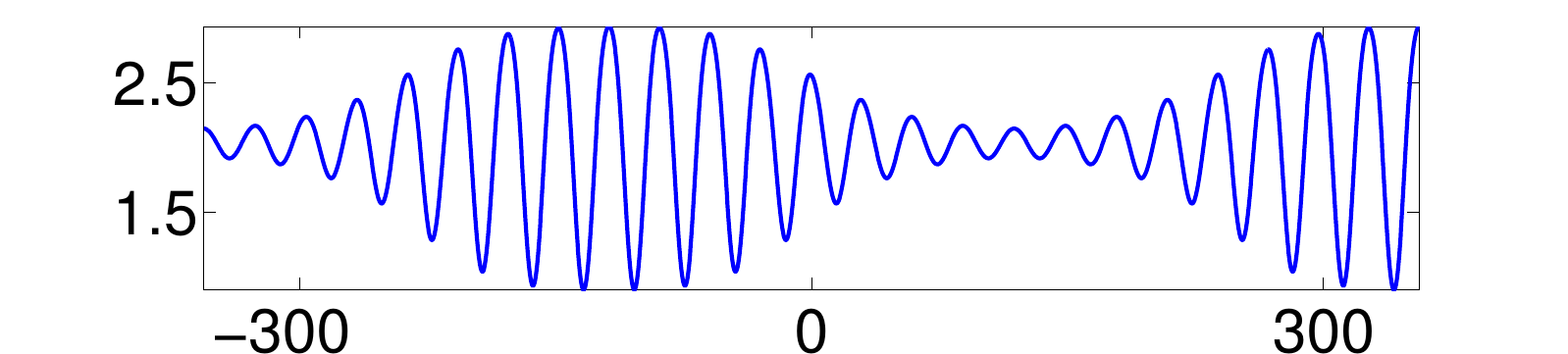}
\includegraphics[width=1\textwidth]{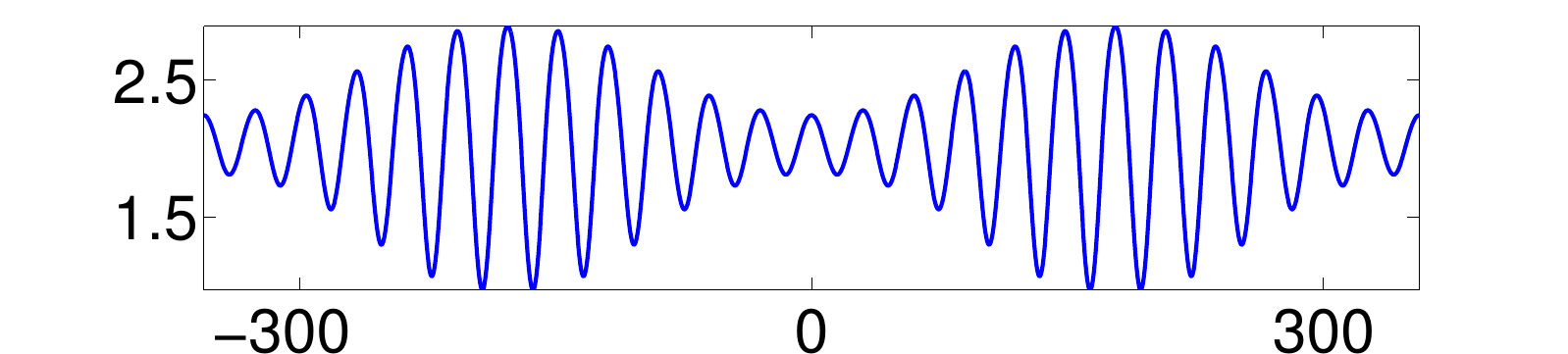}
\includegraphics[width=1\textwidth]{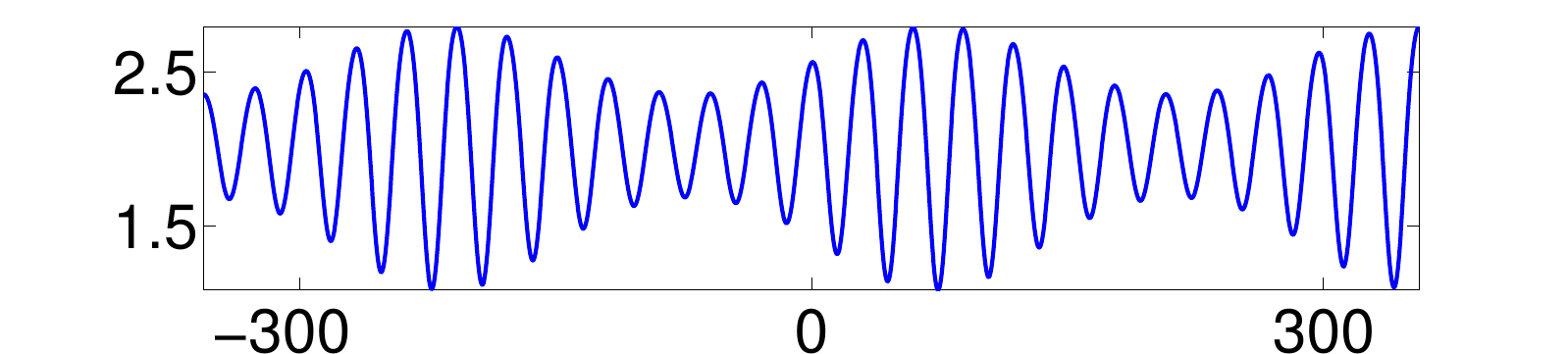}
\end{minipage} 
\caption{(a) Bifurcation diagram for the branches {\tt l3}, {\tt l4}, and {\tt l5}, which bifurcates from the third, fourth, and fifth bifurcation point of the hot stripes, respectively. (b) $u$ on  {\tt l3}, {\tt l4}, and {\tt l5}.}
\label{1dl345}
\end{figure}

There are 10 bifurcation points on {\tt s} on the way from its bifurcation to the fold. The first and 10th, second and 9th, third and 8th, 4th and 7th, 5th and 6th are connected pairwise by branches of stationary states, which we call {\tt l1}, {\tt l2},  {\tt l3}, {\tt l4}, and {\tt l5}, respectively (see \figref{1dl1}, \ref{1dl2}, and \ref{1dl345}).

Solutions on these branches are of the form
\ali{u=u^*+A \cos(k_c x)+\text{h.o.t.}, \qquad v=v^*+B \cos(k_c x)+\text{h.o.t.},\label{Neumanncos}
}
where $A$ and $B$ are space-dependent amplitudes. It holds that $A\in(0,A_s)$ and $B\in(0,B_s)$, where $A_s$ and $B_s$ are the amplitudes of the corresponding stripe solution. We can \replaced{extend}{ continue} such solutions into the right and left spatial direction periodically, since we use Neumann boundary conditions. Thus these solutions move spatially to and fro between the homogeneous and stripe solution on the entire real line. We call such solutions periodic connections. Let $n\in\{1,2,3,4,5\}$. The wavelength of $A$ and $B$ is $2|\Omega|/n$ for solutions on {\tt ln}.
For simplicity we only describe the behavior of $A$ and $u$. 

The amplitude $A$ is constant for solutions on {\tt s}. By following {\tt l1} from the first bifurcation point on {\tt s} to the first fold, which lies near the 171st solution, $A$ transforms to a nonhomogeneous state, which looks on $\Omega$ like a front between $0$ and $A_s$.
After this 
{\tt l1} 'snakes' back and forth by changing its stability. Along the snake the inflection point of the overlying function $A$ moves from the right to the left, and the position of the inflection point of $A$ shifts by $\pi/(2k_c)$ between two successive folds.
Beyond the last fold $A$ starts to grow at the left boundary and {\tt l1} returns to {\tt s}.

A \replaced{shift symmetry}{ reflection} also exists for every solution of {\tt ln} (see \added{\figref{1dl345}(b) and (e))}. We call the corresponding branch {\tt ln'}. Clearly, the illustration of {\tt ln'} is congruent with {\tt ln} and both branches together generate a loop from one bifurcation point to the other and back again.

Splitting the domain $\Omega$ into $\Omega_l=(-24 \pi /k_c,0)$ and $\Omega_r=(0,24 \pi /k_c)$, we can describe solutions on {\tt l2} as two front-like connections between the homogeneous and the striped solution on $\Omega_l$ and $\Omega_r$, where the one on $\Omega_r$ is a reflection of the one on $\Omega_l$. Here the snake is shorter with fewer wiggles, because the lengths of the domains $\Omega_l$ and $\Omega_r$  are shorter than the length of $\Omega$ such that the inflection point of $A$ reaches the boundaries of $\Omega_l$ and $\Omega_r$ earlier.  
We can split the domain $\Omega$ into 3, 4, and 5 parts and use front-like connections as above to describe solutions of {\tt l3}, {\tt l4}, and {\tt l5}, respectively. Here the branches do not show any snaking behavior, because the partitions of $\Omega$ are too small. 
The same branches for connections between cold stripes and the homogeneous solutions can be found on the cold-stripe branch.

For more details on localized patterns and snaking over bounded domains see \cite{bbkm2008,dawes08,dawes09,KAC09,hokno2009}.
Seminal results for localized patterns over unbounded domains can be found in \cite{burke,bukno2007,BKLS09}. For a detailed analysis by using the Ginzburg-Landau formalism and beyond all order asymptotics see \cite{chapk09,dean11}.

With respect to the results of for instance \cite{BKLS09} we expect that additional localized stripes of the form
\ali{u=u^*+A \sin(k_c x)+\text{h.o.t.} \label{Dirichletsin}}
exist for periodic boundary conditions. Let us call the corresponding branches {\tt d1}, {\tt d2}, {\tt d3}, {\tt d4}, and {\tt d5}.  The illustrations of these branches are not congruent with {\tt l1}, {\tt l2}, {\tt l3}, {\tt l4}, and {\tt l5}. The folds of {\tt l1} and {\tt d1},..., {\tt l5} and {\tt d5} are connected pairwise by branches, which are called rungs. 

Furthermore, we can expect from these numerical results that the following solutions exist over the entire real line: a heteroclinic connection between stripes and the homogeneous solution, a homoclinic connection from the homogeneous solution to stripes and back to the homogeneous state, and a homoclinic connection from stripes to the homogeneous solution, and back to stripes.

\section{2D Patterns}\label{2d}
\begin{figure}[h]

\begin{minipage}{0.49\textwidth}
\includegraphics[width=1\textwidth]{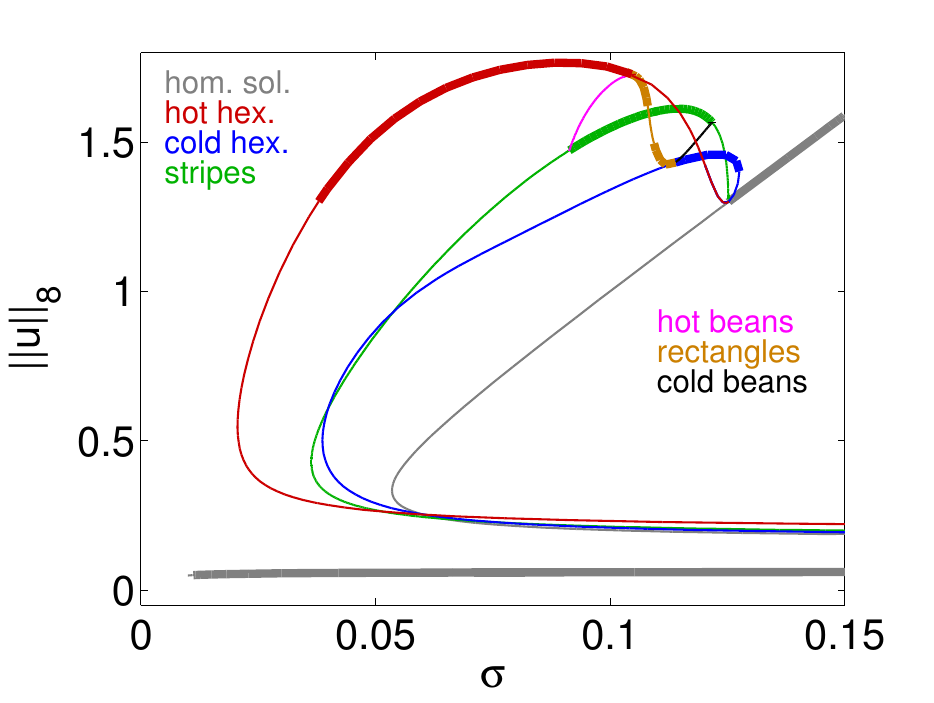}
\end{minipage}
\begin{minipage}{0.24\textwidth}
\tiny{stable cold hex. $\sigma=0.127$  } \\
\includegraphics[width=1\textwidth]{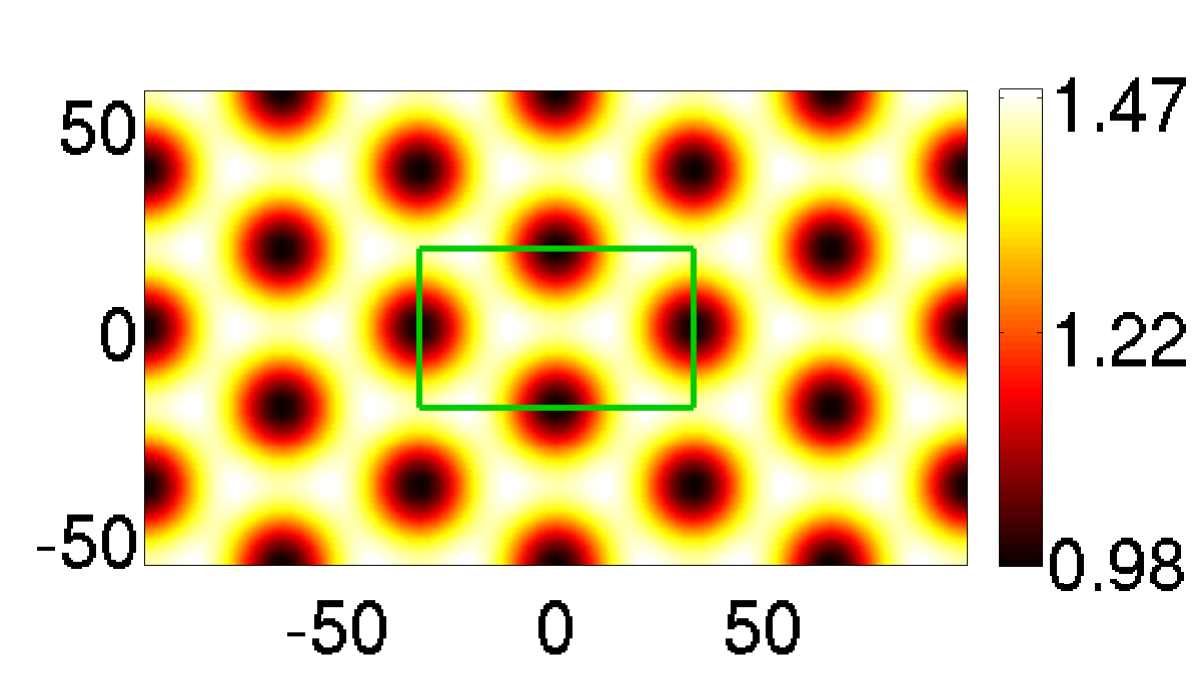}\\ \\
\tiny{unst. cold hex. $\sigma=0.04$}\\
\includegraphics[width=1\textwidth]{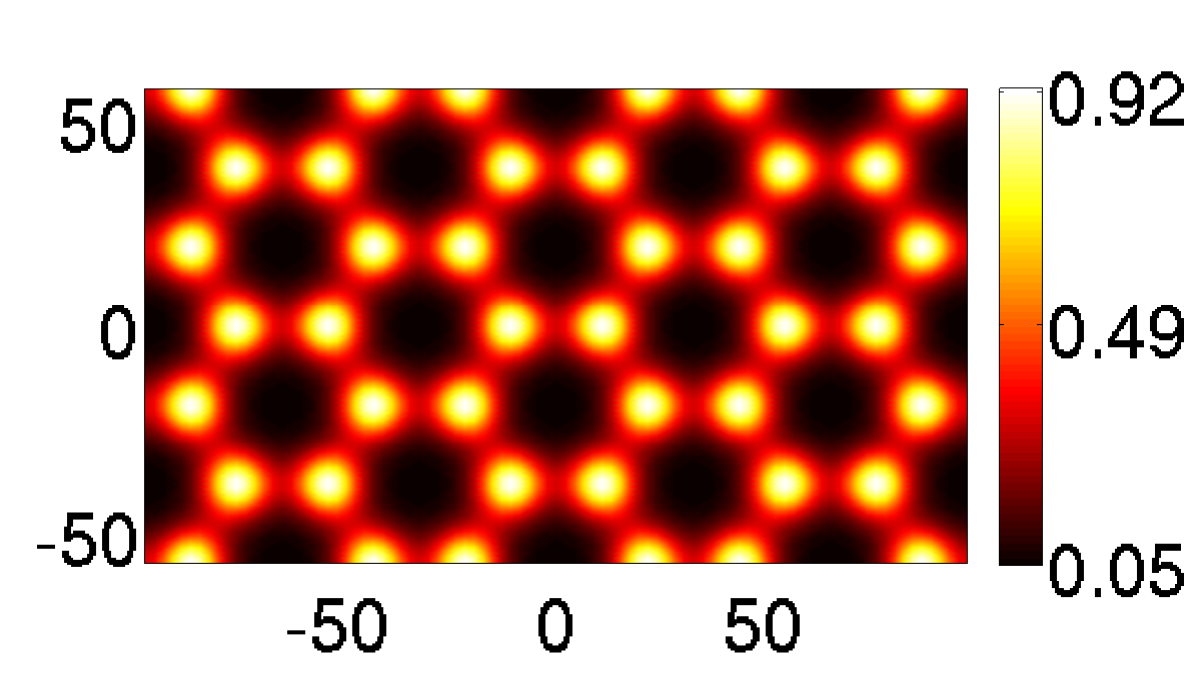}
\end{minipage}
\begin{minipage}{0.24\textwidth}
\tiny{stable cold hex. $\sigma=0.114$}\\
\includegraphics[width=1\textwidth]{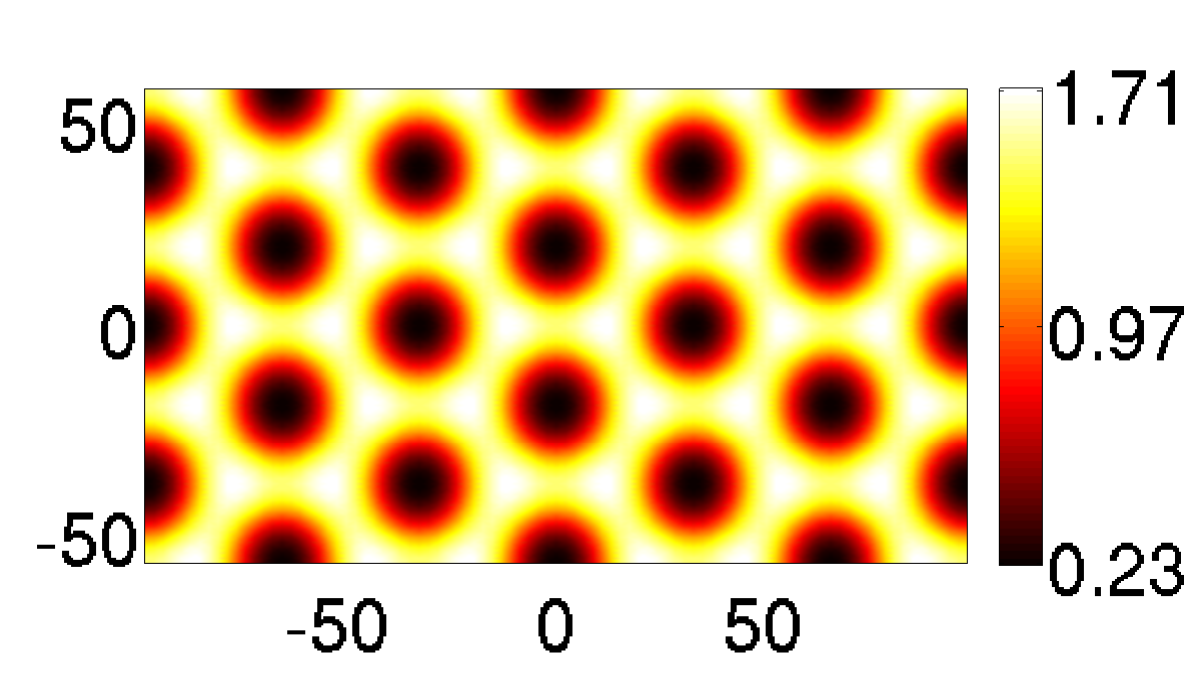}\\ \\
\tiny{cold bean $\sigma=0.12$}\\
\includegraphics[width=1\textwidth]{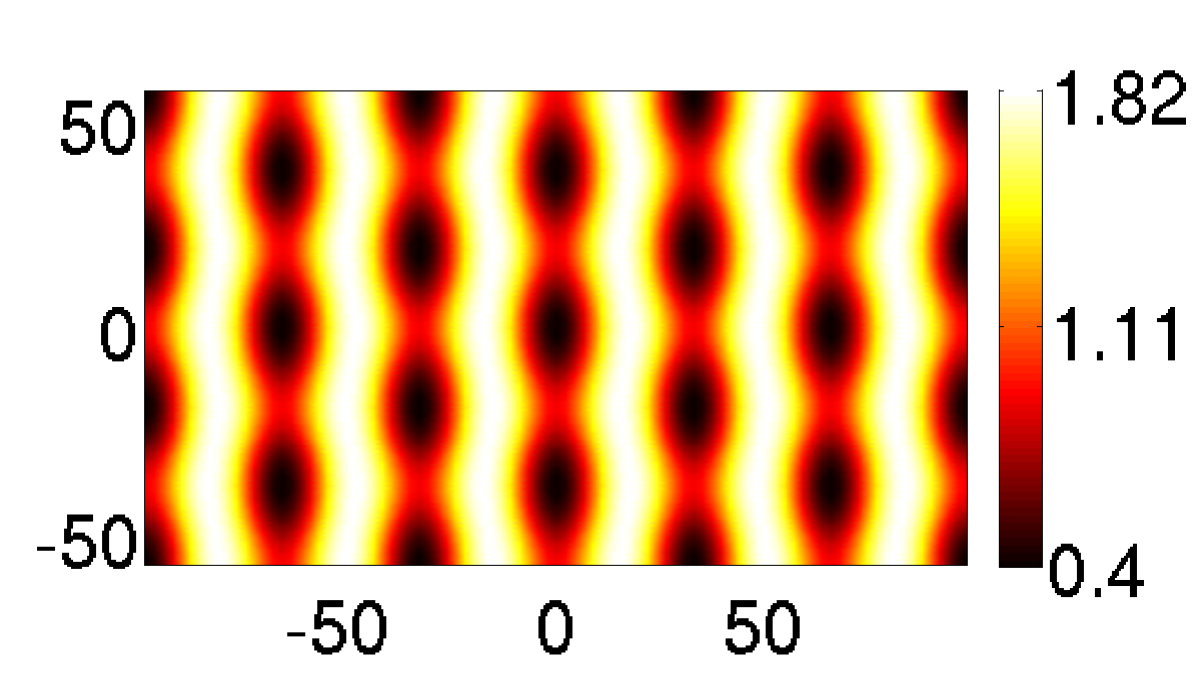}
\end{minipage}

\begin{minipage}{0.24\textwidth}
\tiny{stable stripes $\sigma=0.12$}\\
\includegraphics[width=1\textwidth]{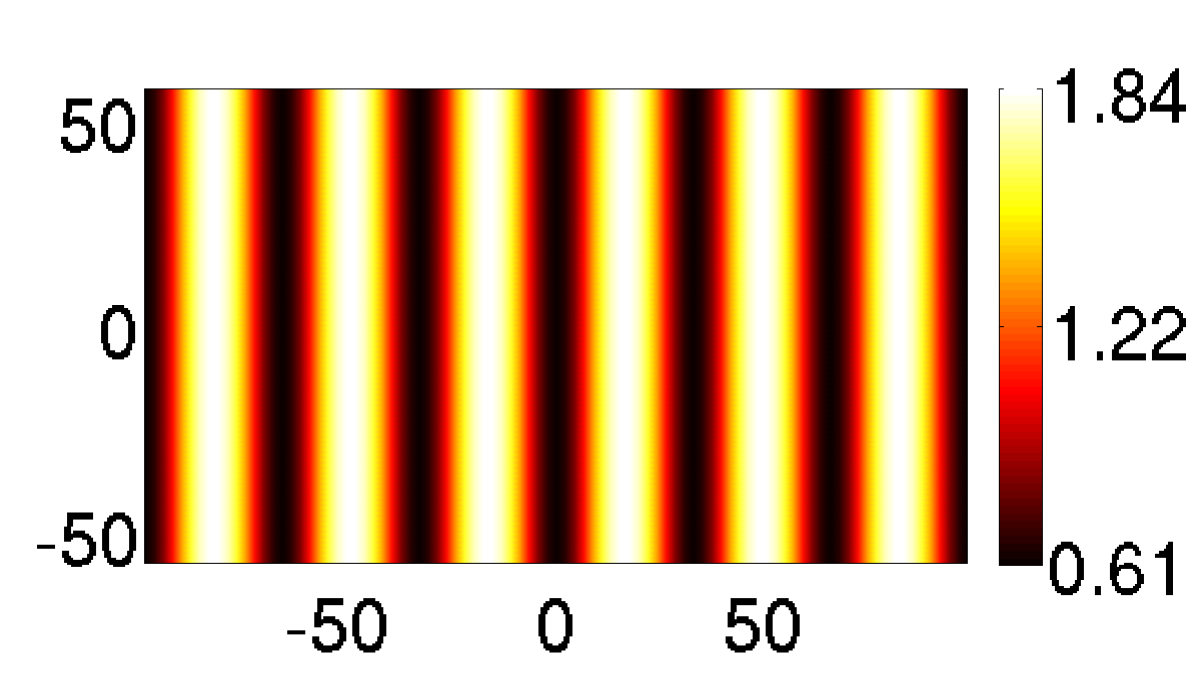}\\ \\
\end{minipage}
\begin{minipage}{0.24\textwidth}
\tiny{stable stripes $\sigma=0.09$}\\
\includegraphics[width=1\textwidth]{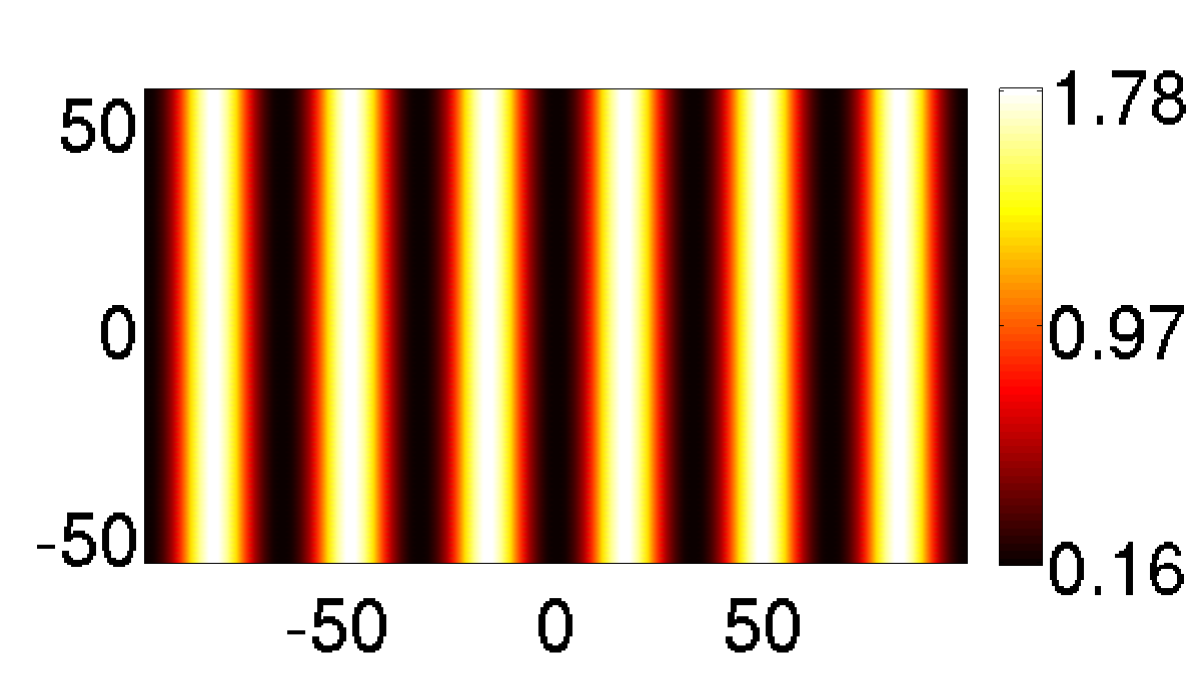}\\ \\
\end{minipage}
\begin{minipage}{0.24\textwidth}
\tiny{stable hot hex. $\sigma=0.1$}\\
\includegraphics[width=1\textwidth]{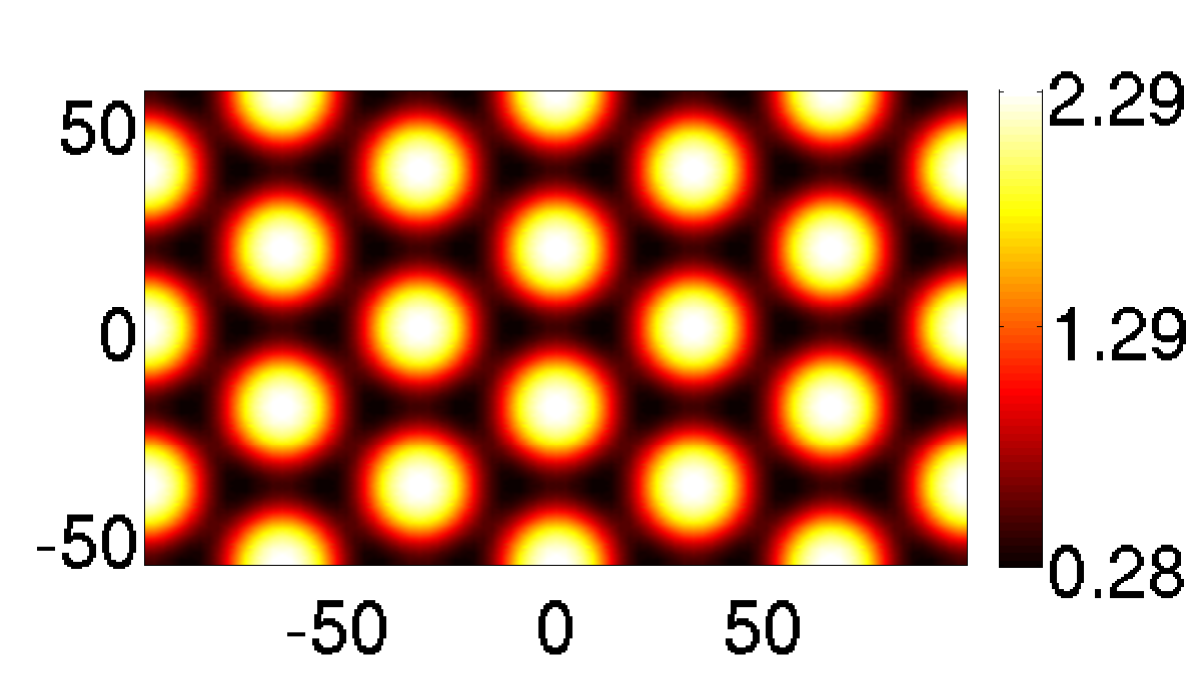}\\ \\
\end{minipage}
\begin{minipage}{0.24\textwidth}
\tiny{stable hot hex $\sigma=0.04$}\\
\includegraphics[width=1\textwidth]{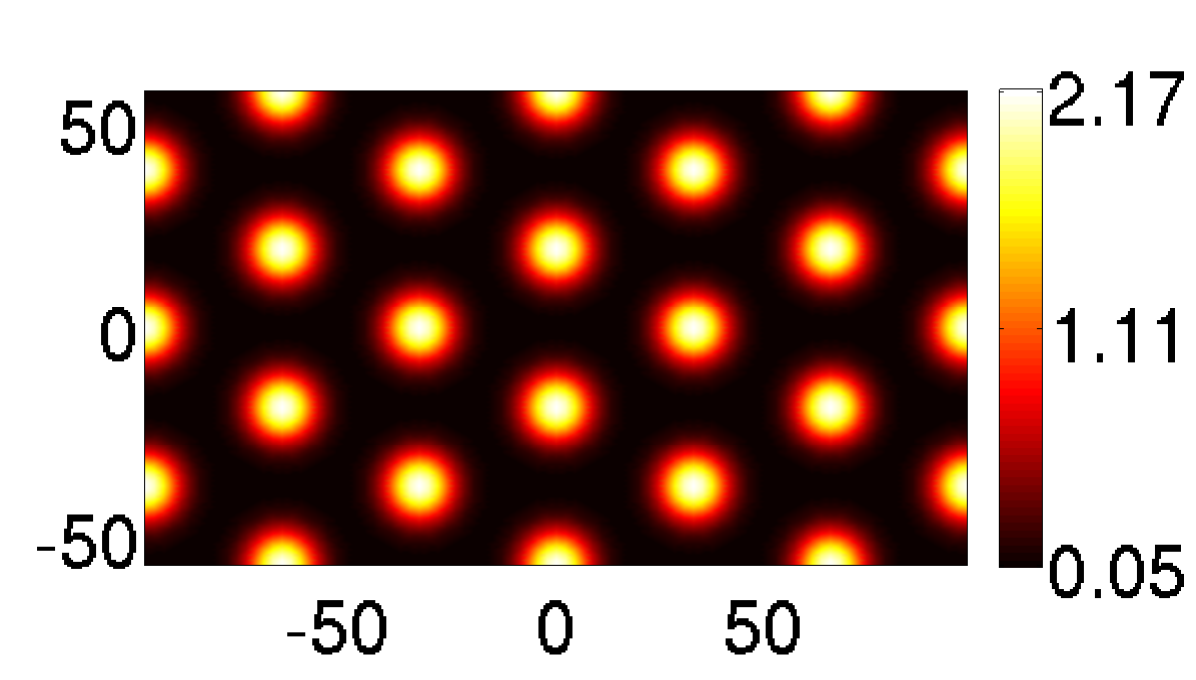}\\ \\
\end{minipage}

\begin{minipage}{0.24\textwidth}
\tiny{stable cold rec. $\sigma=0.11$}\\
\includegraphics[width=1\textwidth]{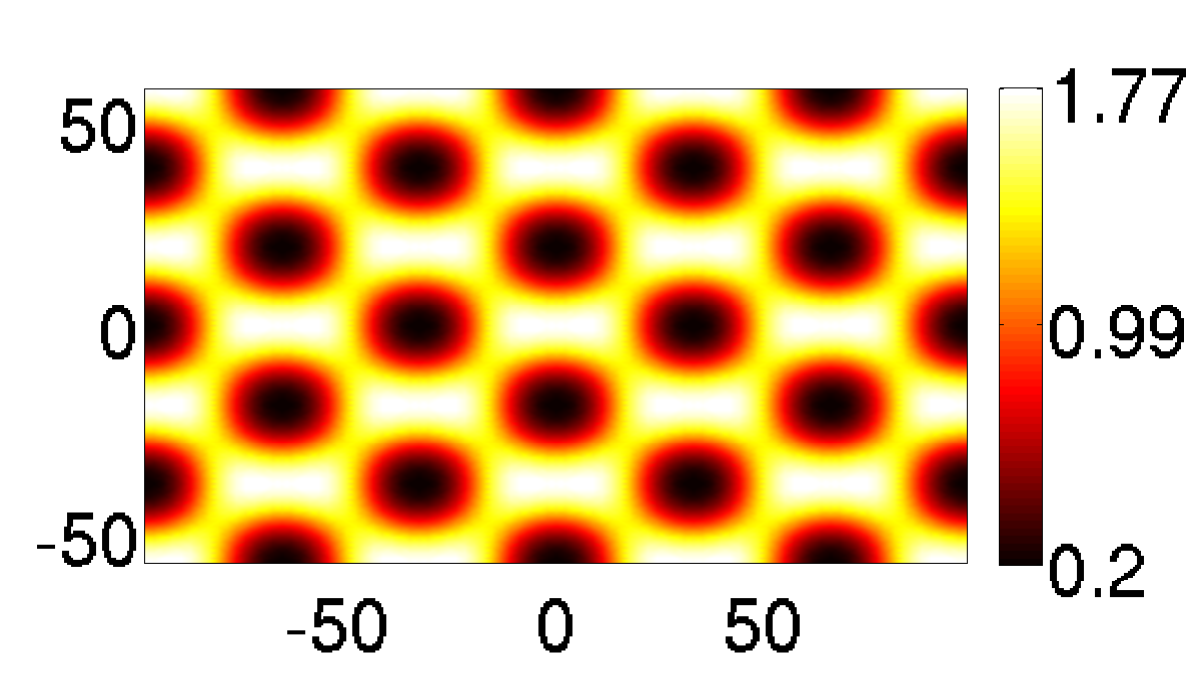}\\ \\
\end{minipage}
\begin{minipage}{0.24\textwidth}
\tiny{unstable deg. rec. $\sigma=0.108$}\\
\includegraphics[width=1\textwidth]{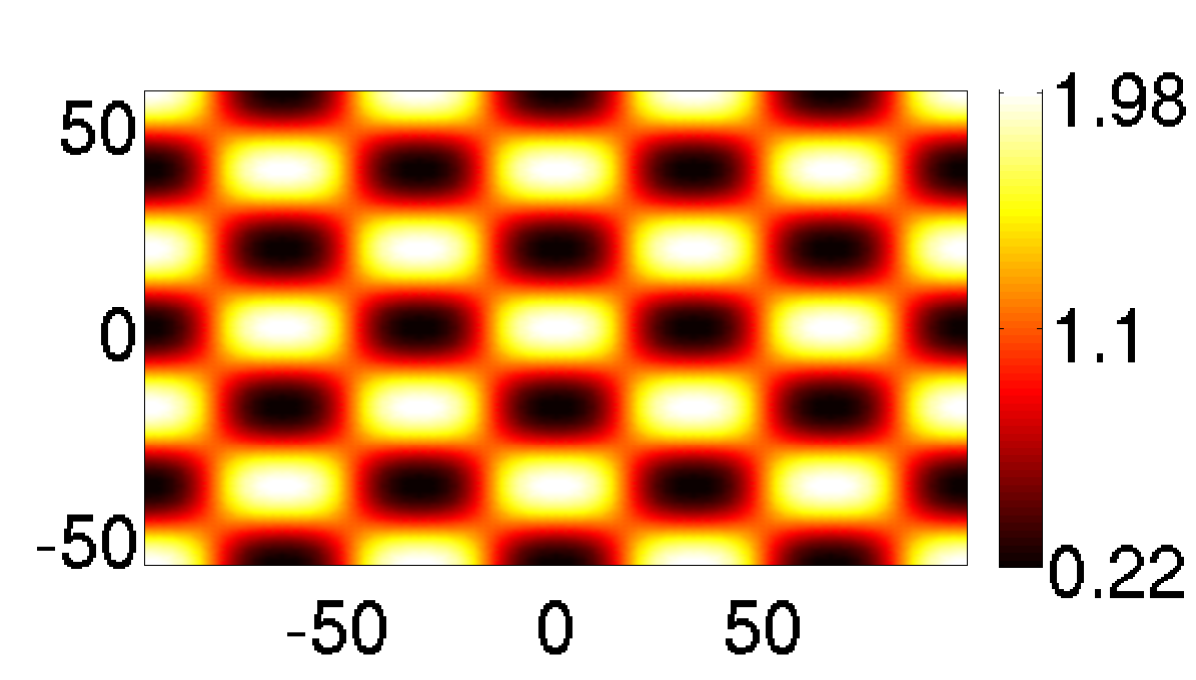}\\ \\
\end{minipage}
\begin{minipage}{0.24\textwidth}
\tiny{stable hot rec. $\sigma=0.107$}\\
\includegraphics[width=1\textwidth]{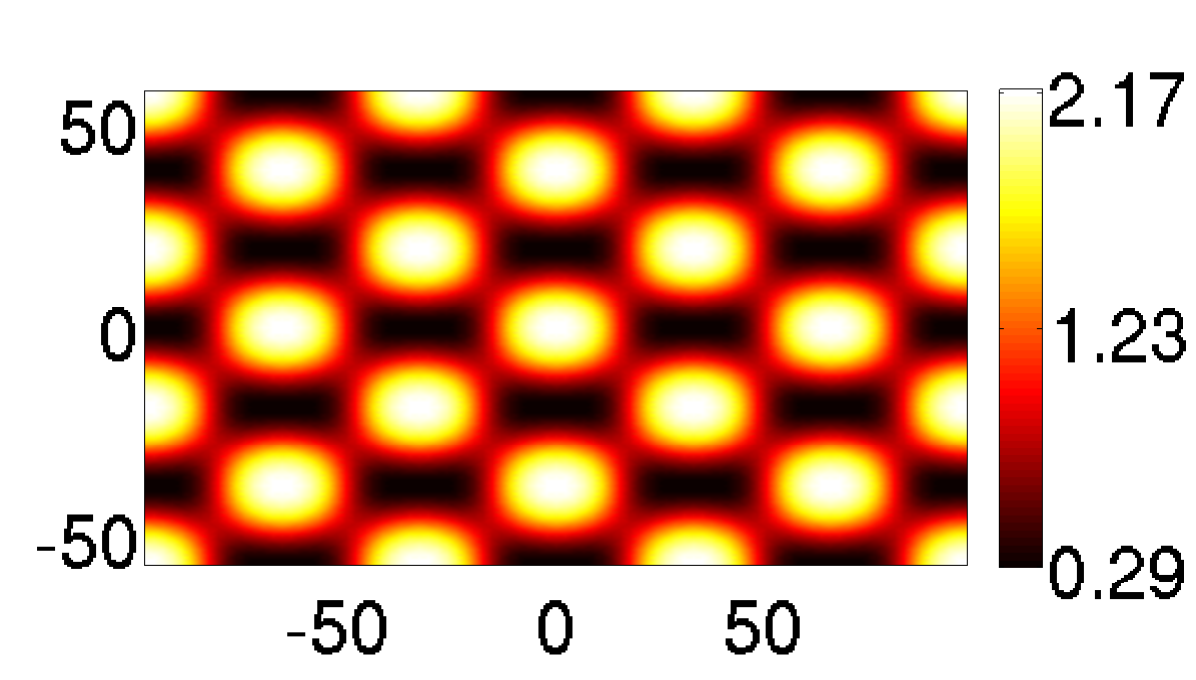}\\ \\
\end{minipage}
\begin{minipage}{0.24\textwidth}
\tiny{hot beans $\sigma=0.1$}\\
\includegraphics[width=1\textwidth]{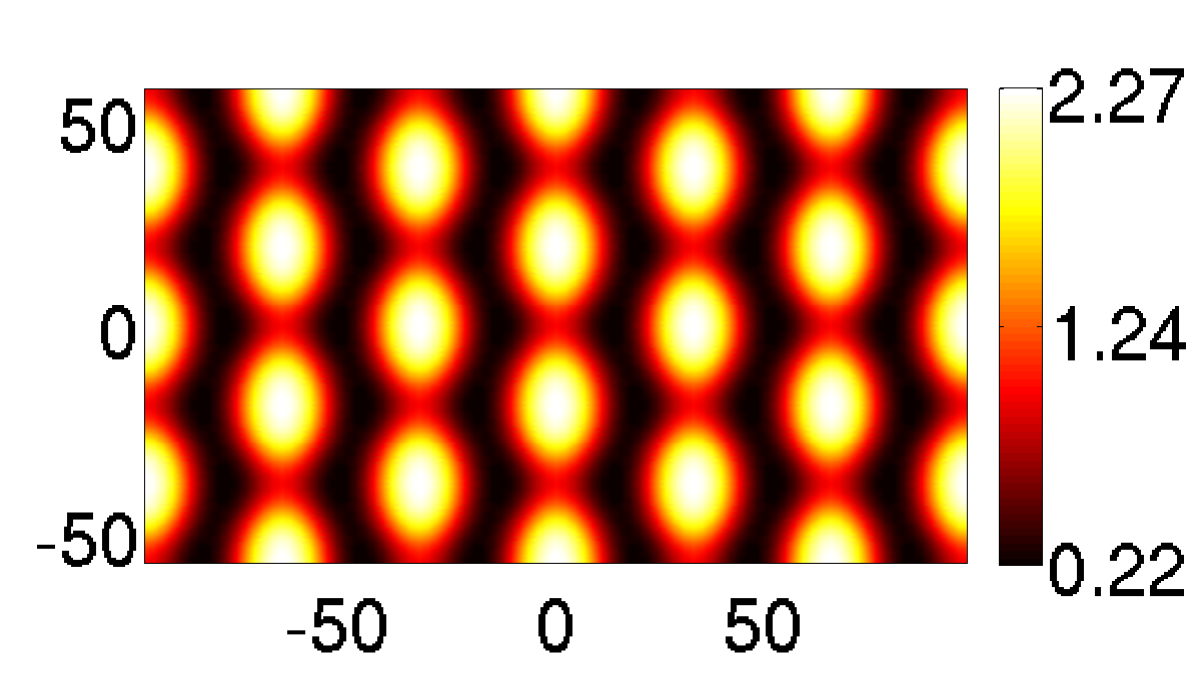}\\ \\
\end{minipage}
\caption{Shown is a bifurcation diagram and density plots of $u$ for some special solutions. Here $\gamma=0.25$. For a better illustration we plot the solutions on a large domain, but the domain, which is used for the numerics, is $\Omega=(-l_x,l_x) \times(-l_y,l_y)$, where $l_x=2\pi/kc$,  $l_y= 2\pi/(\sqrt{3} k_c)$, and $k_c\approx 0.19$. The domain $\Omega$ is marked in the density plot of the cold hexagons for $\sigma=0.127$ (green box).
} \label{hsbbifu}
\end{figure}

The stripe patterns also exist over two dimensional domains. Here we present additionally some genuine 2D patterns, which can be analyzed via the Landau system \eqref{gl}. Setting $A_1=A_2=A_3=:A$, the system \eqref{gl} reduces to
\alinon{\partial_tA=c_1A+c_2\overline{A}^2+(c_3+2c_4)A|A|^2.}
Stationary amplitudes fulfill
\ali{c_1A+c_2\overline{A}^2+(c_3+2c_4)A|A|^2=0.\label{hexpolynom}}
If
\ali{A= H_\pm:=-\frac{c_2}{2(c_3+2c_4)}\pm \sqrt{\frac{c_2^2}{4(c_3+2c_4)^2}-\frac{\mu(k_c)}{c_3+2c_4}} \label{hexamp}}
are real, then \eqref{hexamp} solve\added{s} \eqref{hexpolynom}. These amplitudes generate hexagon patterns. All other solutions of \eqref{hexpolynom} generate phase shifts of these hexagon patterns.
We already mentioned that we classify the hexagons in hot and cold, which means that the hexagon pattern, which corresponds to the bacteria, has maximums and minimums in the center of the hexagonal spots, respectively.  
Inserting $(A_1,A_2,A_3)=(A,B,B)$ into \eqref{gl} one can find mixed mode solutions. A solution, for which $\abs{A}\neq \abs{B}$ and $A\neq0$ holds, is called rectangle in \cite{goswkn84,Hoyle}. We classify them into bean and rectangle patterns, which fulfill $|A|>|B|$ and $|A|<|B|$, respectively. We choose these names because of their 2D density plots (see \figref{hsbbifu}) and their different roles in the bifurcation diagram. These roles will be described below. We classify beans and rectangles into hot and cold in a similar way as for hexagons. \\
There are a lot of studies, which show bifurcation diagrams over 2D domains, which are \replaced{generated}{ developed} by using amplitude equations (see e.g. \cite{goswkn84,pomeau,MNT90,dsss03,Hoyle}). This method works only well near the onset. To obtain good results further away from the onset, one can use numerical path following methods. This is done in a few works for problems over 2D domains (see e.g. \cite{hexsnake,strsnake,p2p,schnaki}).
We used the numerical bifurcation and continuation software {\tt pde2path} to \replaced{generate}{ develop} \figref{hsbbifu}. Shown is a bifurcation diagram, which contains branches of homogeneous solutions, hexagons, stripes, rectangles, and beans. We see that stripes and hexagons bifurcate from the homogeneous solution, when the homogeneous solution becomes unstable. Their branches have stable parts. Decreasing $\sigma$ from 0.15 to 0 we see that cold hexagons are stable before stripes and stripes before hot hexagons. Thus a change  from the homogeneous solution to hot hexagons by passing cold hexagons and stripes is expected if the balancing rate $\sigma$ decreases slowly. We conclude that such a change of patterns can arise if the food influx decreases and thus it can be seen as a possible signal that the bacteria is in danger to die out. This is an often observed order of patterns (see e.g. \cite{meron01,meron05,schnaki}). It is shown on a Landau level on a hexagonal lattice that this is in general a robust pattern sequence in \cite{TransPattern}, but using the Landau formalism means that this holds only near the onset for small amplitudes so that it is unclear if this robustness also holds in general.  \\ 
For stripes and both types of hexagons we show density plots for bacteria near the endpoints of their stable ranges (see \figref{hsbbifu}). We see that the hot part of the solution becomes smaller for all three pattern types if $\sigma$ decreases. This means for cold and hot hexagons that the size of hexagons become larger and smaller, respectively. So \replaced{another}{ an other} indicator for the decrease of the balancing rate (and with that the food influx) is a change of the pattern size itself. Such ideas are already pointed out in \cite{meron05} for a vegetation model.  
\\
There are bistable ranges between the homogeneous state and cold hexagons, cold hexagons and stripes, and stripes and hot hexagons. Branches of hot and cold beans bifurcate at the left and right endpoints of the stable range of stripes, which form a connection to hot and cold hexagons, respectively. All solutions of both bean branches are unstable, which also holds by using the Landau reduction. The rectangles build a connection between cold and hot hexagons in parameter space. It is claimed in \replaced{\cite[p.153]{Hoyle}}{\cite{Hoyle} on page 153} and in \cite{TransPattern} (see table I) that the Landau reduction on hexagonal lattice predicts that the rectangles are always unstable. This is in contrast to the results in \cite{schnaki}. There it is found that rectangles are always stable on hexagonal lattices. For the system, which is considered in the \replaced{present}{ actual} work, the Landau reduction on a hexagonal lattice predicts again that all solutions on the rectangle branch are stable. We see in \figref{hsbbifu} that this does not hold, when we use \replaced{the finite element method}{ numerical methods}. We see that the rectangles are stable at the outer ranges, but unstable in the middle. By increasing the domain size we find that this unstable range become greater. If the rectangle branch is completely unstable or if a stable part remains on an unbounded domain, remains open in the \replaced{present}{ actual} work.   \\
\begin{figure}[H]
\begin{minipage}{0.24\textwidth}
\includegraphics[width=1\textwidth]{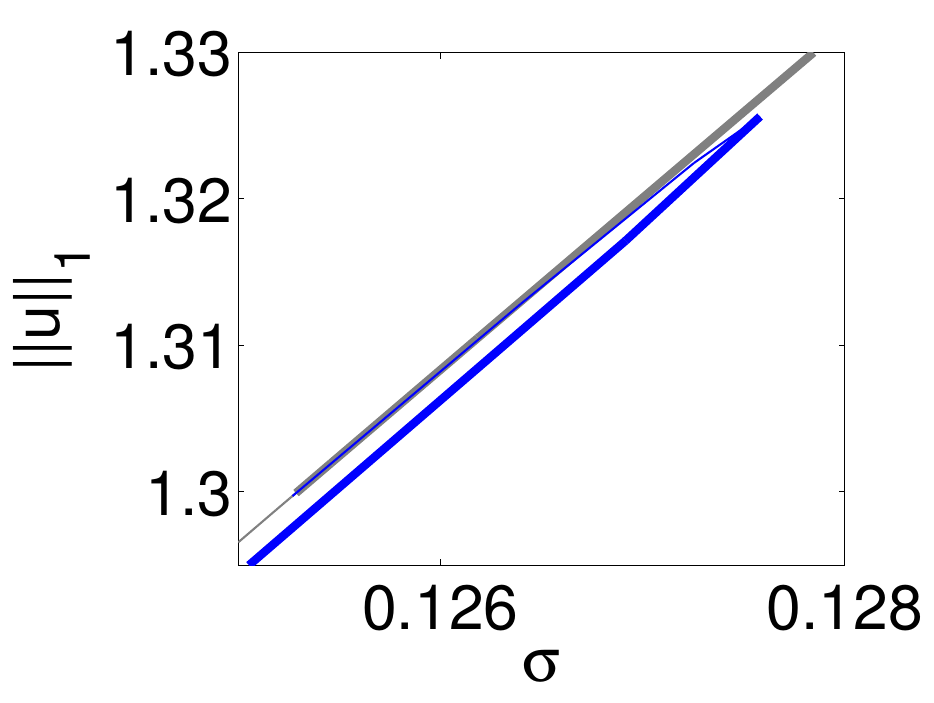}\\
\includegraphics[width=1\textwidth]{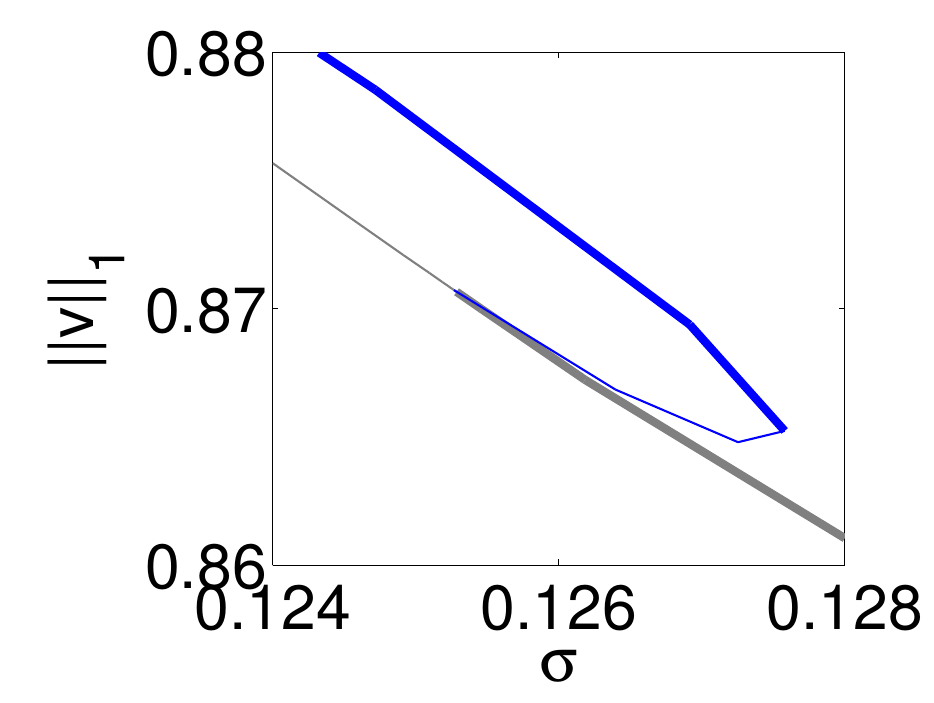}
\end{minipage}
\begin{minipage}{0.24\textwidth}
\includegraphics[width=1\textwidth]{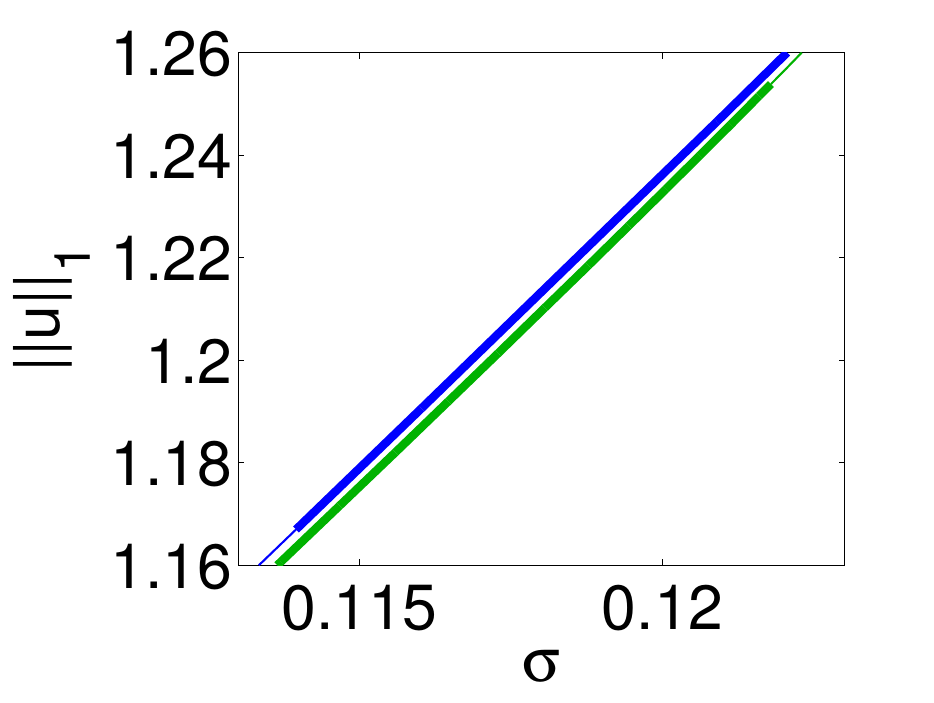}\\
\includegraphics[width=1\textwidth]{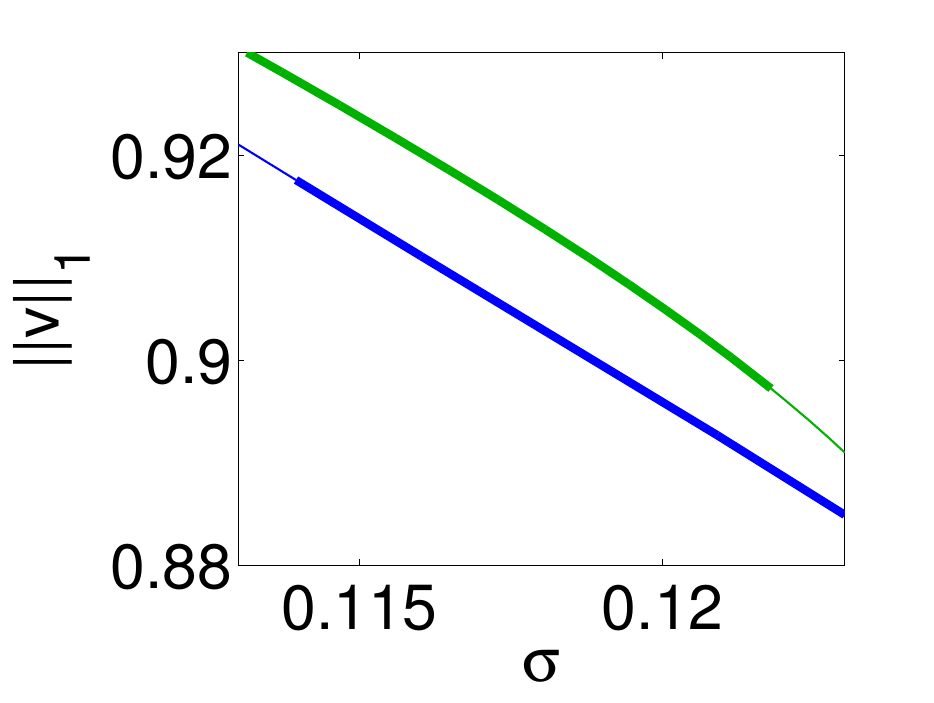}
\end{minipage}
\begin{minipage}{0.24\textwidth}
\includegraphics[width=1\textwidth]{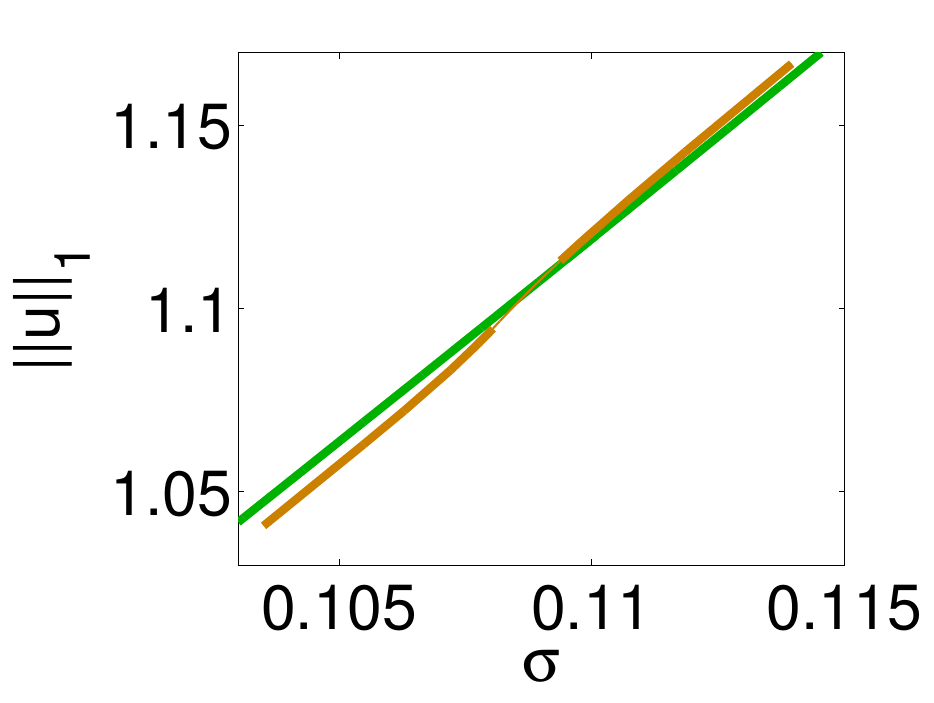}\\
\includegraphics[width=1\textwidth]{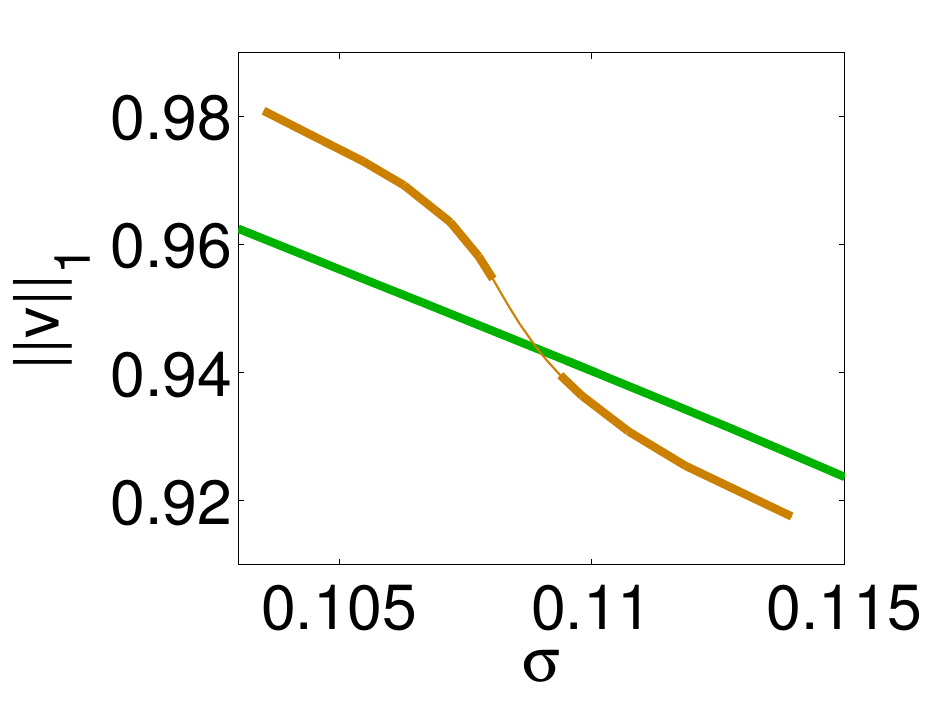}
\end{minipage}
\begin{minipage}{0.24\textwidth}
\includegraphics[width=1\textwidth]{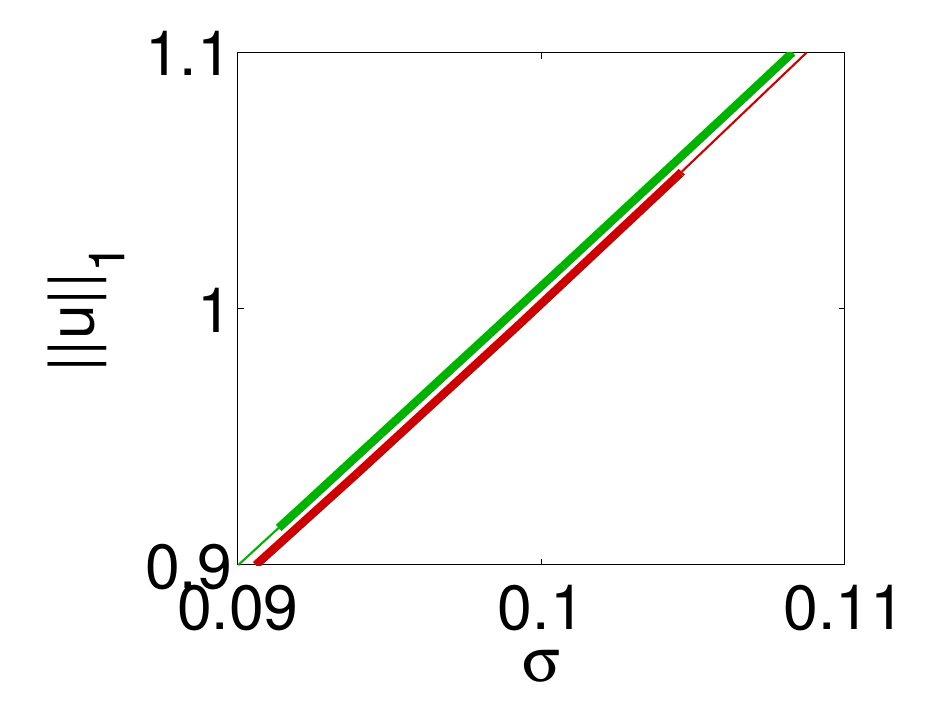}\\
\includegraphics[width=1\textwidth]{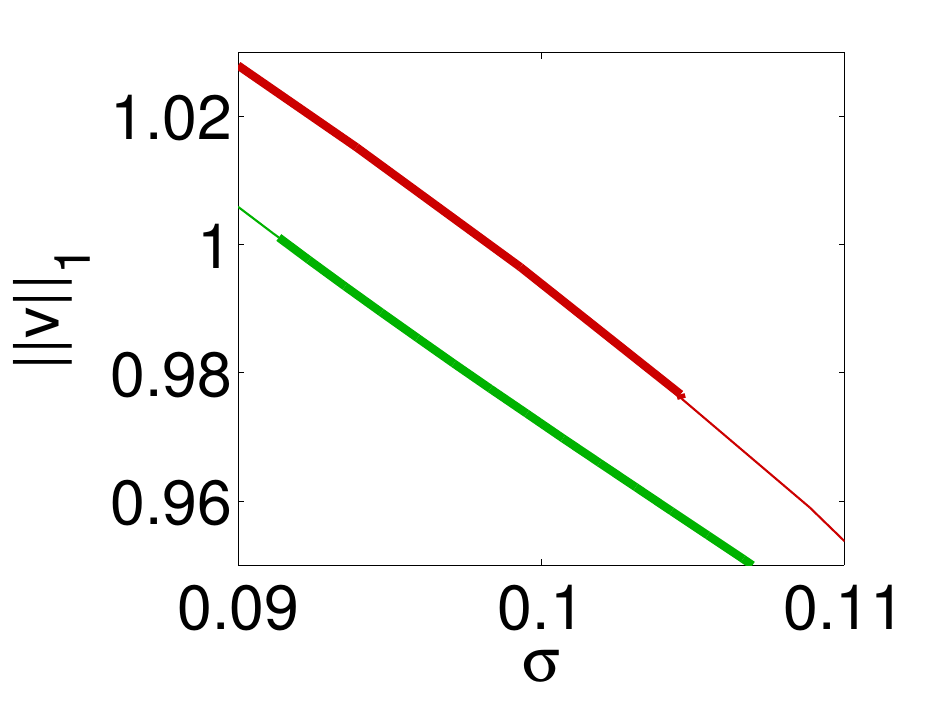}
\end{minipage}
\caption{Shown are the bistable ranges of the bifurcation diagram of \figref{hsbbifu}. Here we denote $\norm{u}_1 $ and $\norm{v}_1$ on the vertical axis instead of $\norm{u}_8$. It holds $\norm{u}_1=\frac{1}{\Omega}\int_\Omega |u(z)| \text{d}z$, where $\Omega$ is the considered domain and $z$ represents the spatial coordinates. Notice that $u(z)$ is alway nonnegative so that $\norm{u}_1$ equals the average of $u$. We use the same colors for the branches \replaced{as}{ like} in \figref{hsbbifu}. 
}\label{L1fig}
\end{figure}
We have seen that two different stable solutions exist in some $\sigma$-ranges. To see which of these is richer in bacteria or nutrient, we show the bistable ranges with respect to the normalized $L_1$-norm in \figref{L1fig}. We have also seen that by decreasing $\sigma$ the different solution types become stable in the following order: Homogeneous solution, cold hexagons, stripes, cold rectangles, hot rectangles, hot hexagons. In \figref{L1fig} we see that the solution type, which becomes stable, has fewer bacteria and more nutrient in contrast to the solution, which was already stable. Similar changes of the average for species by changing the type of solution in bistable ranges can also be seen in \cite{meron01,meron04}. 

\subsection{Localized Patterns and Snaking}
\begin{figure}[h]
\begin{minipage}{0.49\textwidth}
a)\\
\includegraphics[width=\textwidth]{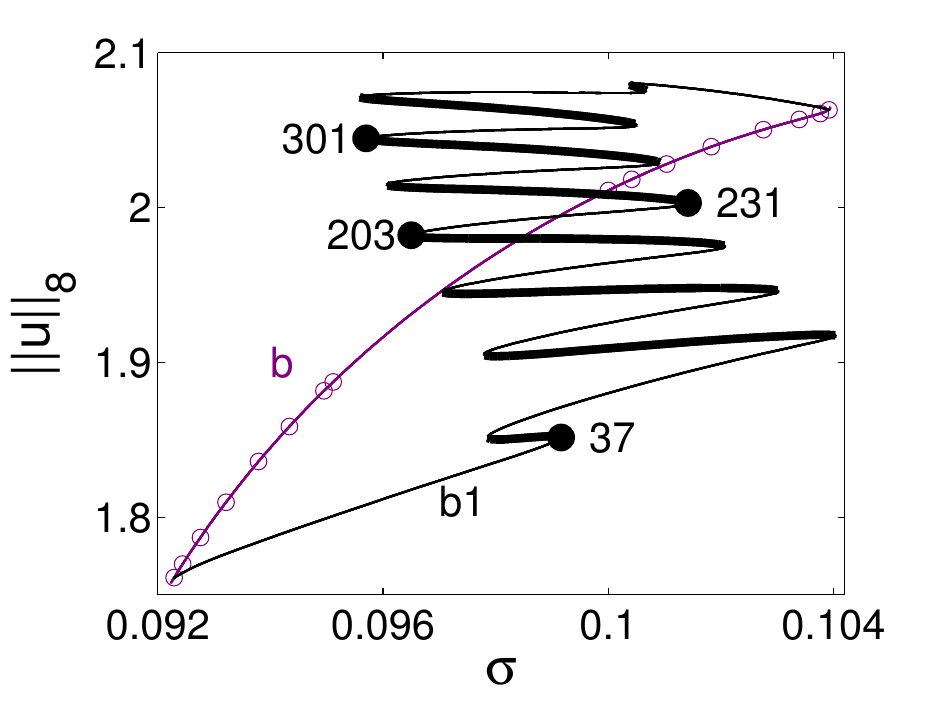}
\end{minipage} 
\begin{minipage}{0.49\textwidth}
b)\\
\includegraphics[width=\textwidth]{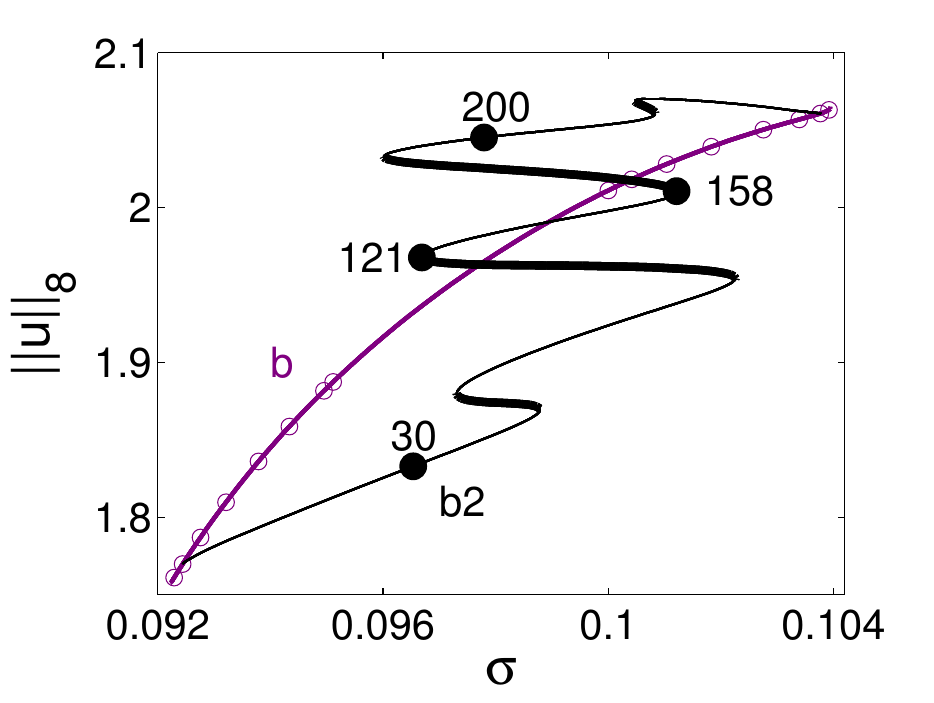}
\end{minipage} 
\begin{minipage}{0.49\textwidth}
c)\\
\includegraphics[width=\textwidth]{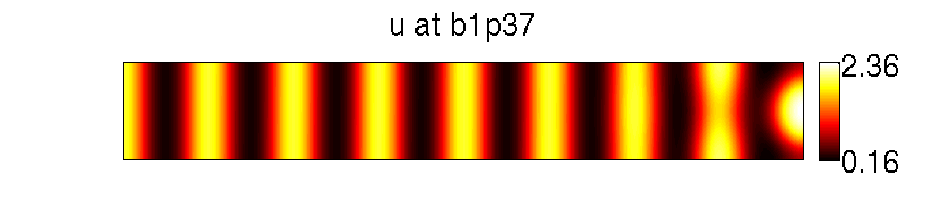}
\includegraphics[width=\textwidth]{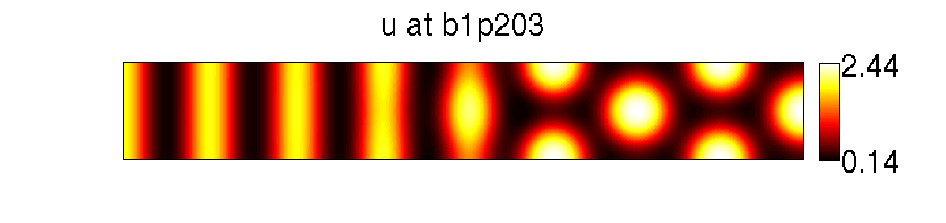}
\includegraphics[width=\textwidth]{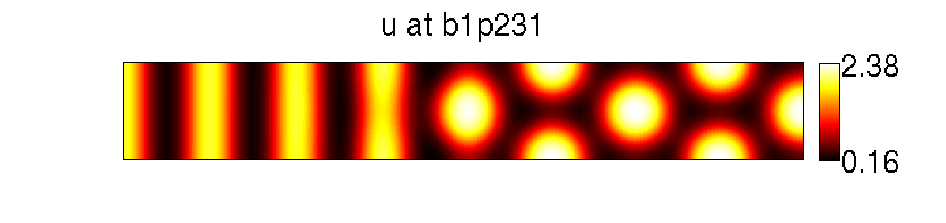}
\includegraphics[width=\textwidth]{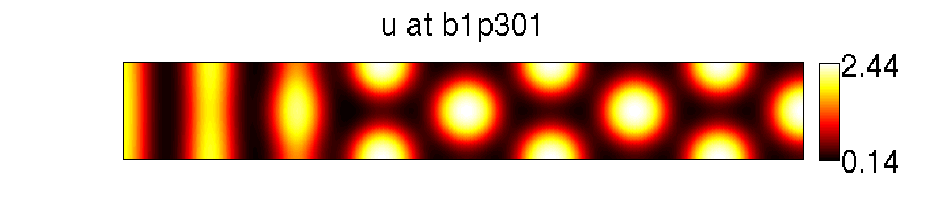}
\end{minipage}  
\begin{minipage}{0.49\textwidth}
d)\\
\includegraphics[width=\textwidth]{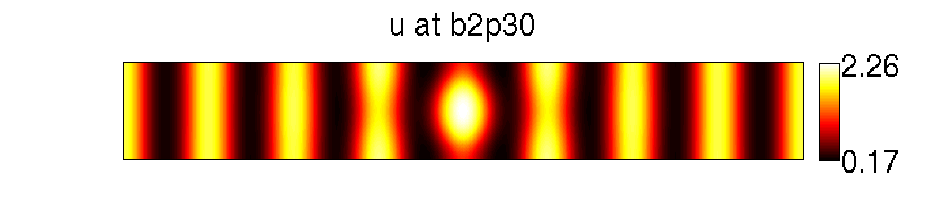}
\includegraphics[width=\textwidth]{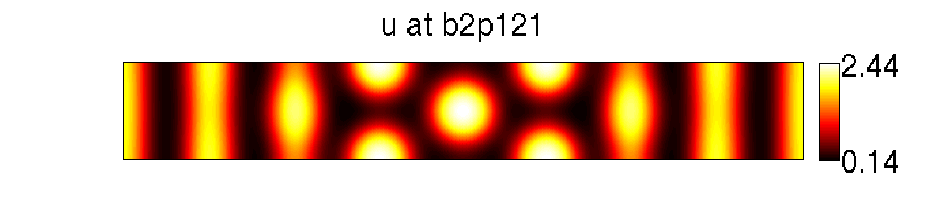}
\includegraphics[width=\textwidth]{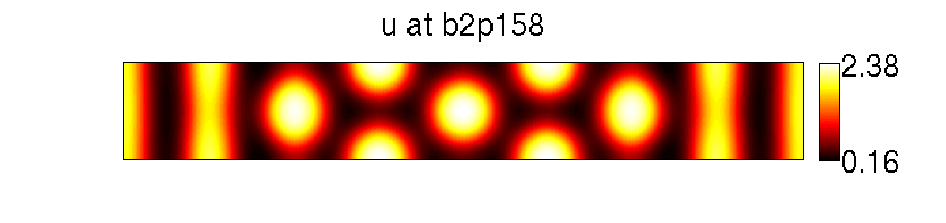}
\includegraphics[width=\textwidth]{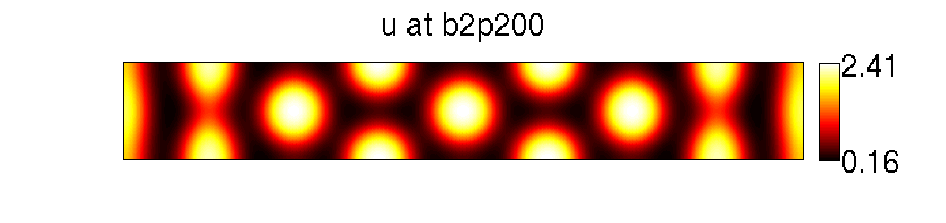}
\end{minipage} 
\caption{Density plots and bifurcation diagrams of solutions of \eqref{dgl2} for $\gamma=0.25$ over the domain $(-l_x,l_x) \times[-l_y,l_y]$, where $l_x=8\pi/kc$ and $l_y= 2\pi/(\sqrt{3} k_c)$. a) and b) {\tt b} is the hot-bean branch (violet). {\tt b1} and {\tt b2} are the branches which bifurcate from the first and second bifurcation points of the hot-bean branch.  c) and d) density plots of solutions labeled in a) and b), respectively.  }
\label{snakehs}
\end{figure}

Over the small domain, which is used for \figref{hsbbifu}, the hot-bean branch has 4 bifurcation points (not indicated). By increasing the horizontal direction by a factor of 4, the number of bifurcation points on the hot-bean branch also increases by a factor of 4 (see \figref{snakehs}). Periodic connections between stripes and hot hexagons branch from these bifurcation points. Here most things are similar to the 1D case described above. The bifurcation points are connected in the same way. Snaking branches of solutions bifurcate from the first and second ones, which look like front connections between hot stripes and hot hexagons and pulses of hexagons on homogeneous backgrounds, respectively. These states are actually periodic connections between stripes and hot hexagons, when we \replaced{extend}{ continue} this solutions periodically.  

\begin{figure}[H]

\begin{minipage}{0.49\textwidth}
\small{(a)}\\
\includegraphics[width=1\textwidth]{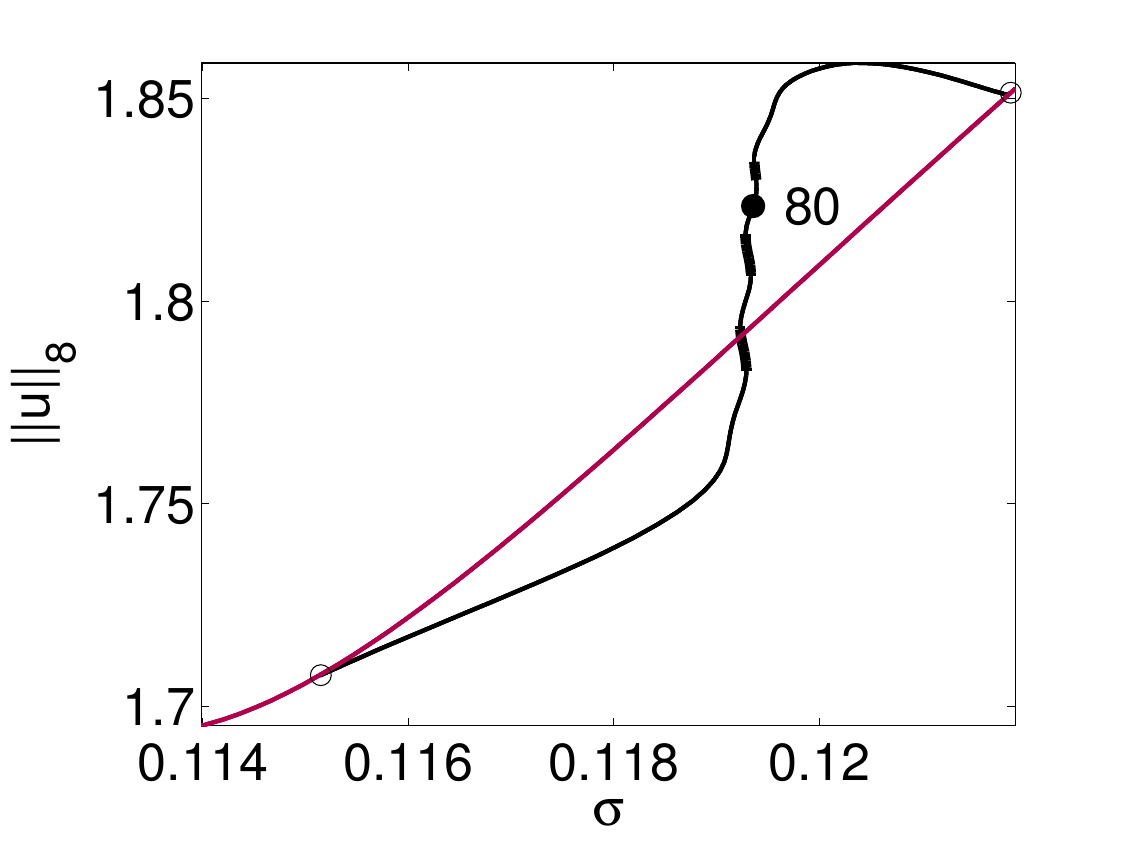}
\end{minipage}
\begin{minipage}{0.49\textwidth}
\small{(b)}\\
\includegraphics[width=1\textwidth]{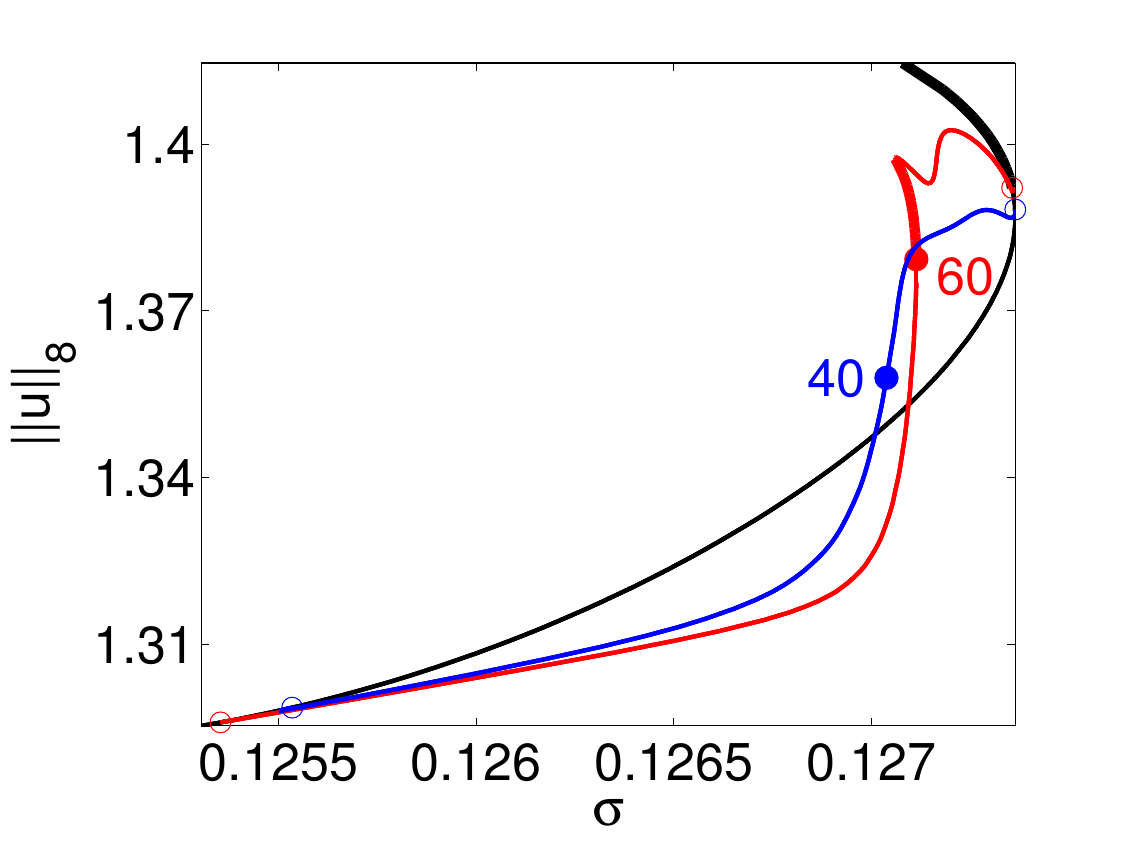}
\end{minipage}

\begin{minipage}{1\textwidth}
 \small{(c) 80th solution of the black branch in (a)}\\
\includegraphics[width=1\textwidth]{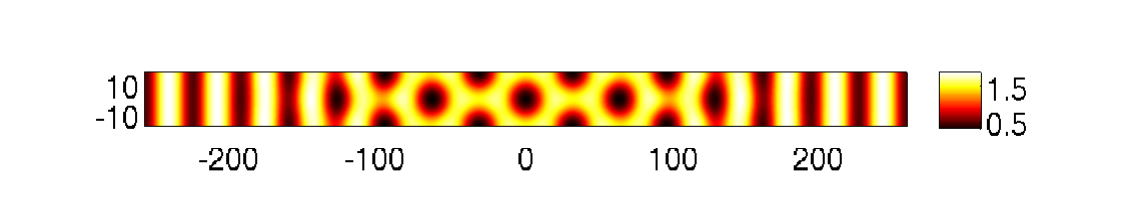}\\
\end{minipage}
\begin{minipage}{0.49\textwidth}
\small{(d) 60th solution of the red branch in (b)} \\
\includegraphics[width=1\textwidth]{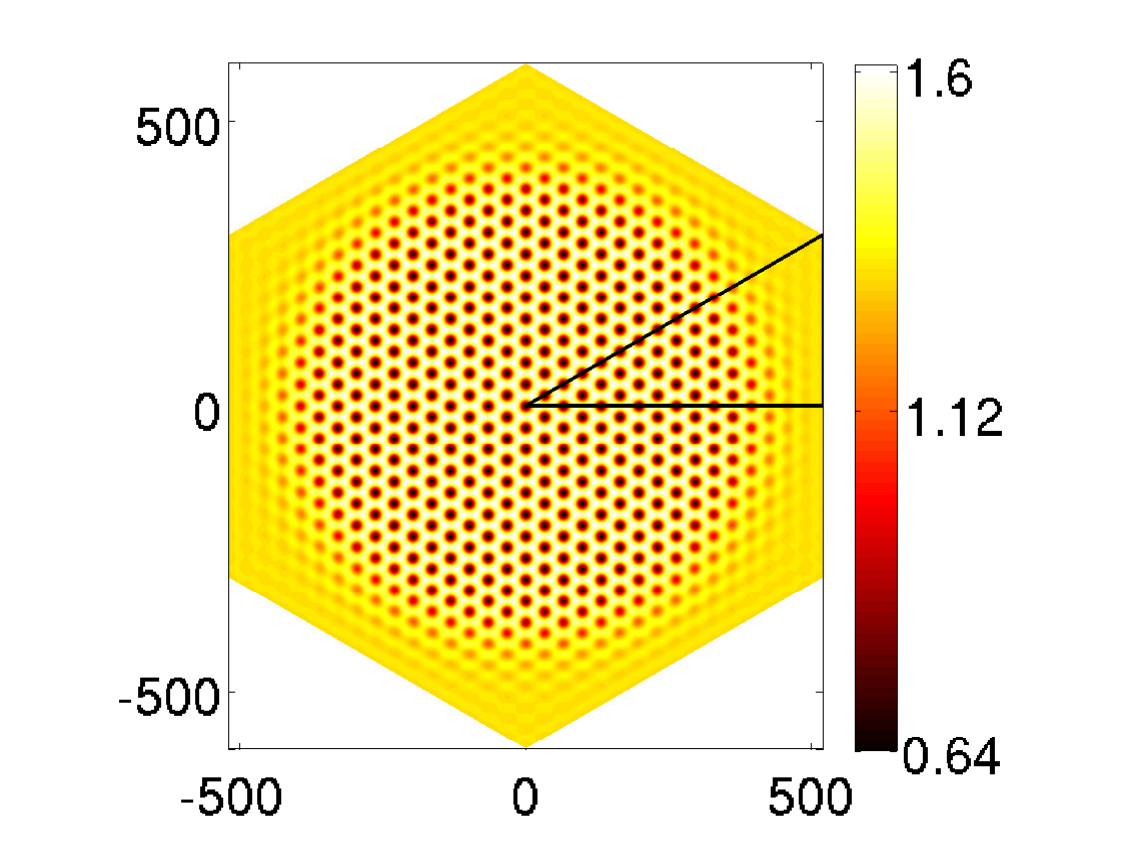}
\end{minipage}
\begin{minipage}{0.49\textwidth}
\small{(e) 40th solution of the blue branch in (b)}\\
\includegraphics[width=1\textwidth]{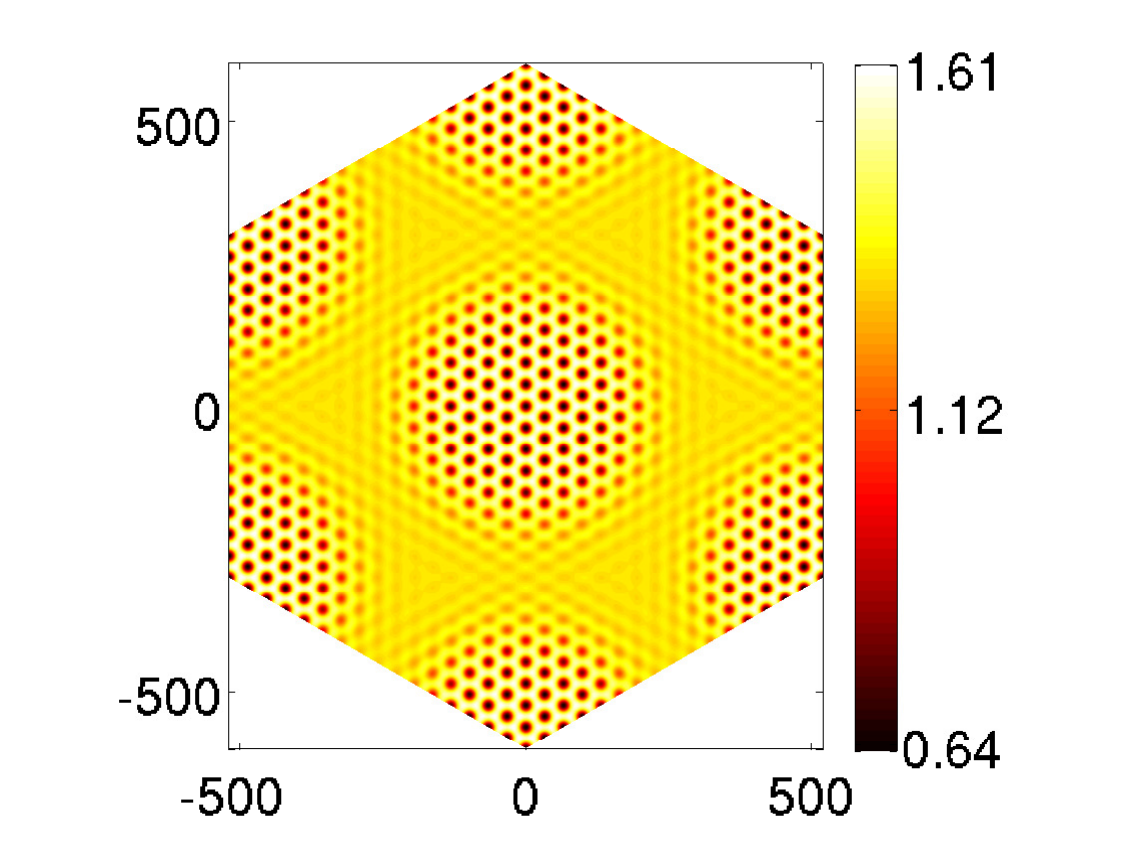}
\end{minipage}
\caption{
All plots are for $\gamma=0.25$. (a) Bifurcation diagram over the domain $(-l_x,l_x) \times(-l_y,l_y)$, where $l_x=16\pi/kc$,  $l_y= 2\pi/(\sqrt{3} k_c)$, and $k_c\approx 0.19$. The violet and black lines represent the branches of cold beans and localized cold hexagons on a striped background. (b) Bifurcation diagram over a triangular domain, which is given by the vertexes $(0,0)$, $(32 \pi /  k_c,0)$, and $(32 \pi /  k_c,32 \pi / (\sqrt{3} k_c))$. The black line represents the branch of cold hexagons, while the red and blue lines are the branches which bifurcate at the first and second bifurcation points after the fold of the cold hexagons, respectively. 
(c), (d), (e) Density plots of the solutions which are labeled in (a) and (b). For (d) and (e) we reflected and rotated the triangular domain to get a hexagonal domain. We call patterns, which are shown in (d) and (e), single and multi patches of hexagons, respectively. 
}
\label{coldbeansnake}
\end{figure}

In contrast to the 1D case the snake does not snake around a vertical line but in a slanted manner. Such localized patterns can also be found on cold bean branches (see \figrefa{coldbeansnake}). It can be seen in \eqref{hexamp} that the hexagons bifurcate in a disturbed pitchfork such that there is a bistable range between the homogeneous and hexagon solution.

We use a triangular domain to calculate branches of localized hexagons on homogeneous backgrounds for $\gamma=0.25$ (see \figref{coldbeansnake}). The solutions which bifurcate from the first and second bifurcation points of the cold hexagon branch are a single localized patch and multi localized patches of hexagons, respectively. Single patches are already observed and studied in \cite{hexsnake}. Multi patches are not mentioned in the literature before. We see that the patches themselves have a hexagonal structure. To understand the difference between single and multi patches one should have in mind that we can \replaced{extend}{ continue} the solutions periodically. If we do this for both, we see that the patches of the single patches lie edge to edge, while the multi patches lie corner to corner. Here only the single patch exhibits a snaking behavior with only one wiggle. 
By increasing the domain size, the number of wiggles should increase for both. 

Under the assumption that the Landau coefficients $c_2, \ c_3,\ c_4$ change much slower than $c_1$ by varying $\sigma$ away from $\sigma_c$, we can see in \eqref{hexamp} that the subcriticality increases if $c_f:=c_2^2/(4(c_3+2c_4)^2)$ increases. $c_f$ evaluated in $\sigma_c$ as function of $\gamma$ is shown in \figrefc{kandc3}. For $\gamma\approx 0.08$ it holds $c_3+2c_4=0$ such that we assume that we can increase the strength of subcriticality and with that the steepness of the connection and the width of snaking ranges by choosing a $\gamma$-value closer to 0.08. A similar prediction of the strength of subcriticality and width of the snaking branch for bean branches on a Landau level can be found in \cite{schnaki}.

\begin{figure}[H]
\centering{
\begin{minipage}{0.49\textwidth}
\small{(a) Bifurcation diagram for regular hexagons} \\
\includegraphics[width=1\textwidth]{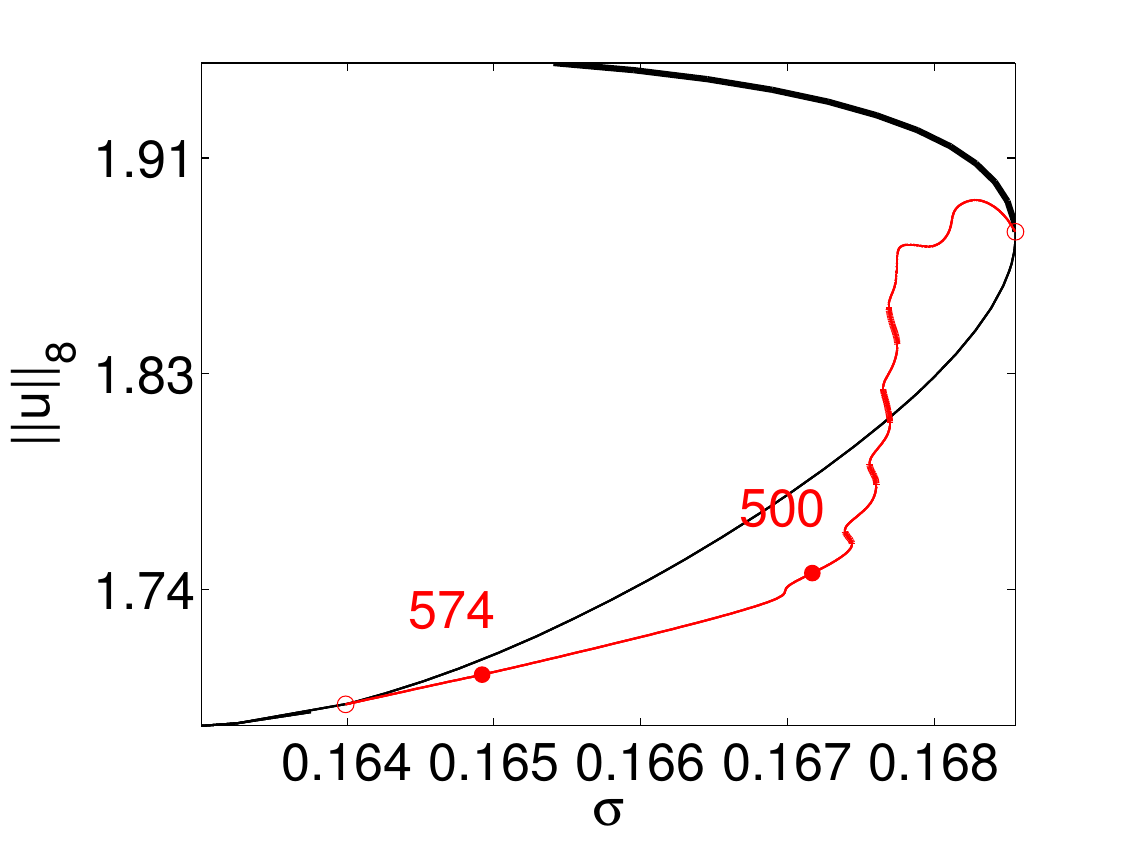}
\end{minipage}
\begin{minipage}{0.49\textwidth}
\small{(b) Bif. diagram for stretched hexagons} \\
\includegraphics[width=1\textwidth]{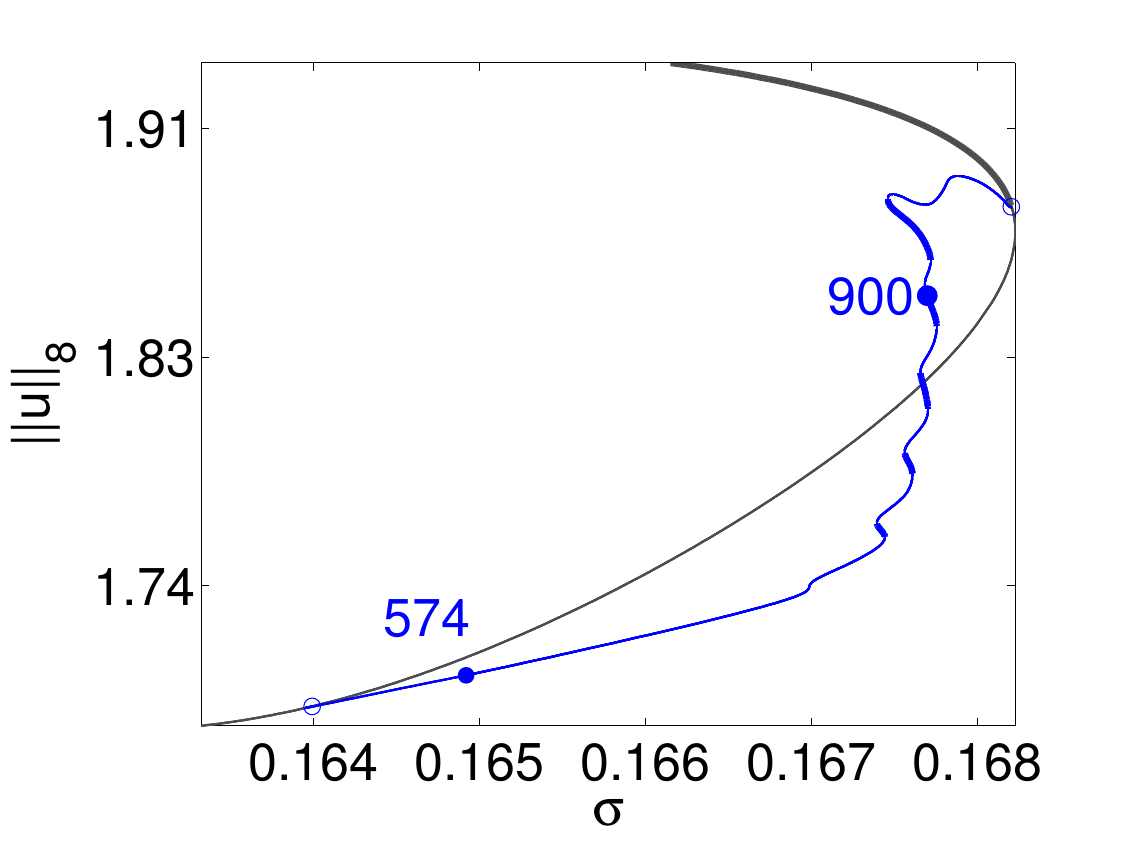}
\end{minipage}

\begin{minipage}{0.49\textwidth}
\small{(c) Fold of the black branch in (a)}\\ 
\centering{\includegraphics[width=0.66\textwidth]{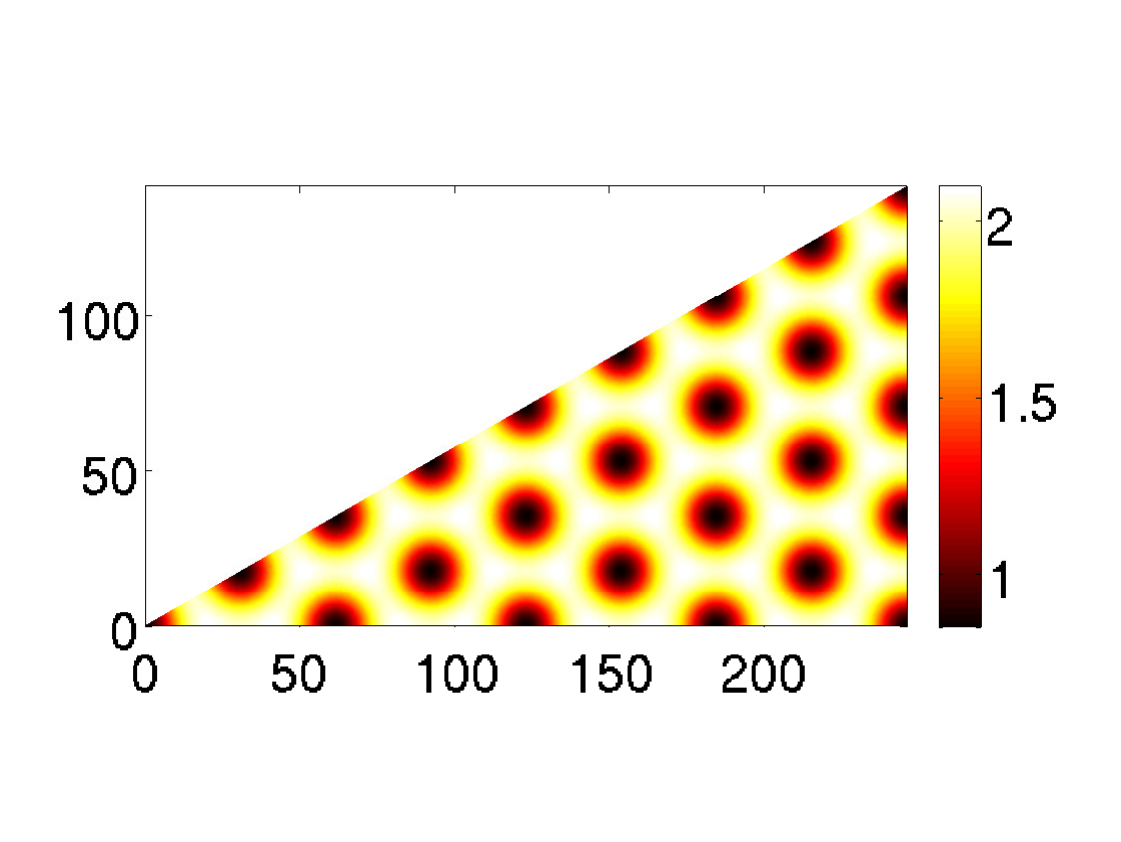}}
\end{minipage}
\begin{minipage}{0.49\textwidth}
\small{(d) Fold of the gray branch in (b)}\\
\centering{\includegraphics[width=0.66\textwidth]{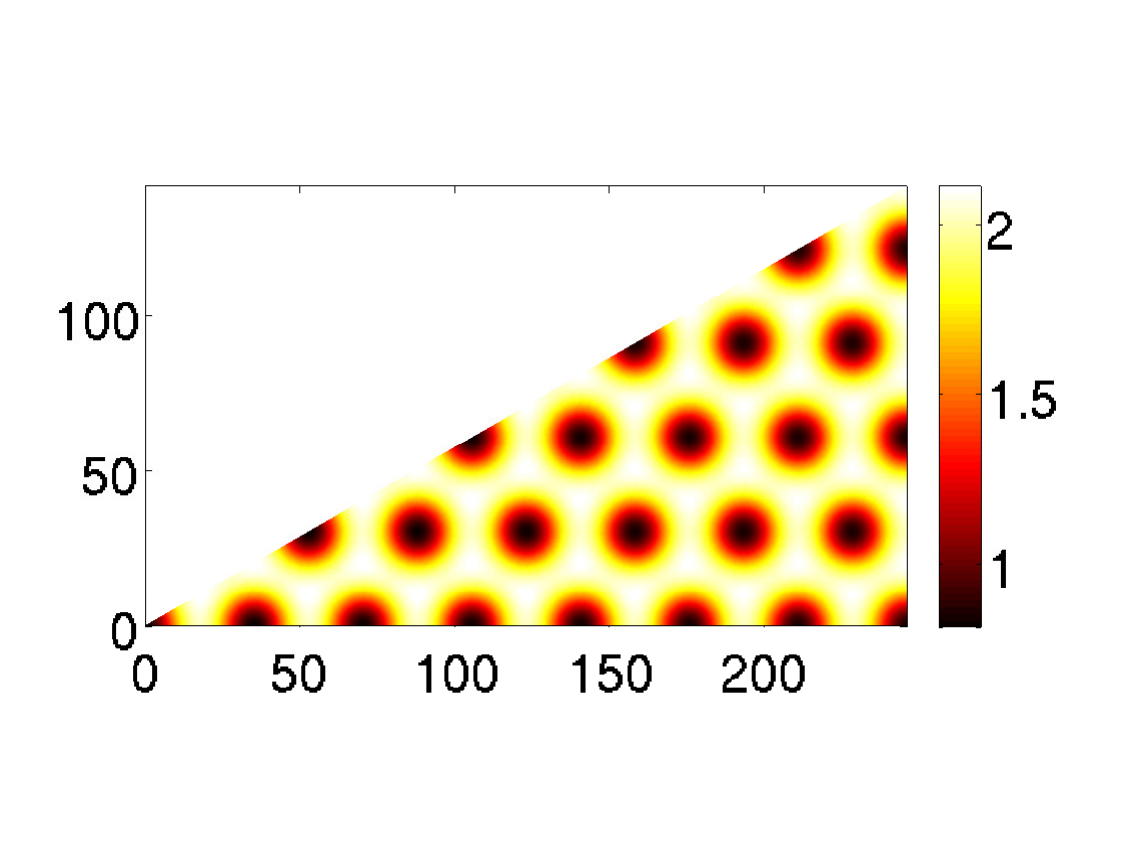}}
\end{minipage}
}

\begin{minipage}{0.32\textwidth}
\small{(e) 574 in (a) and (b)}\\
\includegraphics[width=1\textwidth]{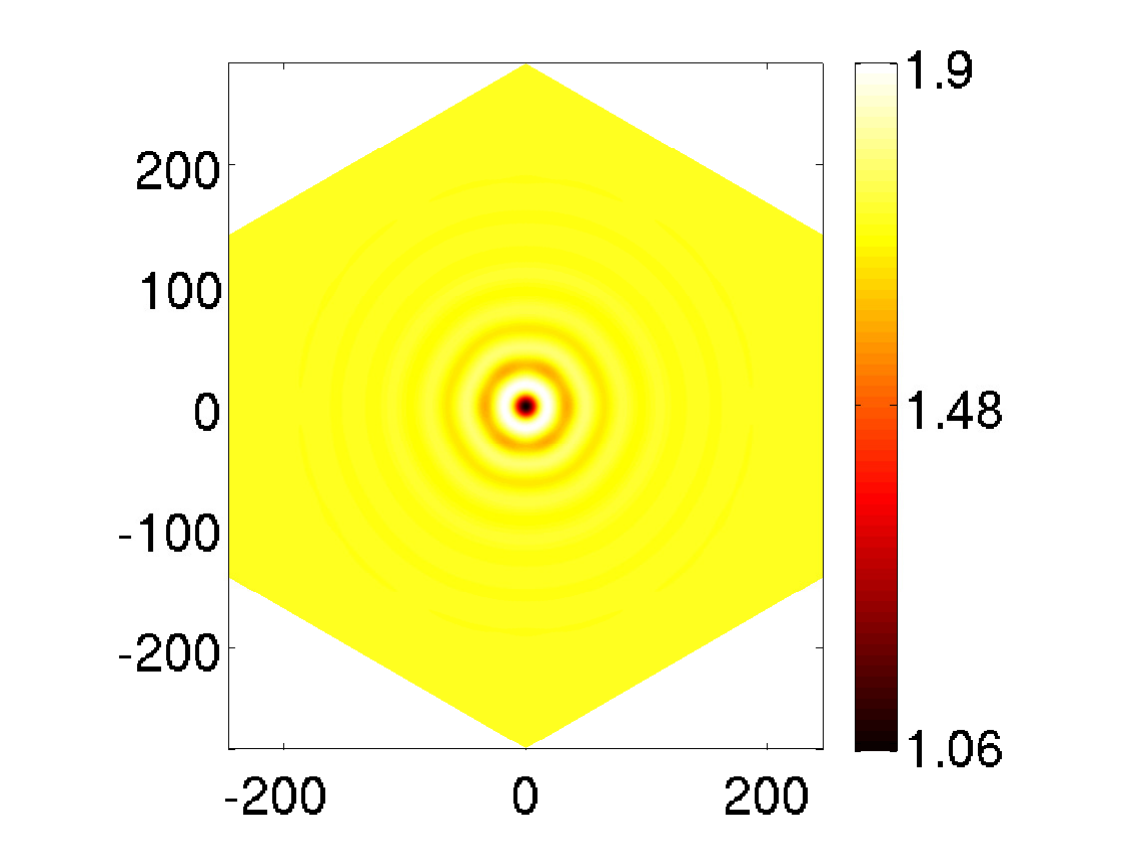}
\end{minipage}
\begin{minipage}{0.32\textwidth}
(f) 500 in (a)\\
\includegraphics[width=1\textwidth]{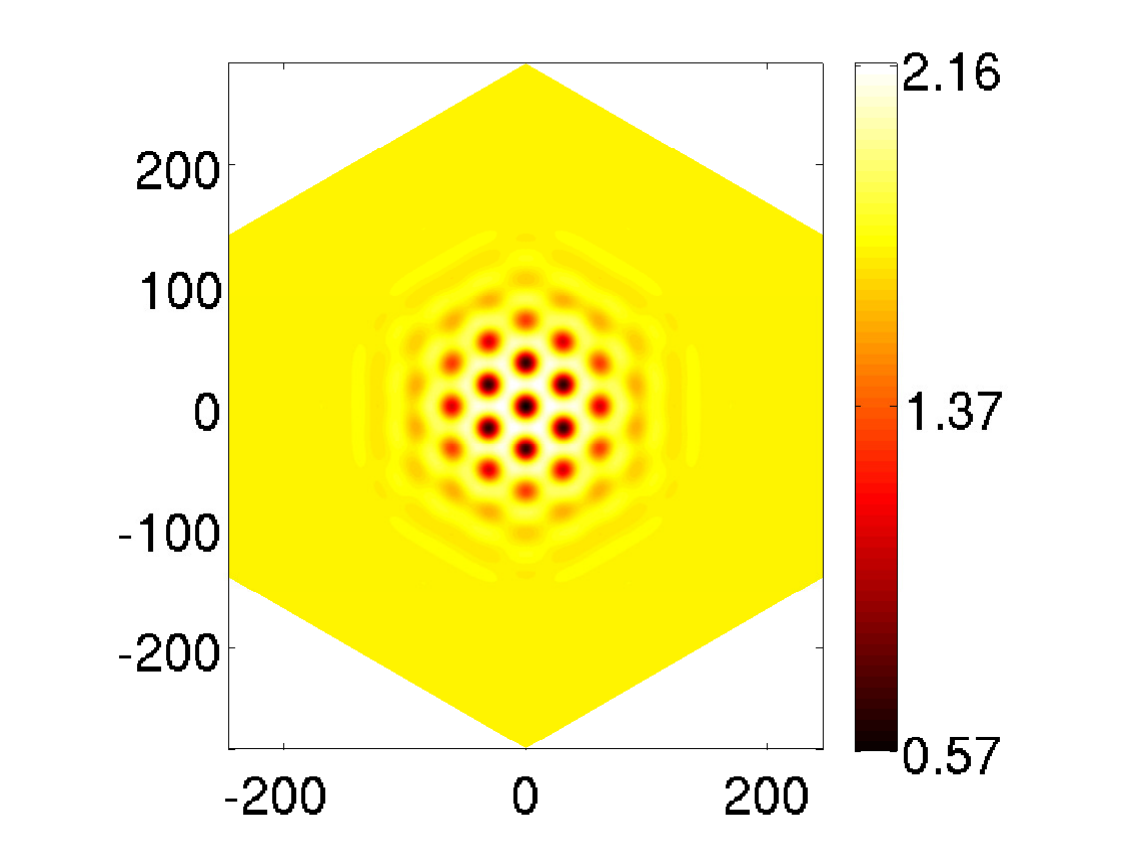}
\end{minipage}
\begin{minipage}{0.32\textwidth}
\small{(g) 900 in (b)}\\
\includegraphics[width=1\textwidth]{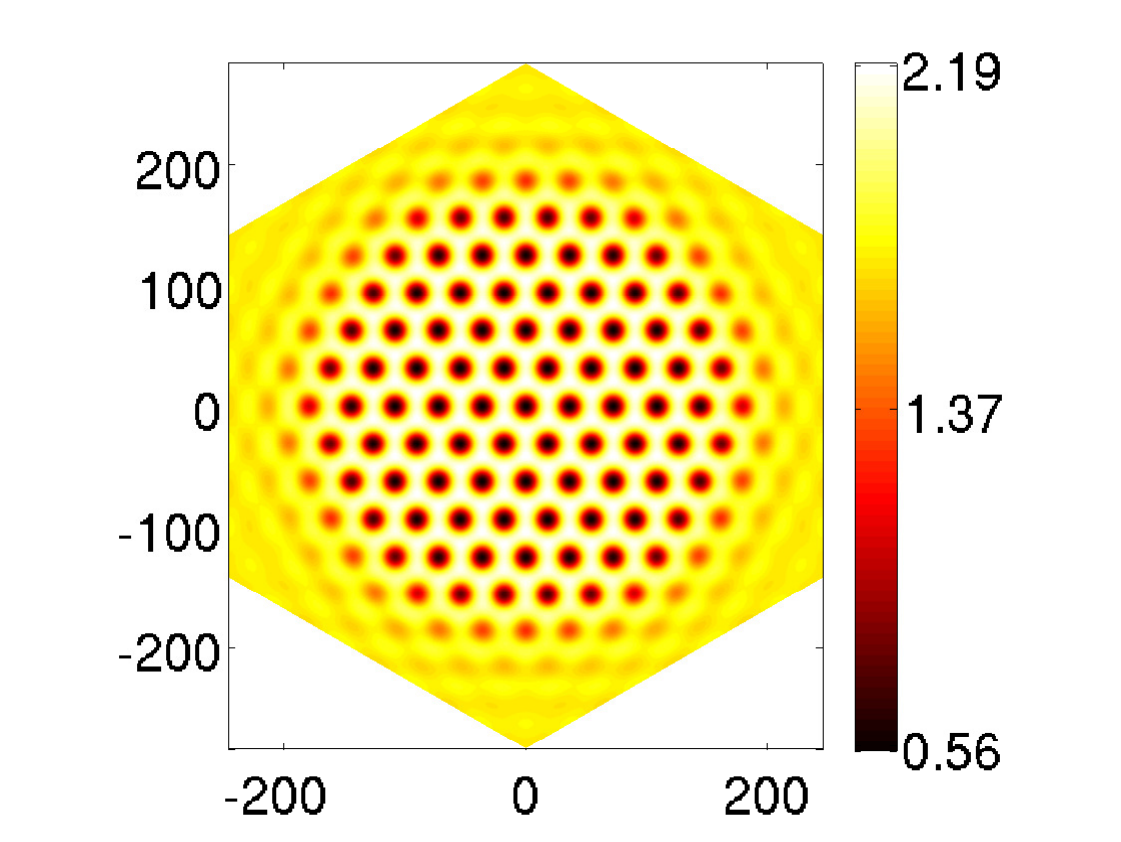}
\end{minipage}
\caption{All plots are for $\gamma=0.12$. (a) Bifurcation diagram over a triangular domain, which is given by the vertexes $(0,0)$, $(16 \pi /  k_c,0)$, and $(16 \pi /  k_c,16 \pi / (\sqrt{3} k_c))$. The black and gray lines in (a) and (b) represent the branches of regular and stretched cold hexagons, while the red and blue lines are branches which bifurcate from the first bifurcation point after the fold of the cold hexagons, respectively. (c), (d) Density plots of regular and stretched hexagons on the triangular domain.  (e), (f), (g) Density plots of solutions which are labeled in (a) and (b). Here we reflected and rotated the triangular domain to obtain a hexagonal domain.  
}
\label{coldhexstrong}
\end{figure}

In \figrefa{coldhexstrong} we see that the branch of single patches already shows a snaking behavior with more than one wiggle for $\gamma=0.12$ over a domain which is smaller than the domain which we used for $\gamma=0.25$. The \replaced{envelop}{overlying} function of the 574th solution becomes so steep that it looks like one single spot. From this point bifurcates a branch of a single hexagon patch which is rotated by $\pi/6$. Between the 574th and 900th solutions the branches in \figrefa{coldhexstrong} and (b) seem to be congruent. Beyond the 900th solution the boundary affects the rotated patches and the branch moves to a branch of stretched hexagons.

Investigations of localized stripes on homogeneous backgrounds over 2D domains can be found in \cite{strsnake}. The main result there is that branches of localized stripes do not snake if the stripes spread into the homogeneous stripe directions.

\section{Layering of patterns} \label{seclayering}
\begin{figure}[h]
\begin{center}
\includegraphics[scale=0.4]{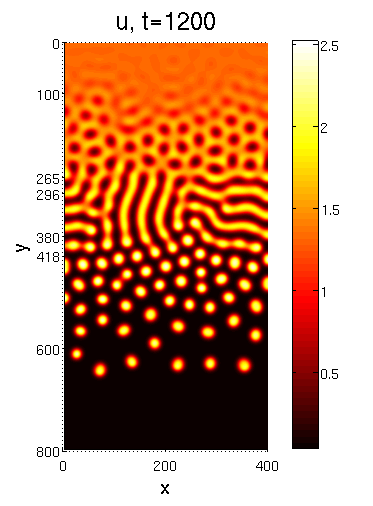}
\includegraphics[scale=0.4]{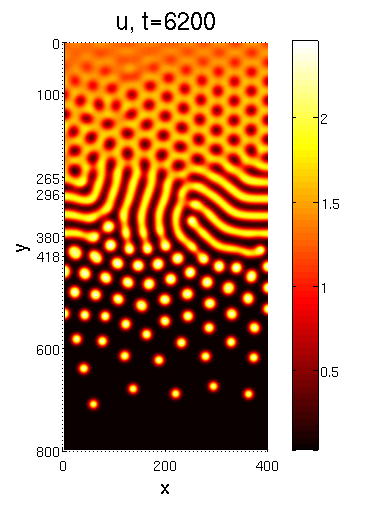}
\end{center}
\caption{Density plots of $u$ at times $t=1200$ and $6200$, which we found by using time-integration methods for a small random perturbation of $(u,v)\equiv(1,1)$. Here we use \eqref{dgl2} for $\gamma=0.25$ and parameterset \eqref{para}. The balancing rate $\sigma$ is described by \eqref{tanh}.}\label{orth}
\end{figure}

It is pointed out in \cite{ulrike} that it is realistic to consider a system for which the balancing rate decreases by increasing depth. This is done in \cite{ulrike} by using depth-dependent balancing-rate functions of the form
\alinon{\tilde{\sigma}(\tilde{y})=\alpha\e^{(-\tilde{y} / \mu)}}
for system \eqref{dgl}. 
Quasi-stationary layering of patterns involving stripes, spots and homogeneous states are observed (see Fig.11 and 12 of \cite{ulrike}). 
We assume that the exponential growth must decay so that a balancing rate function of the form  
 \ali{ \sigma (y)=\frac{0.128}{1+e^{0.011(y-480)}} \label{tanh} }
is more realistic. Using this function we obtain a quasi-stationary solution involving all five pattern types shown in \figref{feudelfig} (see \figref{orth}). Here we used time-iteration methods. \replaced{It is possible}{ One can expect} that there is a stable steady state which looks like the patterns shown in \figref{orth}, since the solution changes radically from $t=0$ to $t=1200$, while it does not change its general pattern from $t=1200$ to $t=6200$.
Using the Landau analysis described below, we predict a transition of stability between cold spots and stripes and between stripes and hot spots in the 
ranges $y\in[265,296]$ and $y\in[380,418]$, respectively. One can see in \figref{orth} that this is a sensible prediction.

\section{3D Patterns}\label{sec3d}
Investigations of 3D Turing patterns via amplitude equations and snaking branches of localized 3D patterns on homogeneous backgrounds can be found in \cite{knobloch3d97,knobloch3d99} and \cite{beaume13}, respectively. 
In \figref{hsbbifu} we show stripe and hexagon patterned solutions of \eqref{dgl}. We already mentioned above that these are the typical patterns, which are stable over 2D domains. Stripes and hexagons, which are \replaced{extended}{ continued} homogeneously into the third dimension, are \replaced{referred to as}{ called} lamellae and hexagonal prisms. We also call them stripes and hexagons if it is clear that we consider a 3D domain. They are also solutions over 3D domains, but can change their stability from stable to unstable. {\tt pde2path} uses MATLAB's PDE-Toolbox for the FEM, which only works for 2D domains, so one cannot consider PDEs over 3D domains with {\tt pde2path} alone. 
We use the continuation and bifurcation methods of {\tt pde2path} and the FEM of U. Pr\"ufert's PDE toolbox {\tt OOPDE} \cite{oopde} to study the stability of stripes and hexagons over bounded 3D domains.

We find that stripes and cold hexagons have the same stable ranges for the domains $\Omega_{3D}=(-l_x,l_x) \times(-l_y,l_y) \times (-l_z,l_z)$ and $\Omega=(-l_x,l_x) \times(-l_y,l_y)$. Here $l_x$ and $l_y$ are as defined in \figref{hsbbifu} and $l_z=200$. The stable range for hot hexagons starts in the same point ($\sigma\approx 0.1$) for $\Omega$ and $\Omega_{3D}$. The endpoint is different. The stability of hot hexagons over the 2D domain $\Omega$ ends at $\sigma\approx 0.04$, while it ends at $\sigma\approx 0.06$ for hot hexagons over $\Omega_{3D}$ with $l_z=10$. 
We tried to answer what kind of solutions bifurcate from this endpoint over 2D domains, but cannot give a precise result so far. By using a time integration method we found that hot hexagons, which correspond to a smaller wavenumber and which are stable on the used bounded domain, exist beyond this endpoint.

\begin{figure}[h]
\begin{minipage}{0.5\textwidth}
(a) Bifurcation diagram\\
\includegraphics[width=1\textwidth]{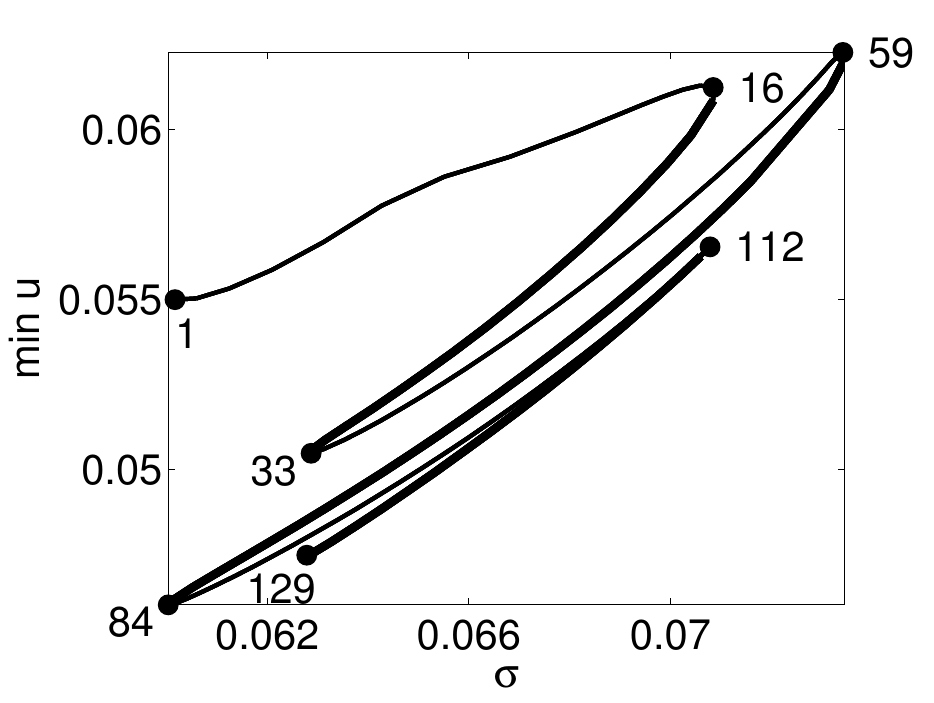}
\end{minipage}
\begin{minipage}{0.4\textwidth}
(b) 1st solution\\
\includegraphics[width=1\textwidth]{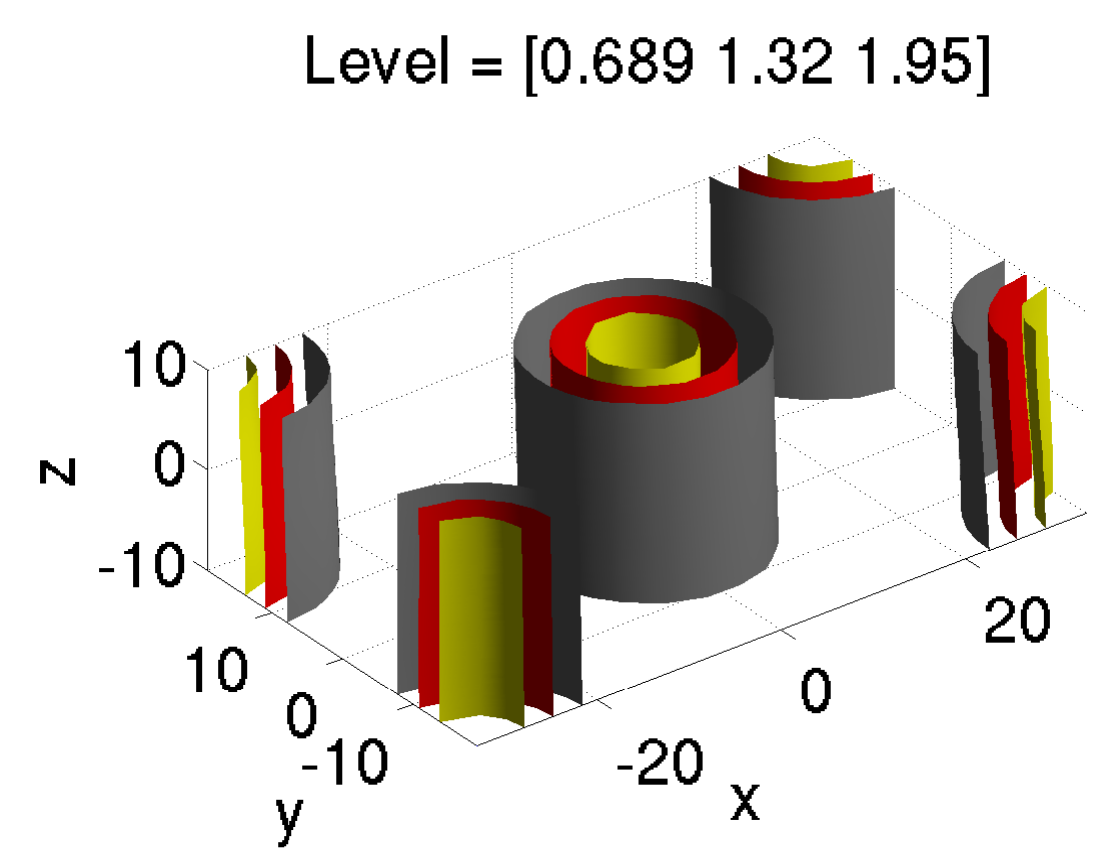}
\end{minipage}
\begin{minipage}{0.33\textwidth}
(c) 16th solution\\
\includegraphics[width=1\textwidth]{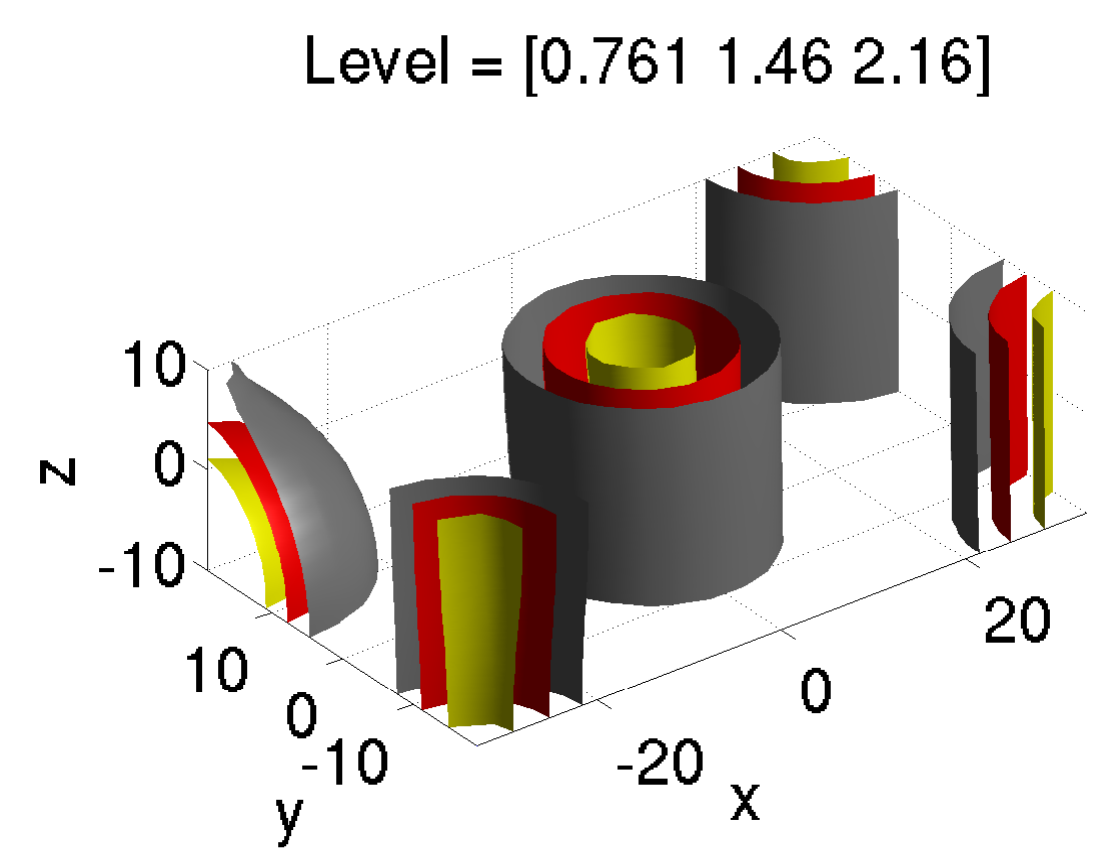}
\end{minipage}
\begin{minipage}{0.33\textwidth}
(d) 33th solution\\
\includegraphics[width=1\textwidth]{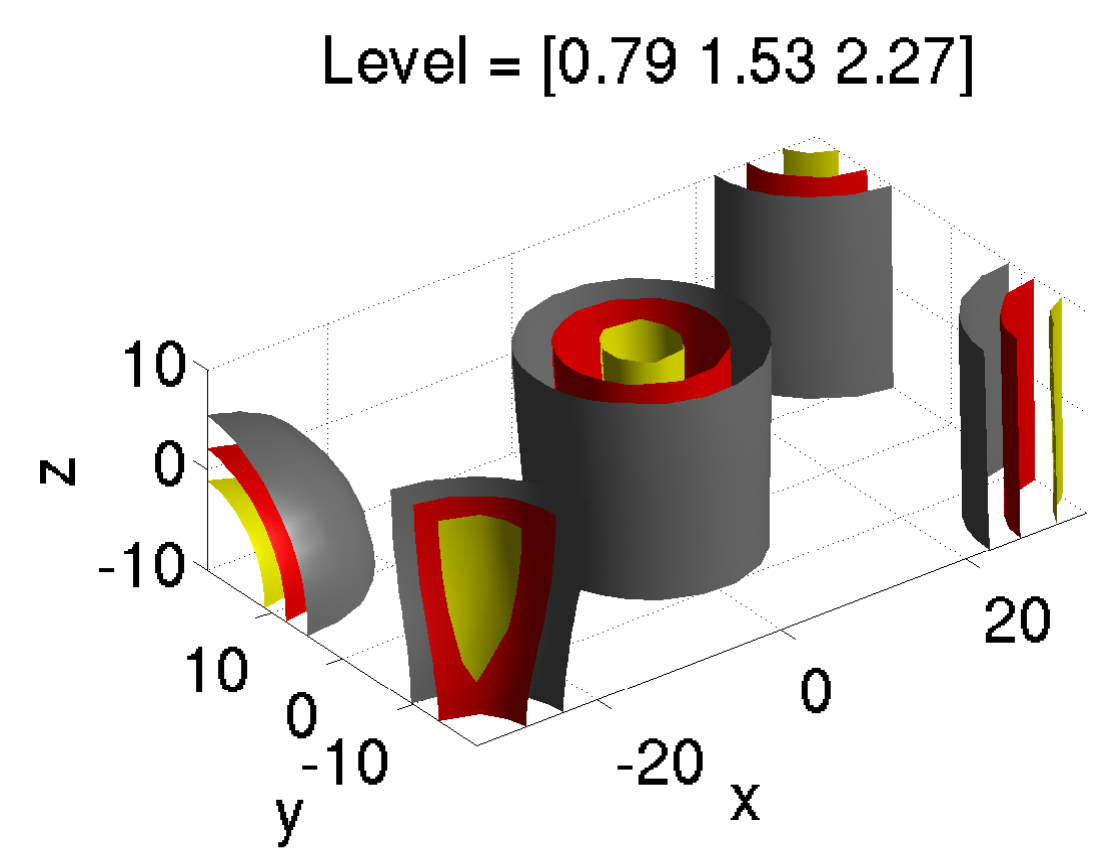}
\end{minipage}
\begin{minipage}{0.33\textwidth}
(e) 59th solution\\
\includegraphics[width=1\textwidth]{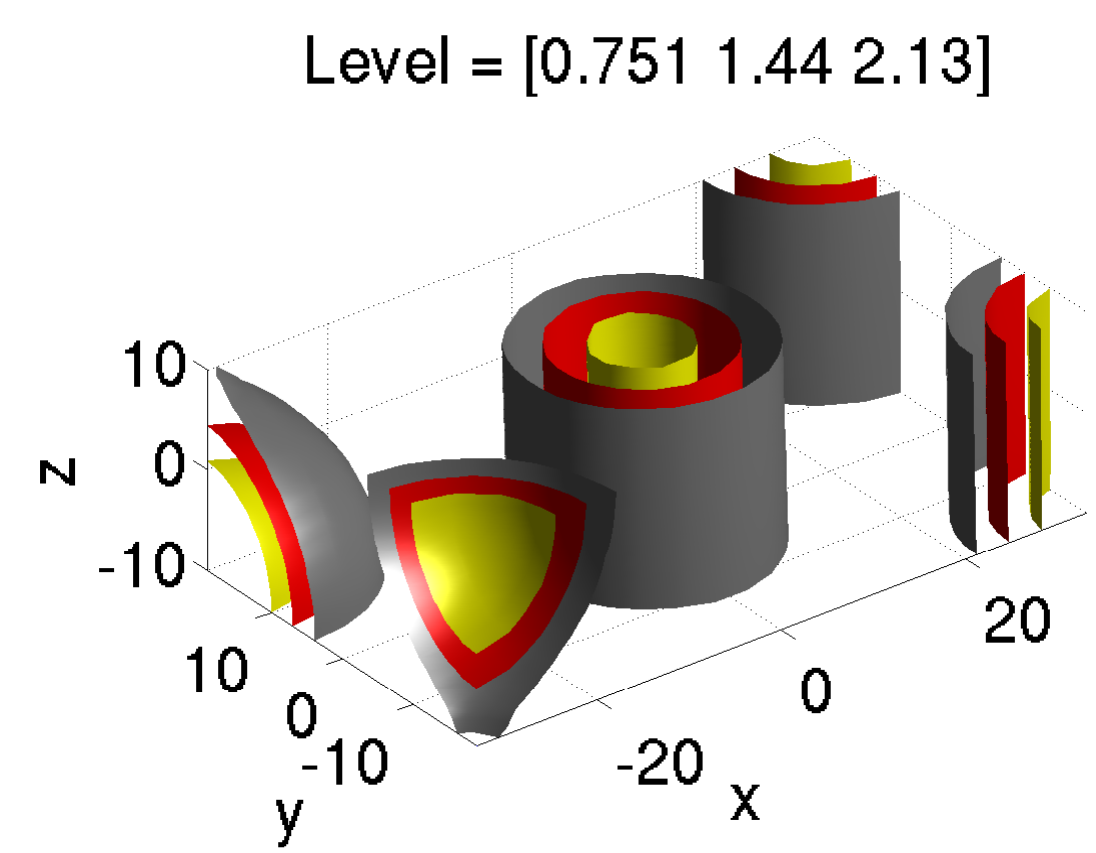}
\end{minipage}
\begin{minipage}{0.33\textwidth}
(f) 84th solution\\
\includegraphics[width=1\textwidth]{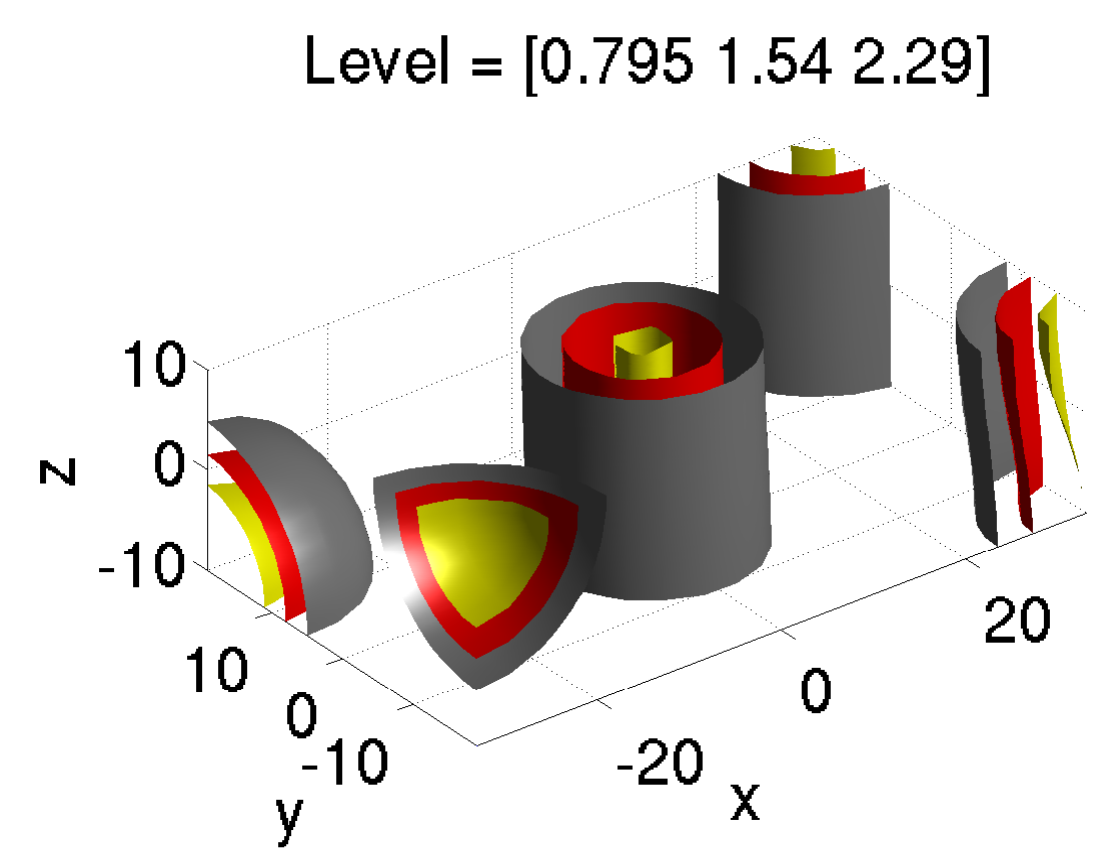}
\end{minipage}
\begin{minipage}{0.33\textwidth}
(g) 112nd solution\\
\includegraphics[width=1\textwidth]{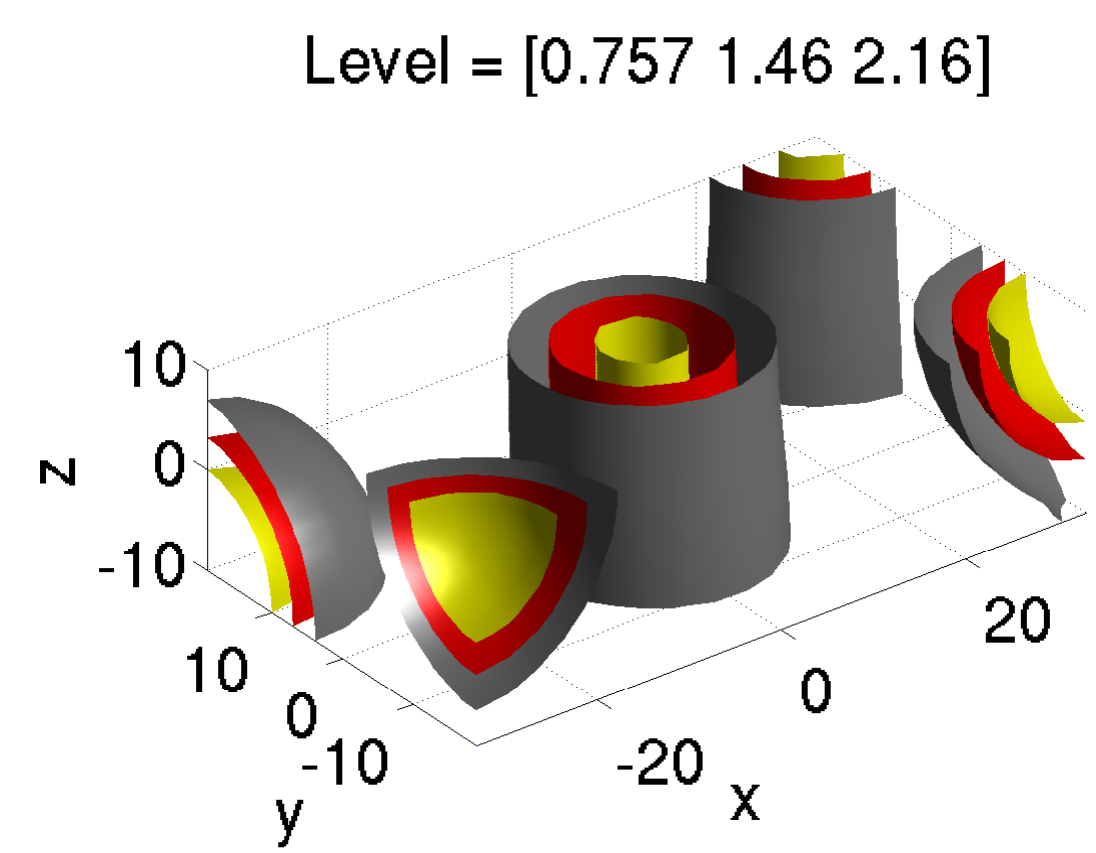}
\end{minipage}
\begin{minipage}{0.33\textwidth}
(h) 129th solution\\
\includegraphics[width=1\textwidth]{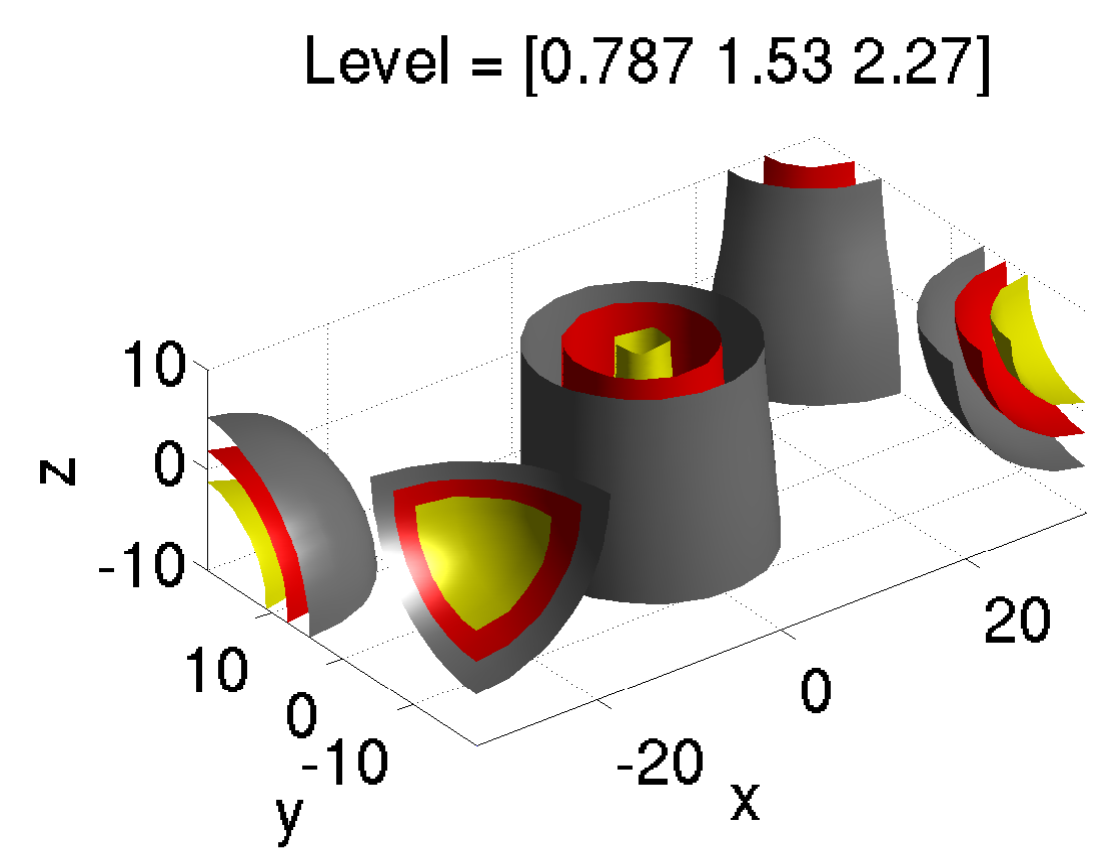}
\end{minipage}
\caption{(a) Shown is a branch, which bifurcates, when the hot hexagons loose their stability over the domain $\Omega_{3D}=(-l_x,l_x) \times(-l_y,l_y) \times (-l_z,l_z)$ with $l_z=10$. (b)-(h) Shown are isoface plots of solutions, which are labeled in (a). The first, second, and third entries of the level vector represent dark gray, red, and yellow faces, respectively. }\label{3d}
\end{figure}

However, we found a branch of solutions, which seems to bifurcate from the left boundary of the stable range for hot hexagons over $\Omega_{3D}$ with $l_z=10$. This branch and some solutions are shown in \figref{3d}. One can see that the branch has a snaking behavior. The 1st solution is the bifurcation point, which also lies on the hexagon branch. Following the snaking branch, one can see that the outer hexagonal prisms deform one by one to genuine 3D patterns, which look like balls. One can see in \figref{3d}(h) that three outer hexagonal prisms are already deformed to balls for the 129th solution. Beyond the 129th solution the branch turns around and its solutions become unstable (not shown). A bit later there is fold and the last outer hexagonal prism is deformed completely to a ball. Beyond this fold the branch does not become stable and the balls deform back to prisms one by one. The corresponding branch has a snaking behavior without stable ranges. 

Here we guess that one can have better results by using another domain. However, this is only an outlook for 3D Turing patterns. It shows that most of the stable 2D patterns are also stable in 3D, that genuine 3D patterns play an important role, and that periodic connections between 2D patterns and 3D patterns exist. 
\section{Discussion} \label{discussion}
In this paper we investigate\added{d} the bacteria-nutrient system \eqref{dgl2} with respect to $\sigma$ and $\gamma$, which are the balancing rate of the nutrient and activity strength of bacteria and which regulate the strength of food supply and ingestion, respectively. To understand the long time behavior of this system, we \replaced{studied}{ study} time independent solutions, which are also called equilibria, steady states, or stationary solutions. 

In \secref{sechom} we \replaced{studied}{ study} such solutions with a spatial homogeneous distribution of bacteria and nutrients and also the stability of those stationary homogeneous states. 
We \replaced{saw}{ see} that there is \replaced{only}{ still} one homogeneous steady state for \replaced{large}{ great} values of $\sigma$ and $\gamma$ and two stable homogeneous steady states if $\gamma$ decreases. One has a high bacteria population density and low nutrient concentration, while this is the opposite case for the other one. This is \replaced{a}{ very} useful \replaced{information}{ to know}, because if the bacteria population density is near the bacteria-poor homogeneous steady state, it can be possible to reach the bacteria-rich homogeneous steady state by putting bacteria into the sediment. Analogously it can be expected that a change from the bacteria-rich homogeneous steady state to the bacteria-poor one can be achieved by taking some bacteria away. These ideas do not hold for great values of $\gamma$ and $\sigma$. If we put bacteria into the system, we cannot expect that the system stays in this state for a long time, but falls back to the stable homogeneous steady state or to another stable steady state, which is not homogeneous.    

\replaced{Our}{ Furthermore, our} analyses \replaced{has shown}{ shows} that there is only one stable positive stationary homogeneous state for all positive $\sigma$ and \replaced{large}{ great} values of $\gamma$. Thus we can conclude that very active bacteria populations track the homogeneous steady state for a slowly decreasing balancing rate $\sigma$ and furthermore that in this case the bacteria population density decreases, while the nutrient concentration increases.
For smaller values of $\gamma$ we cannot predict this behavior from our analysis of homogeneous steady states, because the solution branch of the homogeneous steady state has Turing-unstable ranges. 

It is known that inhomogeneous stationary solutions bifurcate from Turing-unstable states. We use\added{d} the Landau formalism (see \secref{seclandau}) to describe the structure via formulas and to understand the bifurcation behavior of these so-called Turing patterns. The Landau formalism works locally and fails globally \added{in parameter space}. Thus we \replaced{used}{use} numerical path following methods to generate global bifurcation diagrams. 

As an introduction \replaced{to}{ for} Turing patterns we \replaced{studied}{ study} in \secref{secchanging} the system \eqref{dgl2} for $\gamma=0.3$ over the bounded 1D domain $\Omega^r$. Here we \replaced{saw}{ can see} in the bifurcation diagram of the homogeneous state that there is still one homogeneous steady state for all $\sigma$ (see \figref{bifuhomplot} d)), which is stable for small and large values of $\sigma$ and Turing-unstable in the middle. The domain $\Omega^r$ has a length such that a 1D Turing pattern with 4 periods can bifurcate directly, when the homogeneous steady state changes its stability from stable to Turing-unstable by decreasing $\sigma$. This pattern is stable directly after its bifurcation. Decreasing $\sigma$, it becomes unstable, while Turing patterns with 3.5 and 3 periods continue stably. Increasing $\sigma$ a bit more, the Turing pattern with 3.5 periods becomes unstable before the pattern with 3 periods will do the same. 
Thus there is a selection of \replaced{larger}{ greater} wavelengths, when the food supply decreases. Such a behavior is already observed in \cite{jens14} for an extended version of the Klausmeier model.

Let us assume that we are observing a real natural bacteria-nutrient system, which can be modeled by system \eqref{dgl2} for the parameter set \eqref{para} and $\gamma=0.3$ on a \added{quasi 1D} domain for which two spatial directions are very small and the third direction has the same length like $\Omega^r$\deleted{, i.e, it is quasi 1D}. If we observe for a long time that the bacteria population density does not change and is also spatially homogeneous, \replaced{then}{ than} we can assume that the system stays in a stable homogeneous steady state. If we are able to measure the bacteria population density, then we know also, \replaced{whether}{ if} we are on the right or left side of the Turing-unstable range. Let us assume that we are on the right side and after a certain time the spatial distribution of the bacteria population changes to a patterned state with 4 periods, then we can assume that the balancing rate $\sigma$ has decreased and the bacteria population is in danger of dying out. 

If we are not able to measure the bacteria population density and we observe that the state changes from a homogeneous state to periodic states and if we furthermore see that the number of periods decreases, then we can conclude that $\sigma$ has decreased and that the bacteria population is in danger. If the system is \replaced{situated}{ provided} in a homogeneous state on the left side beyond the Turing-unstable range and we try to increase $\sigma$ to reactivate the bacteria population, then we know that we are doing the right thing, when we see that the bacteria distribution becomes periodic and that the number of periods increases after a certain time.

In \secref{2d} we \replaced{studied}{ study} the system \eqref{dgl2} for $\gamma=0.25$ over a small 2D domain. Our analysis of the homogeneous states \added{has} already \replaced{shown}{ shows} that we have two stable homogeneous steady states for \replaced{large}{ great} values of $\sigma$. Decreasing $\sigma$ the bacteria-rich homogeneous steady state becomes Turing-unstable, and Turing-patterns bifurcate, which we call stripes, hot and cold hexagons (see \figref{hsbbifu}).

Let us assume again that we observe a real bacteria-nutrient scenario on a quasi 2D domain, which has the same size as the one we use in \figref{hsbbifu} and let us see what we can learn from our numerical investigations. A change from a homogeneous steady state to cold hexagons indicates the possibility that the balancing rate $\sigma$ has decreased. Furthermore, we see that if the decrease of $\sigma$ continues, the next stable steady states are stripe patterns followed by hot hexagons. For all three patterns of the bacteria population density we have a signal that $\sigma$ decreases, when the hot part of the pattern shrinks. For instance when the size of hexagonal hot spots becomes smaller.  

One main point of this paper is the investigation of bistable ranges between two steady stats for which at least one of the states is a Turing pattern. In all these bistable ranges we \replaced{found}{ find} solutions numerically, which change their patterns along the horizontal spatial direction. The branches of these solutions snake back and forth and have stable and unstable ranges. It is important to know that such stable mixed patterns exist, because if our system is in a bistable range and we want to bring it into another state by using any methods, it can happen that the solution gets caught on a mixed pattern.

In \secref{seclayering} we \replaced{considered}{ consider} layering of patterns, which occur by space dependent parameters. We also learn from this paper that a layering of two different patterns for any reaction-diffusion system is not necessarily an effect of space dependent parameters, but may be due to the system \deleted{is }being in a bistable range.
\section*{Acknowledgments} I am indebted to H. Uecker, J. Rademacher, S. U. Gerbersdorf, P. Harmand, U. Pr\"ufert, M. Herrmann, H. Susanto and U. Feudel for stimulating and enlightening discussions. This work is founded by the DFG, the University of Oldenburg, and the University of Bremen.

\bibliography{snbib}{}
\bibliographystyle{plain}
\medskip
Received xxxx 20xx; revised xxxx 20xx.
\medskip
\end{document}